\documentclass[review,1p,twocolumn]{elsarticle}

\usepackage{hyperref}

\journal{Computers \& Structures} 
 
\usepackage{graphicx}
\usepackage{amsmath}
\usepackage{amssymb}
\usepackage{bm}
\usepackage{placeins} 									
\usepackage{tabularx}
\newcolumntype{Y}{>{\centering\arraybackslash}X}
\newcolumntype{U}{>{\hsize=.080\hsize \centering\arraybackslash}X}
\newcolumntype{X}{>{\hsize=.220\hsize \centering\arraybackslash}X}
\newcolumntype{V}{>{\hsize=.120\hsize \centering\arraybackslash}X}
\newcolumntype{A}{>{\hsize=.100\hsize \raggedright\arraybackslash}X} 
\newcolumntype{B}{>{\hsize=.060\hsize \raggedright\arraybackslash}X} 
\newcolumntype{C}{>{\hsize=.070\hsize \raggedright\arraybackslash}X} 
\newcolumntype{D}{>{\hsize=.060\hsize}X}  
\newcolumntype{E}{>{\hsize=.075\hsize}X}  
\newcolumntype{F}{>{\hsize=.159\hsize}X} 
\newcolumntype{G}{>{\hsize=.068\hsize}X} 
\newcolumntype{H}{>{\hsize=.380\hsize}X} 
\newcolumntype{s}{>{\hsize=.5\hsize \centering}X}
 
\def\tabulkywidth{125mm} 
\def\tabulkywidthB{217mm}
	
\usepackage{caption} 
\usepackage[percent]{overpic}
\usepackage{arydshln}
\usepackage{verbatim}
\usepackage{float}
\usepackage{multirow}  
\usepackage[bottom]{footmisc} 
\usepackage{pgfplots} 
\pgfplotsset{compat=1.5}
\usepackage{multicol} 
\usepackage{adjustbox} 

\def\PART#1#2{\frac{\partial #1}{\partial #2}}

\def\m{\text{m}}
\def\Nd{N_{\text{d}}}
\def\T{\textrm{T}} 
\def\t{\text{t}}
\def\c{\text{c}}
\def\vv{\text{v}}

\def\p{\text{p}}
\def\ep{\dot{\bm{\varepsilon}}}
\def\eep{\dot{\bm{\mathcal{E}}}}

\def\kp{\dot{\bm{\kappa}}}
\def\epp{\dot{\bm{\epsilon}}}
\def\dOmega{\,\text{d}\Omega}

\def\intz{\int_{-\frac{h}{2}}^{\frac{h}{2}}}

\def\dz{\,\text{d}z}

\definecolor{col1}{RGB}{0,125,0}

\makeatletter
\newcommand*{\Biggg}[1]{{\hbox{$\left#1\vbox to28\p@{}\right.\n@space$}}}
\newcommand*{\Bigggg}[1]{{\hbox{$\left#1\vbox to42\p@{}\right.\n@space$}}}
\makeatother

\begin{document}

\begin{frontmatter}

\title{The performance of the DK (Discrete Kirchhoff) and DKM (Discrete Kirchhoff-Mindlin) shell element classes in the upper bound finite-element-based limit analysis and in the sequential finite-element-based limit analysis}

\author[TUaddress]{V\'{\i}t\v{e}zslav \v{S}tembera$\!$\corref{cor1}}
\ead{vitastembera@hotmail.com}
\address[TUaddress]{Institute for Mechanics of Materials and Structures, TU Wien, Karlsplatz 13, 1040 Vienna, Austria} 
\cortext[cor1]{Corresponding author}
 
\begin{abstract}	
	This paper investigates the applicability of the DK and DKM shell element classes for the first time within the framework of the standard (and sequential) finite-element-based limit analysis, which is a direct method used for determining the plastic collapse (and post-collapse) behaviour of structures. Despite these elements not being able to guarantee a strict upper bound, it is shown that they exhibit a significantly lower discretization error compared to the strict upper-bound formulated shell elements used in literature. The higher numerical efficiency of quadrangle elements over triangle elements is observed. The numerical results are also compared with the MITC shell elements, which have been shown to converge only in the pure membrane case.	
\end{abstract}



\end{frontmatter}

\section*{Nomenclature}

\subsubsection*{Abbreviations}
\begin{tabular}{ll}
	BSQ           &Bilinear Strain Quadrangle \\    
	c.s.          &coordinate system \\
	CST           &Constant Strain Triangle \\  
	DKMQ          &Discrete Kirchhoff-Mindlin Quadrangle \\
	DKMT          &Discrete Kirchhoff-Mindlin Triangle \\ 
	DKQ           &Discrete Kirchhoff Quadrangle \\
	DKT           &Discrete Kirchhoff Triangle \\
	dofs          &degrees of freedom \\
	FELA          &Finite-Element-based Limit Analysis \\
	OPQ           &OPtimal Quadrangle \\
	OPT           &OPtimal Triangle \\    
	SFELA         &Sequential Finite-Element-based Limit Analysis \\
	TRIC          &TRIangular Composite \\
\end{tabular} 

\subsubsection*{Latin symbols}
\begin{tabular}{ll}
	$d_{\p}^{z}$                                  & dissipation power density\\
	$d_{\p}$                                      & dissipation power area-density\\
	$D_{\p}$                                      & dissipation power \\
	$E$                                           & Young's modulus \\
	$f$                                           & failure criterion\\ 
	$\sigma_0$                                    & yield strength in the isotropic case\\ 		 
	$f_{\text{t}x},~f_{\text{t}y}$                & yield strengths in tension\\ 		
	$f_{\text{c}x},~f_{\text{c}y}$                & yield strengths in compression\\ 	
	$f_{\text{v}xy}$                              & yield strength in shear in $xy$ plane\\
	$\bm{\mathfrak{f}}$                           & discretized loading vector\\
	$g$                                           & gravitational acceleration \\
	$G$                                           & shear modulus \\	
	$h$                                           & shell thickness\\
	$m_{\p}=\sigma_0\frac{t^2}{4}$                & maximum plastic moment\\
	$n$                                           & number of element nodes\\
	$n_{\p}=\sigma_0t$                            & maximum plastic in-plane force\\
	$N$                                           & number of optimization unknowns\\
	$\Nd$                                         & number of dofs (including zero dofs) \\
	$N_e$                                         & number of elements\\
	$N_q$                                         & number of quadrature points in an element\\
	$N_r$                                         & number of quadrature points in $z$-direction\\
	$p$                                           & pressure\\
	$\dot{\bm{q}} = \left[\dot{\bm{\varphi}},~
		\dot{\bm{u}}
	\right]^{\T}$                                 & vector of velocities at a given point\\
	$\dot{\bm{\mathfrak{q}}}$                     & vector of all unknown velocities  \\ 
	$\dot{u}_x,\dot{u}_y,\dot{u}_z$               & nodal velocities in the local c.s.  \\
	$\dot{u}_X,\dot{u}_Y,\dot{u}_Z$               & nodal velocities in the global c.s.  \\
	$w_i$                                         & $i-$th quadrature weight  \\
	$xyz$                                         & coordinates of the local c.s.  \\
	$XYZ$                                         & coordinates of the global c.s.  \\
\end{tabular} 

\subsubsection*{Greek symbols}
\begin{tabular}{ll}
	$\Gamma$										  &  $\Omega$ boundary  	 \\
	$\Gamma_i$										  &  part of the $\Omega$ boundary at which $q_i=0$ \\	
	$\dot{\epsilon}_x, \dot{\epsilon}_y, \dot{\curlyvee}_{xy}$ &  in-plane plastic strain rates \\ 
    $\epsilon_1,\epsilon_2$                     &  in-plane eigenstrain rates \\
	$\dot{\varepsilon}_x, \dot{\varepsilon}_y, \dot{\gamma}_{xy}$ &  membrane plastic strain rates \\ 
	$\eep=[\dot{\bm{\kappa}},\dot{\bm{\varepsilon}}]^{\T}$ & vector of plastic strain rates\\
	$\dot{\kappa}_x, \dot{\kappa}_y, \dot{\kappa}_{xy}$ &  bending plastic strain rates \\ 
	$\lambda$                                         &  plastic load multiplier \\
	$\nu$											  &  Poisson's ratio \\	
	$\bm{\sigma}=[\sigma_x,\sigma_y,\tau_{xy}]^{\T}$  &  in-plane stress \\
    $\sigma_1,\sigma_2$                               &  in-plane eigenstresses \\
    $\dot{\varphi}_x,\dot{\varphi}_y,\dot{\varphi}_z$ &  nodal rotational velocities in the local c.s.  \\
	$\dot{\varphi}_X,\dot{\varphi}_Y,\dot{\varphi}_Z$ &  nodal rotational velocities in the global c.s.  \\
	$\Omega$										  & shell mid-surface \\	
\end{tabular}

\FloatBarrier
 

\section{Introduction}
The prediction of the collapse load (limit load) of a structure made of a material exhibiting plastic behaviour is often of practical interest. The standard approach to obtain such collapse loads is based on iterative calculation schemes using classical nonlinear finite element methods. However, alternative approaches can be applied: the finite-element-based limit analysis (FELA) or the elastic-reduction method (see, r.g. \cite{Mellati}). The former approach is based on the limit theorems of plasticity, first formulated by A.A.Gvozdev in 1936 and later independently by D.C.Prager in 1952. Thereby, the collapse load is obtained as the minimum of a certain non-smooth convex optimization problem, either considering kinematically compatible velocity fields (upper bound approach) or statically admissible stress fields (lower bound approach) within the structure, at the time instant of collapse. Thus, the whole load history does not need to be taken into account, resulting in a stable and numerically efficient approach compared to the standard scheme based on classical finite element formulations. High computational efficiency is further achieved by the fact, that the problem is formulated as the second-order cone programming (SOCP), for which highly efficient primal-dual interior-point optimizers have been developed in last decades (\cite{mosek}). Further advantages of the FELA include a trivial incorporation of multi-surface-plasticity models (no special treatment of singularity points is needed) and an easy incorporation of transversal shear stresses in the failure criterion in the case of beams and shells. The latter is achieved without any negative impact on calculation stability, unlike in the case of FEM (\cite{Kabelac}).
However, the FELA also has some principal disadvantages: it neglects the influence of material elasticity, assumes small deformations, and presumes ideally-plastic material behaviour.
Moreover, it is important to keep in mind that the FELA does not account for a loss of stability due to elastic buckling, which must be checked separately. The assumptions of small deformations and, to some extent, ideal-plasticity can be overcome by the so-called sequential finite-element-based limit analysis (SFELA), as, e.g., shown in (\cite{SFELA1, SFELA2, Seitzberger}). In this process, the FELA method is called repeatedly, with the geometry and plastic strain being updated after each iteration. Updating the plastic strain allows for the incorporation of non-decreasing hardening curves (with zero slope at infinity). 
However, incorporating softening results in a non-convex optimization problem, which is challenging to handle. For details on this approach, often referred to as the extended FELA, see \cite{Tangaramvong}.

In this paper we focus on the upper bound FELA approach applied on shell structures. Much attention in literature is dedicated to solid elements, plane-strain elements (\cite{Makrodimopoulos1, Makrodimopoulos2, Ming1, Ming2}), plate elements (\cite{Bleyer2013,Bleyer2015, Bleyer2016, Meshfree,MeshfreeCook,Le,Leonetti, Makrodimopoulos, Isogeom}) and axisymmetrical shells (\cite{Seitzberger}). Surprisingly, shell elements, which are the most used finite elements in real engineering applications\footnote{From 110 real customer problems from an area of civil engineering published on websites of ADINA Structures (\cite{ADINA:CP}), ATENA (\cite{ATENA:CP}) and RFEM (\cite{RFEM:CP}) finite element programs, 66 projects (i.e., 60\%) contain shell elements.}, seem not to obtain in the context of the FELA a significant attention. Exception to the rule are papers by Bleyer (\cite{Bleyer2016}), Corradi (\cite{Corradi}), Milani (\cite{Milani}) and Tr$\grave{\hat{\text{a}}}$n (\cite{Tran}). 
Bleyer (\cite{Bleyer2016}) uses a 6-node triangle shell element with straight sides, a quadratic approximation of velocities and a linear approximation of rotational velocities, equipped with discontinuities on element boundaries. Corradi (\cite{Corradi}) uses continuous flat triangular shell elements called TRIC. Milani (\cite{Milani}) uses curved 6-node rigid triangular elements with discontinuities on element sides. Tr$\grave{\hat{\text{a}}}$n (\cite{Tran}) uses continuous flat quad Kirchhoff shell elements.

If certain restrictive conditions are met, then finite elements can be developed that guarantee a strict upper bound (\cite{Bleyer2016, Corradi}, see also Section 2.4). This numerically advantageous property, however, restricts us from using most of the standard shell elements used in FEA, particularly all quadrilateral elements. Thus, the main objectives of this work can be introduced as follows:
\begin{enumerate}
	\item To test the DK, DKM and MITC shell element classes for the first time in the context of FELA and SFELA, and to compare discretization errors with the strict upper-bound shell elements used in literature.
	\item To demonstrate that, although these elements cannot guarantee a strict upper bound, discretization errors are significantly lower compared to the strict upper-bound shell elements.
	\item In particular, to show that the use of quadrilateral elements in the context of FELA is beneficial.
\end{enumerate}
The optimization script file was generated by the finite element program femCalc \cite{femCalc}. The optimization problem was solved by the optimizer MOSEK, version 10.0.20 (\cite{mosek}).

The paper is structured as follows: Section 2 presents the kinematic upper bound formulation of FELA and the failure criteria used. Section 3 describes the tested shell elements. Section 4 investigates the convergence behaviour of the considered shell elements on six benchmark problems. Closing remarks are given in Section 5.

\section{Upper bound formulation of the finite-element based limit analysis for shells} 
\subsection{Continuous formulation} 
The upper bound formulation of the FELA is based on minimization of the internal dissipation power summed over a body domain, keeping the power of external forces equal to one. In this way we get a convex, generally non-smooth, optimization problem defined in terms of velocities in the body at the instance of the plastic collapse, which has the plastic load multiplier $\lambda$ as the objective. We consider a linear distribution of the strain across the shell thickness
$\dot{\bm{\epsilon}}(x,y,z)=\dot{\bm{\varepsilon}}(x,y)+z\dot{\bm{\kappa}}(x,y)$, where 
$\dot{\bm{\varepsilon}}=[\dot{\varepsilon}_x,\,\dot{\varepsilon}_y,\,\dot{\gamma}_{xy}]^{\T}$ is the membrane plastic strain rate,  
$\dot{\bm{\kappa}}=[\dot{\kappa}_x,\,\dot{\kappa}_y,\,\dot{\kappa}_{xy}]^{\T}$ is the bending plastic strain rate and $\dot{\bm{\epsilon}}=[\dot{\epsilon}_x,\,\dot{\epsilon}_y,\,\dot{\curlyvee}_{xy}]^{\T}$ is the in-plane plastic strain rate at coordinate $z$.
The primal upper bound formulation in its continuous form has the form 
\begin{align}
	\min_{\dot{\bm{q}}}& ~D_{\p}(\dot{\bm{q}})                          \label{eq1}\\
\text{s.t.}	&~~ \int_{\Omega} \bm{p}^{\T}\dot{\bm{q}}\,\text{d}	
	\Omega + \int_{\Gamma} \bm{f}^{\T}\dot{\bm{q}} \,\text{d}\Gamma =1, \label{eq2}\\
	&~~ \dot{q_i}\Big|_{\Gamma_i}=0~~\forall i \in J \subset \{1,\dots,6\},                       \label{eq3}
\end{align}
where $D_{\p}$ is the dissipation power in the considered body, $\bm{p}=[0,0,0,\,p_x,\,p_y,\,p_z]^{\T}$ is the pressure vector and $\bm{f}=[0,0,0,\,f_x,\,f_y,\,f_z]^{\T}$ is the line-pressure vector acting at the boundary $\Gamma\!=\!\partial\Omega$. $\Gamma_i\subset\Gamma$ is the part of $\Gamma$, at which $\dot{q}_i=0$, which is $i$-th component of the solution vector $\dot{\bm{q}} = 
[\dot{\bm{\varphi}},\, \dot{\bm{u}}]^{\T}=[\dot{\varphi}_X ,\, \dot{\varphi}_Y ,\, \dot{\varphi}_Z ,\, \dot{u}_X ,\, \dot{u}_Y ,\, \dot{u}_Z]^{\T}$.

In the case of the stress-resultant-based failure criteria we have
\begin{align}
	D_{\p}(\eep) =\int_{\Omega} d_{\p}(	\eep(x,y))\,\dOmega, \label{eq4}
\end{align}  
where $\eep(x,y)=[\begin{array}{c} \kp,~\ep \end{array}]^{\T}$ and
\begin{align}
	 	d_{\p}(\eep)&=\max_{\bm{m},\bm{n}~\text{s.t.}~f(\bm{m},\bm{n})\le1} \bm{m}^{\T}\dot{\bm{\kappa}}+\bm{n}^{\T}\dot{\bm{\varepsilon}}, \label{eq5} 
\end{align}
where $\bm{n}=[n_x,\,n_y,\,n_{xy}]^{\T}$ is the vector of internal in-place forces, $\bm{m}=[m_x,\,m_y,\,m_{xy}]^{\T}$ is the vector of internal moments and $f$ is the failure criterion. In the case of the stress-based failure criteria we have
\begin{align}
	D_{\p}(\dot{\bm{\epsilon}}) =\int_{\Omega}\intz d_{\p}^{z}(\epp(x,y,z))\,\dz \dOmega, \label{eq6}
\end{align}
where $\dot{\bm{\epsilon}}(x,y,z)=\dot{\bm{\varepsilon}}(x,y)+z\dot{\bm{\kappa}}(x,y)$ and
\begin{align}
    d_{\p}^z(\epp)&=\max_{\bm{\sigma}~\text{s.t.}~f(\bm{\sigma})\le1} \bm{\sigma}:\epp, \label{eq7}
\end{align}
where $\bm{\sigma}=[\sigma_x,\sigma_y,\sigma_{xy}]^{\T}$ is the stress vector. All failure criteria used in this article, given either by Eq. (\ref{eq5}) or (\ref{eq7}), have explicit forms defined only in terms of plastic strain rates. Note, that these formulations can be easily extended to the Mindlin case, in which the failure criterion $f$ also depends on transversal shear stresses (\cite{Leonetti}). 
 
\subsection{Failure criteria}
In the framework of the FELA, which yields the second order cone optimization problem (\cite{MOSEKCookBook}), only those failure criteria can be used, that belong to the class of associated plasticity and which can simultaneously be formulated in the form of a second order cone\footnote{An $n$-dimensional second order cone is defined as $Q^n=\left\{\bm{x}\in\mathbb{R}^n:~x_1\ge\sqrt{x_1^2+\dots+x_n^2}\right\}$.}. From another perspective, we can divide failure criteria used in shells into criteria defined in terms of stress-resultants and to criteria defined in terms of stresses. Any failure criterion defined in terms of stresses that is to be used in shells requires integration throughout the shell thickness, which increases the number of optimization unknowns needed. Therefore, it is of great benefit to use criteria defined in terms of stress-resultants. In this section, explicit forms of dissipation powers of the failure criteria used in this article are presented. 

\subsection*{Stress-resultant based failure criteria}
\subsubsection*{Approximative Ilyushin failure criterion} 
A.A.Ilyushin proposed in 1948 exact parametrization of the von Mises failure surface (\cite{Ilyushin}), which however cannot be expressed as a second order cone due to its complicated form. A frequently used approximation of the exact Ilyushin parametrization, also proposed by A.A.Ilyushin, is given by intersection of two second order cones:
\begin{align} 
	f&=\min(f_1,f_2), \label{Ilyushin1}\\
	f_i&=
	\frac{\bm{n}^{\T}\bm{P}\bm{n}}{n_{\p}^2}\!+\!\frac{\bm{m}^{\T}\bm{P}\bm{m}}{m_{\p}^2}\!+\!
	\frac{(-1)^i}{\sqrt{3}}\frac{\bm{m}^{\T}\bm{P}\bm{n}}{n_{\p}m_{\p}} \le 1,~i=1,2,  \\ 
	\bm{P}&=\left[
	\begin{array}{ccc}
		1 & -\frac{1}{2}& 0\\
		-\frac{1}{2}& 1 & 0\\
		0 & 0 & 3
	\end{array}
	\right],\label{Ilyushin3}
\end{align}
where $n_{\p}=\sigma_0t$ is the maximum plastic in-plane force and $m_{\p}=\sigma_0\frac{t^2}{4}$ is the maximum plastic moment. Note, that if one restricts to the pure in-plane situation ($\bm{m}=\bm{0}$) or to the pure bending situation ($\bm{n}=\bm{0}$), then Eqs. (\ref{Ilyushin1}-\ref{Ilyushin3}) coincide with the von Mises failure criterion. We decompose the matrix $\bm{P}$ as follows: 
\begin{align}
	\bm{P}=\bm{R}^{\T}\bm{R},~\bm{R}=
	\left[
	\begin{array}{ccc}
		~~\frac{\sqrt{3}}{2} & 0 & 0\\
		-\frac{1}{2} & 1 & 0\\
		0 & 0 & \sqrt{3}
	\end{array}\right].
\end{align}	
The final expression of the dissipation power area-density has the form (\cite{Bleyer2016}) 
\begin{align} 
	d_{\p}(\eep)&=\sum_{i=1}^2||\mathcal{\tilde{R}}_i^{-\T} \eep_i||, \label{IL1}\\
	~~\eep&=
	\left[ \begin{array}{c} \bm{\dot{\kappa}} \\ \bm{\dot{\varepsilon}}  \end{array} \right]
	=\eep_1+\eep_2=	
	\left[ \begin{array}{c} \bm{\dot{\kappa}}_1 \\ \bm{\dot{\varepsilon}}_1  \end{array} \right]+	
	\left[ \begin{array}{c} \bm{\dot{\kappa}}_2 \\ \bm{\dot{\varepsilon}}_2  \end{array} \right], \\
	\mathcal{\tilde{R}}_{i}^{-\T}&=\left[
	\begin{array}{cc}
		\bm{0}   &   n_{\p} \bm{R}^{-\T}    \\
		\frac{m_{\p}}{\sqrt{1-\alpha^2}} \bm{R}^{-\T} & -\frac{\alpha n_{\p}}{\sqrt{1-\alpha^2}}\bm{R}^{-\T}   \\
	\end{array}
	\right],~\alpha=\frac{(-1)^i}{2\sqrt{3}},~i=1,2, \\
	~\bm{R}^{-\T}&=\left[
	\begin{array}{ccc}
		\frac{2}{\sqrt{3}} & \frac{1}{\sqrt{3}} &  0 \\
		0&1&0\\
		0&0&\frac{1}{\sqrt{3}} \\   
	\end{array}
	\right]. \label{IL2}
\end{align}
\subsubsection*{Von Mises failure criterion for in-plane loading} 
The von Mises failure criterion for the case of in-plane loading is the special case of the Ilyushin criterion for $\bm{m}=0$: $f=\frac{\bm{n}^{\T}\bm{P}\bm{n}}{n_{\p}^2}\le1$.
The final expression of the dissipation power area-density has the form $d_{\p}(\bm{\dot{\varepsilon}})=n_{\p}||\bm{R}^{-\T} \bm{\dot{\varepsilon}}||$.
\subsubsection*{Von Mises failure criterion for pure bending} 
The von Mises failure criterion for pure bending is the special case of the Ilyushin criterion for $\bm{n}=0$: $f=\frac{\bm{m}^{\T}\bm{P}\bm{m}}{m_{\p}^2}\le1$. The final expression of the dissipation power area-density has the form $d_{\p}(\bm{\dot{\kappa}})=m_{\p}||\bm{R}^{-\T} \bm{\dot{\kappa}}||$.
\subsection*{Stress based failure criteria}
\subsubsection*{Plane-stress Tresca failure criterion}
The Tresca in-plane plane stress failure criterion is given by
\begin{align}
	f=\frac{\max(|\sigma_2-\sigma_1|,|\sigma_1|,|\sigma_2|)}{\sigma_0}\le1,
\end{align}
where $\sigma_1, \sigma_2$ are eigenstresses. The final expression of the dissipation power density has the form 
\begin{align}
	d_{\p}^z(\dot{\bm{\epsilon}})&=f_{\text{y}}\max\left[\left|\frac{\dot{\epsilon}_x+\dot{\epsilon}_y-\sqrt{A}}{2}\right|, \left|\frac{\dot{\epsilon}_x+\dot{\epsilon}_y+\sqrt{A}}{2}\right|,|\dot{\epsilon}_x+\dot{\epsilon}_y|\right],  \label{eq447} \\
	A&=(\dot{\epsilon}_x - \dot{\epsilon}_y)^2+\dot{\curlyvee}_{xy}^2. \label{eq447b}
\end{align}
For the derivation and implementation details see the Appendix.
\subsubsection*{Plane-strain Tsai-Wu failure criterion}
In the plane-strain approximation we consider $\dot{\varepsilon}_z=\dot{\gamma}_{xz}=\dot{\gamma}_{yz}=0$. This criterion is extracted from the three-dimensional Tsai-Wu criterion under the assumption of transversal isotropy ($f_{\t y}\!=\!f_{\t z},\,f_{\c y}\!=\!f_{\c z}$), under the assumption of a zero coefficient of the cross-term $\sigma_x\sigma_y$, and under the assumption of a negligible influence of transversal shear stresses: 
\begin{align}
	&f=\sigma_x\left[\frac{1}{f_{\t x}}-\frac{1}{f_{\c x}}\right]+
	(\sigma_y+\sigma_z)\left[\frac{1}{f_{\t y}}-\frac{1}{f_{\c y}}\right]+ \frac{\sigma_x^2}{f_{\t x}f_{\c x}}+\frac{\sigma_y^2+\sigma_z^2}{f_{\t y}f_{\c y}} + \frac{\tau_{xy}^2}{f_{\vv xy}^2}\le 1. 
	\label{functionf}
\end{align}
The final expression of the dissipation power density has the form 
\begin{align}
	d_{\p}^z(\dot{\bm{\epsilon}})=&\frac{f_{\t x}-f_{\c x}}{2}\dot{\epsilon}_x+\frac{f_{\t y}-f_{\c y}}{2}\dot{\epsilon}_y + \sqrt{(\bm{R}^{-\T}\dot{\bm{\epsilon}})^{\T}\bm{R}^{-\T}\dot{\bm{\epsilon}}},
\end{align}
where
\begin{align}
	\bm{R}^{-\T}&=\left[
	\begin{array}{ccc}
		\sqrt{\chi f_{\t x}f_{\c x}} & 0 & 0\\
		0 & \sqrt{\chi f_{\t y}f_{\c y}} & 0\\
		0 & 0 &\sqrt{\chi}f_{\vv xy}\\
	\end{array}\right],\\
\chi&=1+\frac{1}{4}\frac{(f_{\t x}-f_{\c x})^2}{f_{\t x}f_{\c x}}+\frac{1}{2}\frac{(f_{\t y}-f_{\c y})^2}{f_{\t y}f_{\c y}}.
\end{align}
For the derivation see the Appendix.
\subsection{Discretized formulation}
We discretize at first the formulation based on the Ilyushin failure criterion, i.e. Eqs. (\ref{eq1}-\ref{eq4},\ref{IL1},\ref{IL2}). We consider continuous shell elements with six degrees of freedom per node. The vector of unknown velocities is denoted by  $\dot{\bm{\mathfrak{q}}}$. The integral in Eq. (\ref{eq4}) is evaluated by means of the Gauss quadrature in $N_q$ quadrature point within each element. The strain matrix in element $e$ and quadrature point $q$ of dimensions $6\times \Nd$ is denoted by $\bm{B}_{eq}$ ($\Nd$ is a total number of dofs; the matrix $\bm{B}_{eq}$ has however at most $6n$ nonzero columns, where $n$ is a number of element nodes). The corresponding strain vector of length $6$ is denoted by $\eep_{eq}$. The primal formulation in the discretized form yields
\begin{align}
	\min_{\dot{\bm{\mathfrak{q}}}}&~\,\sum_{e=1}^{N_e} \sum_{q=1}^{N_q} J_{eq} w_q \sum_{i=1}^2 ||\mathcal{\tilde{R}}_i^{-\T}\bm{\eep}_{ieq}||   , \label{eq10}\\
	\text{s.t.}&~~\sum_{i=1}^2 \eep_{ieq} =\bm{B}_{eq}\dot{\bm{\mathfrak{q}}},~~e=1,\dots,N_e,~q=1,\dots,N_q, \label{eq11}\\ 
	&~~ \bm{\mathfrak{f}}^{\T}\dot{\bm{\mathfrak{q}}}=1, \label{eq12} \\
	&~~ \bm{H}\dot{\bm{\mathfrak{q}}}=0, \label{eq13}
\end{align} 
where $J_{eq}$ is the Jacobian in the element $e$ and in the Gauss point $q$, $w_q$ is the Gauss quadrature weight and $\bm{H}$ is the matrix prescribing zero boundary conditions. Equivalently
\begin{align}
	\min_{\dot{\bm{\mathfrak{q}}}}&~\,\sum_{e=1}^{N_e} \sum_{q=1}^{N_q} J_{eq} w_q \sum_{i=1}^2 s_{ieq}, \label{eq14}\\
\text{s.t.}&~~ ||\mathcal{\tilde{R}}_i^{-\T}\bm{\eep}_{ieq}|| \le s_{ieq},~~i=1,2,~~e=1,\dots,N_e,~q=1,\dots,N_q,\\ 
&~~\sum_{i=1}^2\eep_{ieq}=\bm{B}_{eq}\dot{\bm{\mathfrak{q}}},~~e=1,\dots,N_e,~q=1,\dots,N_q, \label{eq15}\\
&~~ \bm{\mathfrak{f}}^{\T}\dot{\bm{\mathfrak{q}}}=1, \label{eq16} \\
&~~ \bm{H}\dot{\bm{\mathfrak{q}}}=0. \label{eq17}
\end{align}
By substitution $\bm{e}_{ieq}=\mathcal{\tilde{R}}_i^{-\T}\bm{\eep}_{ieq}$, we get the final form of the primal formulation, denoted by (P):
\begin{align}
	\min_{\dot{\bm{\mathfrak{q}}}}&~\,\sum_{e=1}^{N_e} \sum_{q=1}^{N_q} J_{eq} w_q \sum_{i=1}^2 s_{ieq}, \label{eq18}\\
	\text{s.t.}&~~\sum_{i=1}^2\mathcal{\tilde{R}}_i^{\T}\bm{e}_{ieq}-\bm{B}_{eq}\dot{\bm{\mathfrak{q}}}=\bm{0},~~e=1,\dots,N_e,~q=1,\dots,N_q, \label{eq19}\\ 
	&~~ ||\bm{e}_{ieq}|| \le s_{ieq},~~i=1,2,~~e=1,\dots,N_e,~q=1,\dots,N_q, \label{eq20}\\
	&~~ \bm{\mathfrak{f}}^{\T}\dot{\bm{\mathfrak{q}}}=1, \label{eq16b} \\
	&~~ \bm{H}\dot{\bm{\mathfrak{q}}}=0. \label{eq21}
\end{align}
The dual formulation of the primal problem, denoted by (D), is computationally more efficient (\cite{Makrodimopoulos}):
\begin{align}
	\max_{\bm{m}}~&\beta \label{D0}\\
	\text{s.t.}&~~\mathcal{\tilde{R}}_i\bm{m}_{eq}+\bm{s}_{ieq}=\bm{0},~~i=1,2,~~e=1,\dots,N_e,~~q=1,\dots,N_q,\\
	&~~||\bm{s}_{ieq}||\le s_{ieq}^*:=J_{eq}w_q,~~i=1,2,~~e=1,\dots,N_e,~~q=1,\dots,N_q, \label{eq_yield} \\
	&~~ \sum_{e=1}^{N_e} \sum_{q=1}^{N_q} \bm{B}^{\T}_{eq} \bm{m}_{eq} - \bm{H}^{\T}\bm{r} - \beta\bm{\mathfrak{f}} =\bm{0}, \label{eq_sily} 
\end{align}
where Eq. (\ref{eq_yield}) defines yield surfaces at each element, and each quadrature point, and Eq. (\ref{eq_sily}) equates nodal force components at each degree of freedom (dof). The term $\bm{H}^{\T}\bm{r}$ represents reaction forces in fixed dofs. If this information is not needed, the fixed dofs can be excluded from the calculation. The term $\bm{\mathfrak{f}}$ represents a vector of external nodal forces. In the case of a single yield surface, the internal force $\bm{m}_{eq}$ can be eliminated from the system (by using $\bm{m}_{eq}=-\mathcal{\tilde{R}}^{-1}\bm{s}_{eq}$), which further simplifies the dual formulation (D):
\begin{align}
	\max_{\bm{s}}~&\beta \label{D1}\\
 	 \text{s.t.}&~~||\bm{s}_{eq}||\le s_{eq}^{*}:=J_{eq} w_q,~~e=1,\dots,N_e,~~q=1,\dots,N_q,\\
	&~~ \sum_{e=1}^{N_e} \sum_{q=1}^{N_q} \bm{B}^{\T}_{eq} \mathcal{\tilde{R}}^{-1}\bm{s}_{eq}  + \bm{H}^{\T}\bm{r} + \beta\bm{\mathfrak{f}} =\bm{0}. \label{DN}
\end{align} 
Let us compare the number of dofs in the primal and dual formulations for the clamped pressure-loaded square plate problem. Please note that the number of dofs varies in every problem due to different boundary conditions. We consider only one quarter of the plate with $N_e$ shell elements due to symmetry, 4-node quad elements, the von Mises failure criterion, and 3 dofs per node. The comparison given in Table \ref{Tab:PD} shows a smaller number of dofs for the dual formulation.
 \begin{center}
	\begin{table}
		\centering
		\begin{tabular}{ccc}
			\hline
			Formulation  & \#dofs & Comment  \\  
			\hline
			\hline
			(P) &  $\underbrace{4N_qN_e}_{\#\bm{e}_{eq} + \#\bm\dot{\bm{\mathfrak{q}}}}+\underbrace{3N_e-2\sqrt{N_e}}_{\#\bm\dot{\bm{\mathfrak{q}}}}$ & zero dofs of $\bm\dot{\bm{\mathfrak{q}}}$ excluded \\
		   \hline
			(D) & $\underbrace{4N_eN_q}_{\#\bm{s}_{eq} + \#s_{eq}^{*}}+\underbrace{1}_{\beta}$ &  $-$\\
			\hline  
		\end{tabular} \nonumber
		\captionof{table}{\#dofs for the clamped pressure loaded square plate problem \label{Tab:PD}}
	\end{table}
\end{center}
The primal formulation based on the stress-based failure criterion is discretized as follows:
\begin{align}
	\min_{\dot{\bm{\mathfrak{q}}}}&~\,\frac{h}{2}\sum_{e=1}^{N_e} \sum_{q=1}^{N_q} \sum_{r=1}^{N_r} J_{eq} w_q w_r d_{\p}^{z}(\dot{\bm{\epsilon}}_{eqr}) , \label{eq110}\\
	\text{s.t.}&~~\dot{\bm{\epsilon}}_{eqr}=\dot{\bm{\varepsilon}}_{eq}+z_r\dot{\bm{\kappa}}_{eq},~~e=1,\dots,N_e,~~q=1,\dots,N_q,~~r=1,\dots,N_r,\\
	&~~\eep_{eq}=\bm{B}_{eq}\dot{\bm{\mathfrak{q}}},~~e=1,\dots,N_q,~q=1,\dots,N_q, \label{eq111}\\ 
	&~~ \bm{\mathfrak{f}}^{\T}\dot{\bm{\mathfrak{q}}}=1, \label{eq112} \\
	&~~ \bm{H}\dot{\bm{\mathfrak{q}}}=0, \label{eq113}
\end{align} 
where $w_r$ is the Gauss quadrature weight for integration across the shell thickness and $\frac{h}{2}$ is the corresponding Jacobian.

\subsection{Strict upper-bound versus approximative upper bound}

The strict upper bound formulation needs to avoid numerical integration across the element, as numerical integration has no bounding effect. The solution to this problem is to consider elements with exactly quadratic velocity interpolation, i.e., having a constant strain approximation, and at the same time having a constant Jacobian, which necessitates the use of triangle elements with straight sides. Then, the integrand (\ref{eq4}) or (\ref{eq6}) is constant and the one-point Gauss integration rule is exact. Such shell elements guarantee the strict upper bound formulation. However, this solution restricts us from using quadrangle elements or different numerical vehicles, like bubble functions for example, which can all increase numerical efficiency. Moreover, the quadratic velocity field cannot be used to construct conforming Kirchhoff elements, which require $\mathcal{C}^1$ continuity between elements. As a result, conforming Kirchhoff elements can only guarantee an approximate upper bound formulation. In this article, we compare different shell elements, which have no bounding effect but exhibit smaller numerical errors and therefore higher efficiency in comparison with the strict upper bound shell elements. 

The first-order optimality condition of the maximization (\ref{eq5}) or (\ref{eq7}) immediately yields the flow rule\footnote{For example, in the case of Eq.(\ref{eq7}) we get the Lagrangian $\mathcal{L}=\bm{\sigma}:\dot{\bm{\varepsilon}}-\lambda(f(\bm{\sigma}) - 1)$,
and the first order optimality condition yields
$\PART{\mathcal{L}}{\bm{\sigma}}=\dot{\bm{\varepsilon}}-\lambda\PART{f(\bm{\sigma})}{\bm{\sigma}}=0$.}. If the strict upper bound of the formulation is not guaranteed, then the flow rule is not fulfilled point-wise but only on average for each element.

\subsection{Sequential finite-element-based limit analysis (SFELA)}
Two basic assumptions of the FELA (the assumption of small deformations and the assumption of ideal plasticity) are both restrictive for real engineering applications. The SFELA method is a modification of the FELA method, first introduced by Yang (\cite{SFELA1}), which overcomes, to some extent, both of these problems. The idea of the method is based on the time integration of nodal velocities  
\begin{align}
	\triangle \bm{u}_i=\triangle t \dot{\bm{u}}_i,~~i=1,\dots,\Nd,
\end{align}
for some time step $\triangle t$. In this way, we get a new deformed mesh on which the FELA can be applied again to obtain the new velocity field. By repeating this procedure, we get the time-dependent displacement field of the plastic collapse. By updating the yield stress at each time step according to the current value of the plastic strain, hardening or even softening can be incorporated into the calculation. For further details, see \cite{SFELA1,SFELA2,StemberaMilovy}.
\subsection{Practical numerical experiences with optimization} 
The FELA formulation results in a non-smooth convex second-order cone optimization problem, which is solved by the MOSEK optimizer (\cite{mosek}). As an input to MOSEK a script file in the \emph{opf} file format is used. The following suggestions help to achieve calculation efficiency and accuracy:
\begin{enumerate}
	\item{We eliminate zero dofs of $\dot{\bm{\mathfrak{q}}}$ from the primal formulation (P), rather than just setting them to zero. }
	\item{We eliminate as many variables as possible from the formulation if their values are not explicitly needed. For example, eliminating $\bm{m}_{eq}$ in Eqs. (\ref{D0}--\ref{eq_sily}) leads to Eqs.(\ref{D1}--\ref{DN}).} 
	\item We use parallelization for creating the \emph{opf} script file. 
	\item {We use parallelization in MOSEK with the setting \mbox{[p MSK\_IPAR\_NUM\_THREADS]} \mbox{16} \mbox{[/p]}. Additionally, we use the setting \mbox{[p MSK\_IPAR\_INTPNT\_SOLVE\_FORM]} \mbox{MSK\_SOLVE\_PRIMAL} \mbox{[/p]}.}
	\item 
	It is important to remember that round-off errors can disrupt convergence if the primal variables differ by many orders of magnitude, which can easily occur in practice. MOSEK automatically handles the scaling between primal and dual variables, but the scaling of the primal variables is left to the user. Therefore, we scale the problem variables by reducing plastic strengths by $n$ orders of magnitude to keep the underlying optimization problem well-conditioned. Consequently, the resulting plastic multiplier $\lambda$ must be increased by $n$ orders of magnitude. In our case, we use $n=9$ in example 4.5 and $n=6$ in all other examples.
\end{enumerate}

\section{Finite elements used}  
We use continuous shell elements used for the standard FEM analysis, which have six degrees of freedom per node (three displacements and three rotations). In particular, we use plate elements, membrane elements and shell elements, all listed in Table \ref{Tab:Elements}. The triangular plate elements (the Kirchhoff plate elements DKT, Specht and the Mindlin plate element DKMT) must be equipped with membrane elements (CST or OPT) to function as shells. Analogously, the quad plate elements (the Kirchhoff plate elements DKQ and the Mindlin plate element DKMQ) must be equipped with membrane elements (CSQ or BSQ) to function as shells. The DKQ and DKMQ elements only work in the planar case, so they must be further equipped with the warping matrix (the rigid links correction) to also work in the non-planar case (\cite{warp}). No special treatment of the drilling rotations is used, i.e. no stabilization energy term is used. 

Regarding the numerical integration of the stiffness matrix, we use 3-point quadrature for triangles and 4-point quadrature for quadrangles, except for the MITC9 and MIC16 shell elements. In the case of MITC9 shell element, we use an 8-point economical quadrature for its efficiency (\cite{Solin}), instead of the standard composite $3\!\times\!3$ quadrature. Similarly, in the case of the MITC16 element, we use a 12-point economical quadrature (\cite{Solin}), instead of the standard composite $4\!\times\!4$ quadrature.

The Reisner-Mindlin elements (DKMT, DKMQ, DKMQ24, DKMQ24+) need to be adapted for the plastic-rigid calculation. In the derivation of these shell elements, there are two places where elastic properties are used: 
\begin{enumerate}
	\item Poisson's ratio is used in the shear approximation.
	\item The elastic properties are used in the membrane stiffness matrix, which is used in the static condensation of the DKMQ24+ element.
\end{enumerate}
Let us discuss the choice of these (purely numerical) parameters in two different cases: isotropic and orthotropic. We start with the isotropic in-plane von Mises criterion for which the dissipation energy takes the following form:   
\begin{align} 
	D_{\p}(\dot{\bm{\varepsilon}})&=\int_A ~h\sigma_0\sqrt{\left(\bm{\dot{\varepsilon}}_{\m}\right)^{\text{T}}\bm{Q}~\dot{\bm{\varepsilon}}_{\m}}~\text{d}A,~\bm{Q}=\left[\begin{array}{ccc}
		\frac{4}{3}&\frac{2}{3}& 0\\
		\rule[2mm]{0mm}{2mm}
		\frac{2}{3}&\frac{4}{3}& 0\\
		0 & 0 & \frac{1}{3}
	\end{array}\right],~\bm{\dot{\varepsilon}}_{\m}=\left[\begin{array}{c}
	\dot{\varepsilon}_{x} \\
	\dot{\varepsilon}_{y} \\
	\dot{\gamma}_{xy}\\
\end{array}\right],
\end{align}	
which operates on the plastic strain rates. In the framework of the linear FEM the elastic stiffness matrix has the form:
\begin{align}
	 E&=\frac{1}{2}\int_A (\bm{\varepsilon}_{\m})^{\T}\bm{D}\bm{\varepsilon}_{\m}\,\text{d}A,~~\bm{D}=hE\left[
	\begin{array}{ccc}
		\frac{1}{1-\nu^2} & \frac{\nu}{1-\nu^2} & 0 \\
		\frac{\nu}{1-\nu^2} & \frac{1}{1-\nu^2} & 0 \\
		0 &                     0 & \frac{1}{2(1+\nu)} \\
	\end{array}\right],~\bm{\varepsilon}_{\m}=\left[\begin{array}{c}
	\varepsilon_x \\
	\varepsilon_y \\
	\gamma_{xy}\\
\end{array}\right],
\end{align}
which operates on the elastic strains. This leads us to set 
\begin{align}
	E=\sigma_0,~\nu=0.5,
\end{align}
which makes matrices $\bm{Q}$ and $\bm{D}$ identical. The von Mises criterion, defined for the case of pure bending, would yield the same setting. In the orthotropic case, which is represented by the Tsai-Wu failure criterion in the plane-strain approximation where transversal isotropy is considered, we obtain
\begin{align}
		D_{\p}(\dot{\bm{\varepsilon}})&=h\int_A \left[\frac{f_{\text{t}x} - f_{\text{c}x}}{2}\dot{\varepsilon}_x+\frac{f_{\text{t}y} - f_{\text{c}y}}{2}\dot{\varepsilon}_y+ \sqrt{\left(\bm{\dot{\varepsilon}}_{\m}\right)^{\text{T}}\bm{Q}~\dot{\bm{\varepsilon}_{\m}}}\right]~\text{d}A,\\
	\bm{Q}&=h\chi\left[
	\begin{array}{ccc}
		f_{\text{t} x} f_{\text{c} x} & 0 & 0\\	
		0 & f_{\text{t}y} f_{\text{c} y} & 0\\
		0 & 0 & f_{\text{v}xy}^2
	\end{array}\right],~\chi=1+\frac{1}{4}\frac{(	f_{\text{t} x}-f_{\text{c} x})^2}{	f_{\text{t} x} f_{\text{c} x}}+\frac{1}{2}\frac{(	f_{\text{t} y}-f_{\text{c} y})^2}{	f_{\text{t} y} f_{\text{c} y}}.
\end{align}
The elastic stiffness matrix in the orthotropic case has the form: 
\begin{align}
	E&=\frac{1}{2}\int_A (\bm{\varepsilon}_{\m})^{\T}\bm{D}\bm{\varepsilon}_{\m}\,\text{d}A,~~\bm{D}=h\left[
	\begin{array}{ccc}
 \frac{E_x^2(\nu_y-1)}{E_x(\nu_y-1)+2\nu_x^2E_y} & \frac{-E_yE_x\nu_x}{E_x(\nu_y-1)+2\nu_x^2E_y}                     & 0 \\
	\text{symm}.                                  & \frac{-E_y(E_x-\nu_x^2E_y)}{[E_x(\nu_y-1)+2\nu_x^22E_y](1+\nu_y)} & 0 \\
      0                                           &  0                                                                & G_{xy} \\
	\end{array}\right]. 
\end{align}
This leads us to set 
 \begin{align}
 \nu_x=\nu_y=0,~E_x=\sqrt{\chi f_{\text{t}x} f_{\text{c}x}},~E_y=\sqrt{\chi f_{\text{t}y} f_{\text{c}y}},~G_{xy}=f_{\text{v}xy}\sqrt{\chi},
 \end{align}  
which makes matrices $\bm{Q}$ and $\bm{D}$ identical. 

\section{Convergence tests}

All the test problems presented in this section are listed in Table \ref{Tab:Benchmarks}, while all finite elements used in this article are listed in Table \ref{Tab:Elements}. The standard method for evaluating convergence behavior is to depict a relative error with respect to the number of dofs. However, this can be partially misleading because optimization time also depends on a number of auxiliary variables (for example, whether the strain rate is eliminated from the formulation or not, and on the number of Gauss points used, which varies for different elements). Therefore, we also provide diagrams of the relative error with respect to the actual optimization time. It is important to note that this optimization time refers to the time needed by the MOSEK optimizer and does not include the time required for generating the script file. A computer equipped with an i7-10875H CPU and 32GB of RAM was used for the calculation.
\begin{center}
	\begin{table}
		\centering
		\begin{tabular}{clcccc}
			\hline
			 \# & Test problem &  Element type    & Material model         & Analysis   & Ref.\\ 
			\hline
			\hline
1. & Clamped plate                    & plate     & von Mises               & FELA & \cite{Bleyer2013, Makrodimopoulos, Leonetti}\\
2. & Cook's membrane                  & membrane  & von Mises (plane stress)& FELA & \cite{Ciria}\\
3. & Simply support. spherical cap    & shell     & Tresca                  & FELA & \cite{Bleyer2016} \\
4. & Clamped cylindrical tube         & shell     & von Mises               & FELA & \cite{Bleyer2016} \\
5. & Block with asymmetric holes      & membrane  & Tsai-Wu (plane strain)  & FELA & \cite{Zouain, Makrodimopoulos1, Ming2} \\
6. & Simply supported plate           & plate     & approximated Ilyushin   & SFELA &  \cite{Corradi}\\
			\hline 
		\end{tabular} \nonumber
		\captionof{table}{Test problems. \label{Tab:Benchmarks}}
	\end{table}
\end{center}


\begin{center}
	\begin{table}
		\centering
		\begin{adjustbox}{angle=270}
			\begin{tabularx}{\tabulkywidthB}{ABCDEFGH} 
				\hline
				Element \newline notation   &  Element shape   & Element type     & \#qua-\newline drature \newline  points & Transv. \newline ~~~shear \newline approx. & Implementation & Reference  & Comment \\ 
				\hline
				\hline
				Ciria  & triangle  & membrane& ~~~1 & ~~~~---     &  from Reference & \cite{Ciria}                      & \scalebox{0.97}[1.0]{strict upper-bound; continuous 3-node triangle} \\	
				\scalebox{0.75}[1.0]{Makrodimopoulos} & triangle & plate & ~~~1 & Kirchhoff &  from Reference & \cite{Makrodimopoulos}  & \scalebox{0.97}[1.0]{strict upper-bound; continuous 6-node triangle} \\	
				T6   & triangle  & plate   & ~~~1 & Kirchhoff   &  from Reference & \cite{Bleyer2013} & \scalebox{0.97}[1.0]{strict upper-bound; continuous 6-node triangle} \\ 
                \hline		
				T6b   & triangle  & plate   & ~~~1 & Kirchhoff   &  from Reference & \cite{Bleyer2013} & 6-node triangle with a cubic bubble\\   
				TRIC   & triangle  & shell   & ~~~3 & Mindlin     &  from Reference &\cite{TRIC}                        &  =TRIangular Composite\\	
				\hline		
				Specht & triangle  & plate   & ~~~3   & Kirchhoff & our implementation & \cite{Specht} &  \\
				DKT    & triangle  & plate   & ~~~3   & Kirchhoff & our implementation & \cite{DKT}    & =Discrete Kirchhoff Triangle  \\  		
				DKQ    & quad      & plate   & ~~~3   & Kirchhoff & our implementation & \cite{DKQ}    & =Discrete Kirchhoff Quadrangle,  \mbox{uses the warping matrix (\cite{warp})} \\  			
				\hline	
				CST    & triangle  & mebrane & ~~~3   &~~~~---    & our implementation & \cite{CST} & =Constant Strain Triangle      \\  			
				OPT    & triangle  & mebrane & ~~~3   &~~~~---    & our implementation & \cite{OPT} & =OPtimal Triangle   \\  
				BSQ    & quad      & mebrane & ~~~4   &~~~~---    & our implementation & \cite{BSQ} & =Bilinear Strain Quadrangle  \\  
				OPQ    & quad      & mebrane & ~~~4   &~~~~---    & our implementation & \cite{OPQ} & =OPtimal Quadrangle \\ 			
				\hline	
				DKMT   & triangle  & plate   & ~~~3   & Mindlin   & our implementation & \cite{DKMT}   & =Discrete Kirchhoff-Mindlin Triangle   \\  
				DKMQ   & quad      & plate   & ~~~3   & Mindlin   & our implementation & \cite{DKMQ}   & =Discrete Kirchhoff-Mindlin Quadrangle, \mbox{uses the warping matrix (\cite{warp})}  \\  DKMQ24     & quad & shell    & ~~~4      & Mindlin          & our implementation & \cite{Katili2015} & \\ 
				DKMQ24+& quad      & shell   & ~~~4   & Mindlin   & our implementation & \cite{Stembera2020}  & DKMQ24$_2$+ uses the $2\times2$ Gauss quadrature, without nodal moment corrections \\  
				\hline	
				MITC3  & triangle  & shell   & ~~~3   & Mindlin   & our implementation & \cite{Chapelle} & =Mixed Interpolation Tensorial Components \\
				MITC4  & quad      & shell   & ~~~4   & Mindlin   & our implementation & \cite{MITC4}    & \\
				MITC9  & quad      & shell   & ~~~8   & Mindlin   & our implementation & \cite{Chapelle} & uses 8-point economical quadrature (\cite{Solin}) \\
				MITC16 & quad      & shell   & ~~12   & Mindlin   & our implementation & \cite{Chapelle} & uses 12-point economical quadrature (\cite{Solin}) \\
				\hline  
			\end{tabularx} \nonumber
		\end{adjustbox}
		\captionof{table}{Used finite elements. All elements in our implementation were implemented in the program femCalc (\cite{femCalc}) and tested on standard linear elastic benchmarks. Only first two elements  guarantee the strict upper-bound formulation. \label{Tab:Elements}}
	\end{table}
\end{center}
\subsection{Clamped pressure loaded square plate}
\label{b1}
 
\begin{figure}[H]
	\centering
	\begin{overpic}[width=0.35\textwidth]{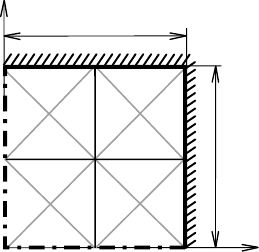} 
		\put(35,87.5){$\frac{L}{2}$}
		\put(86,35){$\frac{L}{2}$}
		\put(39.5,36.7){	\scalebox{0.85}[1.0]{thickness\,$h$}}
		\put(95.5,5){$X$}
		\put(-6,92.5){$Y$}
	\end{overpic}
	\caption{A quarter of the clamped pressure loaded square plate with quad mesh $2\times 2$, a triangular mesh depicted in a grey colour. Parameters: $L=10\,\text{m},\,h=0.001\,\text{m},\,f_{\text{y}}=4\,\text{MPa}$, loading pressure $\,p=1\,\text{Pa}$.}
	\label{ClampedSquarePlate}
\end{figure}
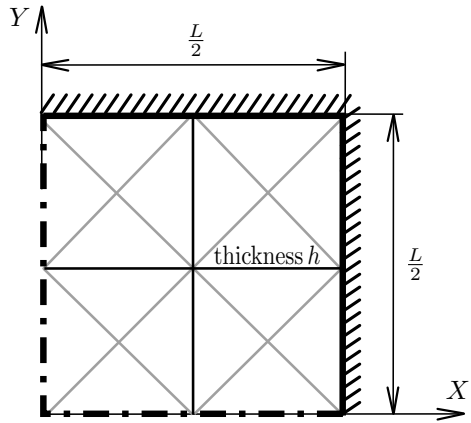

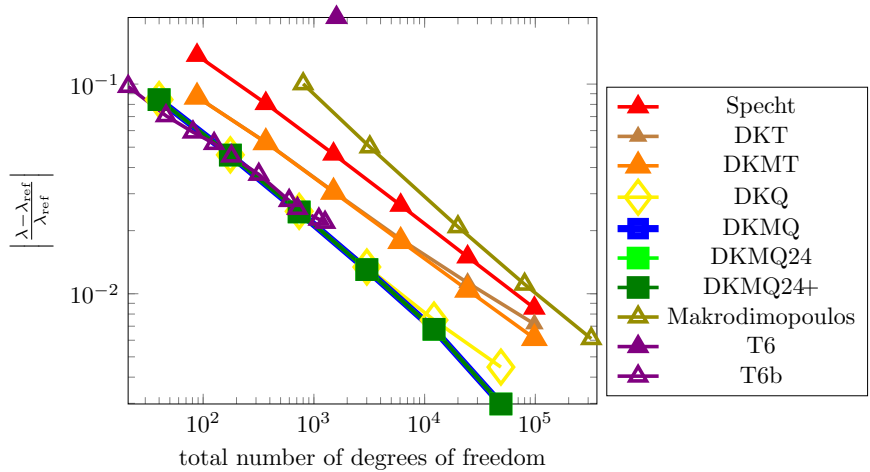
\begin{figure}
	\centering
	\scalebox{.9}
	{
		\begin{tikzpicture}
			\begin{loglogaxis}[
				xlabel=\text{total number of degrees of freedom},
				ylabel= $\left|\frac{\lambda-\lambda_{\text{ref}}}{\lambda_{\text{ref}}}\right|$,
				x tick label style={
					/pgf/number format/1000 sep={}  
				},
				y tick label style={
					/pgf/number format/1000 sep={},
					/pgf/number format/fixed,       
					/pgf/number format/precision=5  
				},
				scaled y ticks=false, 
				xmin=21,
				xmax=360000,
				ymin=0.00297978262349399,
				ymax=0.2078237,
				legend style={at={(1.02,0.02)},anchor=south west}
				]
				\addplot[line width=1.5pt, color=red, mark=triangle*, mark size=4pt, solid] coordinates {
					(88, 0.137119920128883)
					(368, 0.0807743414036442)
					(1504, 0.0463529078078555)
					(6080, 0.0264232488484492)
					(24448, 0.014984139229879)
					(98048, 0.00853447845522001)
				};
				\addplot[line width=1.5pt, color=brown, mark=triangle*, mark size=3pt, solid] coordinates {
					(88, 0.086875541739466)
					(368, 0.0527735245399469)
					(1504, 0.0307212906446417)
					(6080, 0.0181958203807492)
					(24448, 0.0112857434594178)
					(98048, 0.00717126908851621)
				};
				\addplot[line width=1.5pt, color=orange, mark=triangle*, mark size=5pt, solid] coordinates {
					(88, 0.0868184520432938)
					(368, 0.0526607746590729)
					(1504, 0.0305427378548253)
					(6080, 0.0177923124049829)
					(24448, 0.010420843638674)
					(98048, 0.00612856526967849)
				};
				\addplot[line width=1.5pt, color=yellow, mark=diamond, mark size=7pt, solid] coordinates {
					(40, 0.0843418120759684)
					(176, 0.0459650790769441)
					(736, 0.0247313879875655)
					(3008, 0.0134116539220803)
					(12160, 0.0074913435138753)
					(48896, 0.00446788137323869)
				};
				\addplot[line width=3pt, color=blue, mark=square, mark size=3pt, solid] coordinates {
					(40, 0.0842890109142067)
					(176, 0.0458649225114021)
					(736, 0.0245538108960542)
					(3008, 0.0130362142905765)
					(12160, 0.0067667627237866)
					(48896, 0.00297978262349399)
				};
				\addplot[line width=1.5pt, color=green, mark=square*, mark size=4pt] coordinates {
					(40, 0.0842883528851173)
					(176, 0.0458649225114021)
					(736, 0.0245538108960542)
					(3008, 0.0130362142905765)
					(12160, 0.0067667627237866)
					(48896, 0.00297978262349399)
				};
				\addplot[line width=1.5pt, color=col1, mark=square*, mark size=4pt] coordinates {
					(40, 0.0842883528851173)
					(176, 0.0458649225114021)
					(736, 0.0245538108960542)
					(3008, 0.0130362142905765)
					(12160, 0.0067667627237866)
					(48896, 0.00297978262349399)
				};
				\addplot[line width=1.5pt, color=olive, mark=triangle, mark size=4pt] coordinates {
					(800, 0.100270018833246)
					(3200, 0.0501236640874951)
					(20000, 0.0208527149372604)
					(80000, 0.0110957318871821)
					(320000, 0.00610378707086293)
				};
			    \addplot[line width=1.5pt, color=violet, mark=triangle*, mark size=4pt] coordinates {								
				(1600, 0.2078237)
			    };
				\addplot[line width=1.5pt, color=violet, mark=triangle, mark size=4pt] coordinates {					
					(168.22/8, 0.0974178030904677)
					(369.46/8, 0.070611059426834)
					(646.86/8, 0.0590115949263689)
					(1004.46/8, 0.0520523700392549)
					(1439.80/8, 0.0453518186562592)
					(2543.35/8, 0.0371037643802047)
					(4781.20/8, 0.0278255542193279)
					(5686.25/8, 0.0255042998797395)
					(8869.04/8, 0.0226702366635657)
					(10089.3/8, 0.0218964852170362)
				};
				\legend{Specht,DKT,DKMT,DKQ,DKMQ,DKMQ24,DKMQ24\!+~,Makrodimopoulos,T6, T6b}
			\end{loglogaxis}
		\end{tikzpicture}
	}
	\caption{B01: Clamped pressure loaded square plate. Note, that the DKMQ24 and DKMQ24+ graphs coincide. The results for the T6b element are given in \cite{Bleyer2013} in terms of optimization unknowns $N$. The recalculation to the number of dofs $\Nd$ is done as follows: $N=\Nd+12N_e+6N_{\text{edges}}$ (\cite{Bleyer2013}), the number of edges $N_{\text{edges}}=\frac{3}{2}N_e$ ergo $N=\Nd+21N_e$. For the mesh used in \cite{Bleyer2013}, we get $\Nd=3N_e$, also finally we obtain $\Nd=\frac{N}{8}$. The results for the Makrodimopoulos element are given in \cite{Makrodimopoulos} in terms of element number $N_e$, for the mesh used in \cite{Makrodimopoulos} we obtain $\Nd=2N_e$. }
	\label{fig:B01}
\end{figure}

\begin{figure}
	\centering
	\scalebox{.9}
	{
		\begin{tikzpicture}
			\begin{loglogaxis}[
				xlabel=\text{time [s]},
				ylabel= $\left|\frac{\lambda-\lambda_{\text{ref}}}{\lambda_{\text{ref}}}\right|$,
				x tick label style={
					/pgf/number format/1000 sep={}  
				},
				y tick label style={
					/pgf/number format/1000 sep={},
					/pgf/number format/fixed,       
					/pgf/number format/precision=5  
				},
				scaled y ticks=false, 
				xmin=0.01,
				xmax=186.06,
				ymin=0.00297978262349399,
				ymax=0.137119920128883,
				legend style={at={(1.02,0.02)},anchor=south west}
				]
				\addplot[line width=1.5pt, color=red, mark=triangle*, mark size=4pt, solid] coordinates {
					(0.05, 0.137119920128883)
					(0.22, 0.0807743414036442)
					(0.53, 0.0463529078078555)
					(3.09, 0.0264232488484492)
					(25.69, 0.014984139229879)
					(162.13, 0.00853447845522001)
				};
				\addplot[line width=1.5pt, color=brown, mark=triangle*, mark size=3pt, solid] coordinates {
					(0.05, 0.086875541739466)
					(0.22, 0.0527735245399469)
					(0.59, 0.0307212906446417)
					(3.75, 0.0181958203807492)
					(27.53, 0.0112857434594178)
					(166.52, 0.00717126908851621)
				};
				\addplot[line width=1.5pt, color=orange, mark=triangle*, mark size=5pt, solid] coordinates {
					(0.05, 0.0868184520432938)
					(0.2, 0.0526607746590729)
					(0.63, 0.0305427378548253)
					(4.13, 0.0177923124049829)
					(28.34, 0.010420843638674)
					(186.06, 0.00612856526967849)
				};
				\addplot[line width=1.5pt, color=yellow, mark=diamond, mark size=7pt, solid] coordinates {
					(0.01, 0.0843418120759684)
					(0.06, 0.0459650790769441)
					(0.16, 0.0247313879875655)
					(1.05, 0.0134116539220803)
					(4.61, 0.0074913435138753)
					(22.47, 0.00446788137323869)
				};
				\addplot[line width=3pt, color=blue, mark=square, mark size=3pt, solid] coordinates {
					(0.05, 0.0842890109142067)
					(0.06, 0.0458649225114021)
					(0.27, 0.0245538108960542)
					(1.16, 0.0130362142905765)
					(4.83, 0.0067667627237866)
					(21.23, 0.00297978262349399)
				};
				\addplot[line width=1.5pt, color=green, mark=square*, mark size=4pt] coordinates {
					(0.05, 0.0842883528851173)
					(0.06, 0.0458649225114021)
					(0.27, 0.0245538108960542)
					(1.17, 0.0130362142905765)
					(4.84, 0.0067667627237866)
					(21.28, 0.00297978262349399)
				};
				\addplot[line width=1.5pt, color=col1, mark=square*, mark size=4pt] coordinates {
					(0.05, 0.0842883528851173)
					(0.06, 0.0458649225114021)
					(0.27, 0.0245538108960542)
					(1.17, 0.0130362142905765)
					(4.84, 0.0067667627237866)
					(21.28, 0.00297978262349399)
				};
		\addplot[line width=1.5pt, color=olive, mark=triangle, mark size=4pt] coordinates {
				(0.65*0.1, 0.100270018833246)
				(0.65*0.3, 0.0501236640874951)
				(0.65*2.3, 0.0208527149372604)
				(0.65*9.6, 0.0110957318871821)
				(0.65*53, 0.00610378707086293)
			};
				\legend{Specht,DKT,DKMT,DKQ,DKMQ,DKMQ24,DKMQ24\!+~,Makrodimopoulos}
			\end{loglogaxis}
		\end{tikzpicture}
	}
	\caption{B01: Clamped pressure loaded square plate. Note, that the DKMQ24 and DKMQ24+ graphs coincide. Time results for the Makrodimopoulos were calculated on the i7-2670QM, which is according to the SuperPI-32M test $0.65\times$ slower then i7-10875H, the time values were therefore decreased accordingly.}
	\label{fig:B01TIME}
\end{figure}
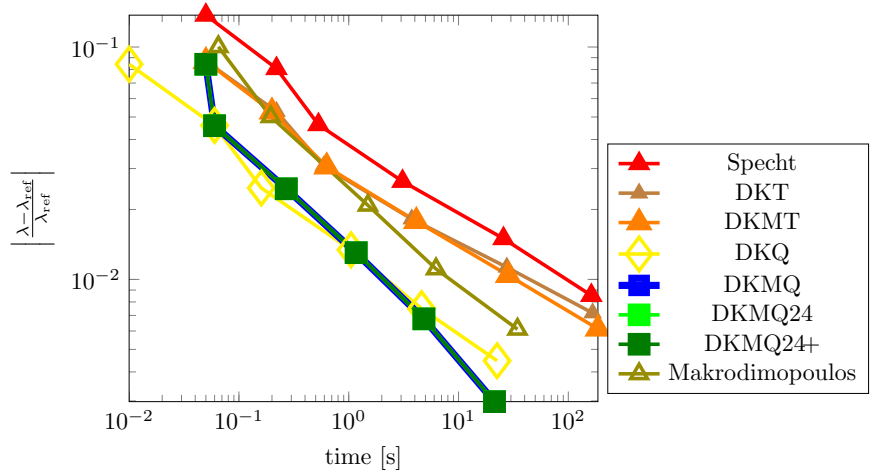

A fully clamped pressure loaded square plate is considered (see Fig. \ref{ClampedSquarePlate}). It is a standard benchmark for limit analysis (\cite{Bleyer2013, Makrodimopoulos, Leonetti}). Only one quarter is computed due to the symmetry. Symmetry and boundary conditions are as follows: $u_X=\varphi_Y=\varphi_Z=0$ at $X=0$, $u_Y=\varphi_X=\varphi_Z=0$ at $Y=0$ and $u_X=u_Y=u_Z=\varphi_X=\varphi_Y=\varphi_Z=0$ at $X=\frac{L}{2}$ or $Y=\frac{L}{2}$. The von Mises plasticity model is used. The analytical solution for this example is not available, therefore the reference solution is obtained numerically using the DKQ element on the refined mesh $256\!\times\!256$. Numerical results are summarized in Tables \ref{Tab:B01TRI}, \ref{Tab:B01QUAD} and \ref{Tab:B01HQUAD}. Let us discuss the results depicted in Fig. \ref{fig:B01}. 
None of the MITC elements converge. The lowest numerical error is achieved by all quadrangle elements (DKQ, DKMQ, DKMQ24, DKMQ24+) and by the T6b triangle element, which, however, seems to have a slightly lower convergence rate. The results of the strict upper-bound plate elements, i.e. T6 and the element by Macrodimopoulos, exhibit significantly higher numerical error than all other tested elements. When we compare two different implementations of the Kirchhoff plate elements, i.e., Specht and DKT, the Specht element exhibits the larger error. When we focus on the time-dependent results in Fig. \ref{fig:B01TIME}, we see that quadrangle elements run faster than the triangle elements. All the tested elements converge monotonically from above to the reference solution. 

We checked for the presence of shear locking in the case of the Mindlin elements of the DKM family (DKMT, DKMQ, DKMQ24, DKMQ24+). When the plate thickness $h$ is limited to zero (increasing at the same time plastic strengths, so that the resulting load factor $\lambda$ keeps the same order), the load factor converges to a fixed value, i.e., no presence of shear locking is observed, see Table \ref{Tab:ShearLocking}. However, we have encountered numerical problems in the case of moderate plate thickness, where the resulting plastic multiplier loses converge with respect to the number of degrees of freedom on finer meshes (see lines with marked with a star (*) in Table \ref{Tab:ShearLocking}) (in our case only for $L=0.1\,\text{m}$). This problem occurs only in this example and only for clamped boundary conditions. In the limit the plastic load multiplier goes to zero. The FE formulation of all these four shell elements is based on the following interpolation formula (\cite{DKMT,DKMQ,Stembera2020}):
\begin{align}
	\bm{B}_{\text{b}} = \bm{B}_{\text{b}\varphi}+\bm{B}_{\text{b}\triangle\beta}\bm{A}_n.
\end{align}
We found out that if the formula is adopted as follows:
\begin{align}
	\bm{B}_{\text{b}} = \bm{B}_{\text{b}\varphi}+K\bm{B}_{\text{b}\triangle\beta}\bm{A}_n,
\end{align}
where $K>1$ is to be understood as a numerical parameter, then convergence is recovered for sufficiently high values of $K$ for any plate thickness $h$. We can also set a higher $K$ for smaller thicknesses and finer meshes, i.e., for example, as $K\!=\!1\!+\!10\frac{h^2}{A}$, see Table \ref{Tab:ShearLocking}. However, the interpretation of this modification remains unclear. In conclusion, the DKMT, DKMQ, DKMQ24, DKMQ24+ elements always work well for the small shell thickness (in comparison to the shell characteristic length), but for moderate thickness, convergence with  respect to the number of degrees of freedom should always be checked. We also note that the tested Kirchhoff elements (DKT, Specht, DKQ) do not exhibit this problem.
 

\begin{center}
	\begin{table}
		\centering
		\begin{tabularx}{\textwidth}{sssslssl}
		\hline
		$h$       & $f_{\t}$   & \#dofs & \scalebox{0.75}[1.0]{DKMT} & \scalebox{0.75}[1.0]{DKMT} &  \#dofs & \multicolumn{2}{c}{\scalebox{0.75}[1.0]{DKMQ,DKMQ24,DKMQ24\!+~}}  \\
		&    &   &~ $K=1$  &  $K\!=\!1\!\!+\!\!10\frac{h^2}{A}$ & &  $K=1$    &  $K\!=\!1\!\!+\!\!10\frac{h^2}{A}$ \\
		$\text{[m]}$  &  [MPa]  &     & ~ $\lambda$ [-]  &  ~$\lambda$ [-] &  &  ~$\lambda$ [-] &   ~$\lambda$ [-] \\
		\hline \hline
		&                & $6080$   & $43.737$  & $45.450$ & $3008$   & $42.982$  &  $45.008$ \\
		$0.1$&$0.04$     & $24448$  & $42.098*$ & $44.793$ & $12160$  & $41.074*$ &  $44.600$ \\
		&                & $98048$  & $38.845*$ & $44.462$ & $48896$  & $37.821*$ &  $44.370$ \\
		\hline
		&                & $6080$   & $44.855$  & $45.007$  & $3008$  & $44.646$  & $44.679$ \\
		$0.01$ & $4$     & $24448$  & $44.529$  & $44.734$  & $12160$ & $44.369$  & $44.432$ \\
		&                & $98084$  & $44.339$  & $44.461$  & $48896$ & $44.202$  & $44.309$ \\
		\hline
		&                & $6080$   & $44.873$  & $44.874$  &$3008$   & $44.662$  & $44.662$ \\
		$0.001$&$400$    & $24448$  & $44.567$  & $44.571$  &$12160$  & $44.401$  & $44.401$ \\
		&                & $98084$  & $44.385$  & $44.392$  &$48896$  & $44.266$  & $44.269$ \\
		\hline
		&                & $6080$   & $44.873$  & $44.873$  &$3008$   & $44.662$  & $44.662$ \\
		$0.0001$&$40000$ & $24448$  & $44.567$  & $44.569$  &$12160$  & $44.401$  & $44.401$ \\
		&                & $98084$  & $44.385$  & $44.387$  &$48896$  & $44.266$  & $44.268$ \\
		\hline
		&                & $6080$   & $44.873$  & $44.873$  &$3008$   & $44.662$  & $44.662$ \\
		$0.00001$ &$4000000$&$24448$& $44.567$  & $44.569$  &$12160$  & $44.401$  & $44.401$ \\
		&                & $98084$  & $44.385$  & $44.388$  &$48896$  & $44.266$  & $44.268$ \\
		\hline				
	\end{tabularx}
	\captionof{table}{Clamped plate, dependence of the plate thickness $h$ going to $0$. $*$---a loss of convergence. \label{Tab:ShearLocking}}
\end{table}
\end{center}

\begin{center}
	\begin{table}
		\centering
		\begin{tabularx}{\tabulkywidth}{c *{6}{Y}}
			\hline
			\#elements&\#dofs& \scalebox{0.75}[1.0]{Specht}& \scalebox{0.75}[1.0]{DKT}& \scalebox{0.75}[1.0]{DKMT}& \scalebox{0.75}[1.0]{MITC3} \rule[-2mm]{0mm}{6mm} \\
			&& $\lambda$ [-] & $\lambda$ [-] & $\lambda$ [-] & $\lambda$ [-]  \\
			\hline \hline
			$4\!\!\times\!\!4\!\!\times\!\!4    $ & $88       $ & $50.114$ & $47.900$ & $47.897$ & $0.000$  \\
			$8\!\!\times\!\!8\!\!\times\!\!4    $ & $368      $ & $47.631$ & $46.397$ & $46.392$ & $0.000$  \\
			$16\!\!\times\!\!16\!\!\times\!\!4  $ & $1504     $ & $46.114$ & $45.425$ & $45.417$ & $0.000$  \\
			$32\!\!\times\!\!32\!\!\times\!\!4  $ & $6080     $ & $45.235$ & $44.873$ & $44.855$ & $0.000$  \\
			$64\!\!\times\!\!64\!\!\times\!\!4  $ & $24448    $ & $44.731$ & $44.568$ & $44.530$ & $0.000$  \\
			$128\!\!\times\!\!128\!\!\times\!\!4$ & $98048    $ & $44.447$ & $44.387$ & $44.341$ & $0.000$  \\
			\hline
			\multicolumn{6}{c}{Reference solution: $\lambda=44.071$ (DKMQ24 element at mesh $256\!\!\times\!\!256$)} \\
			\hline
		\end{tabularx}
		\captionof{table}{Clamped plate (triangle elements). The MITC3 element does not converge. \label{Tab:B01TRI}}
	\end{table}
\end{center}

\begin{center}
	\begin{table}
		\centering
		\begin{tabularx}{\tabulkywidth}{c *{7}{Y}}
			\hline
			\#elements&\#dofs& \scalebox{0.75}[1.0]{DKQ}& \scalebox{0.75}[1.0]{DKMQ}& \scalebox{0.75}[1.0]{DKMQ24}& \scalebox{0.75}[1.0]{DKMQ24\!+~}& \scalebox{0.75}[1.0]{MITC4} \rule[-2mm]{0mm}{6mm} \\
			&& $\lambda$ [-] & $\lambda$ [-] & $\lambda$ [-] & $\lambda$ [-] & $\lambda$ [-]  \\
			\hline \hline
			$4\!\!\times\!\!4                   $ & $40       $ & $47.788$ & $47.786$ & $47.786$ & $47.786$ & $0.000$  \\
			$8\!\!\times\!\!8                   $ & $176      $ & $46.097$ & $46.092$ & $46.092$ & $46.092$ & $0.000$  \\
			$16\!\!\times\!\!16                 $ & $736      $ & $45.161$ & $45.153$ & $45.153$ & $45.153$ & $0.000$  \\
			$32\!\!\times\!\!32                 $ & $3008     $ & $44.662$ & $44.646$ & $44.646$ & $44.646$ & $0.000$  \\
			$64\!\!\times\!\!64                 $ & $12160    $ & $44.401$ & $44.369$ & $44.369$ & $44.369$ & $0.000$  \\
			$128\!\!\times\!\!128               $ & $48896    $ & $44.268$ & $44.202$ & $44.202$ & $44.202$ & $0.000$  \\
			\hline
			\multicolumn{7}{c}{Reference solution: $\lambda=44.071$ (DKMQ24 element at mesh $256\!\!\times\!\!256$)} \\
			\hline
		\end{tabularx}
		\captionof{table}{Clamped plate (quad elements). The MITC4 element does not converge. \label{Tab:B01QUAD}}
	\end{table}
\end{center}

\begin{center}
	\begin{table}
		\centering
		\begin{tabularx}{110mm}{c *{5}{Y}}
			\hline
			\#elements&\#dofs&\scalebox{0.75}[1.0]{MITC9}&\#dofs&\scalebox{0.75}[1.0]{MITC16}  \rule[-2mm]{0mm}{6mm} \\
			&& $\lambda$ [-] && $\lambda$ [-]  \\
			\hline \hline
			$4\!\!\times\!\!4                   $ & $176      $ & $0.000$& $408      $ & $0.000$  \\
			$8\!\!\times\!\!8                   $ & $736      $ & $0.000$& $1680     $ & $0.000$  \\
			$16\!\!\times\!\!16                 $ & $3008     $ & $0.000$& $6816     $ & $0.000$  \\
			$32\!\!\times\!\!32                 $ & $12160    $ & $0.000$& $27456    $ & $0.000$  \\
			$64\!\!\times\!\!64                 $ & $48896    $ & $0.004$& $110208   $ & $0.000$  \\
			$128\!\!\times\!\!128               $ & $196096   $ & $0.001$& $441600   $ & $0.000$  \\
			\hline
			\multicolumn{5}{c}{Reference solution: $\lambda=44.071$ (DKMQ24 element at mesh $256\!\!\times\!\!256$)} \\
			\hline
		\end{tabularx}
		\captionof{table}{Clamped plate (high order quad elements). The MITC9 and MITC16 element do not converge. \label{Tab:B01HQUAD}}
	\end{table}
\end{center}

\begin{figure}[H]
	\centering
	\begin{overpic}[width=0.4\textwidth]{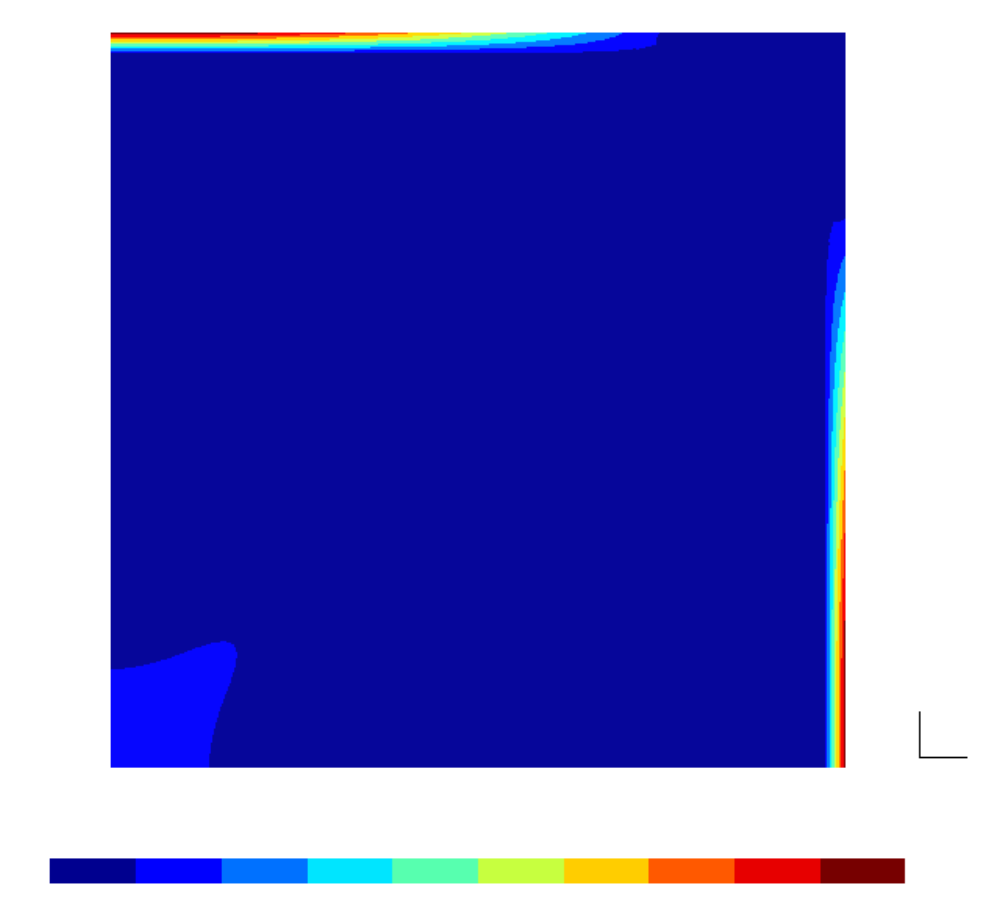} 
		\put(94,10.5){$X$}
		\put(87,18){$Y$}
		\put(3,-3){$0$}
		\put(44.5,-3){$19.2$}
		\put(87,-3){$38.4$}
		\put(44.5,8){$D_{\p}$}		
	\end{overpic}
	\caption{Clamped pressure loaded square plate -- dissipation power at collapse.}
	\label{ClampedDissipation}
\end{figure}

\begin{figure}[H]
	\centering
	\begin{overpic}[width=0.4\textwidth]{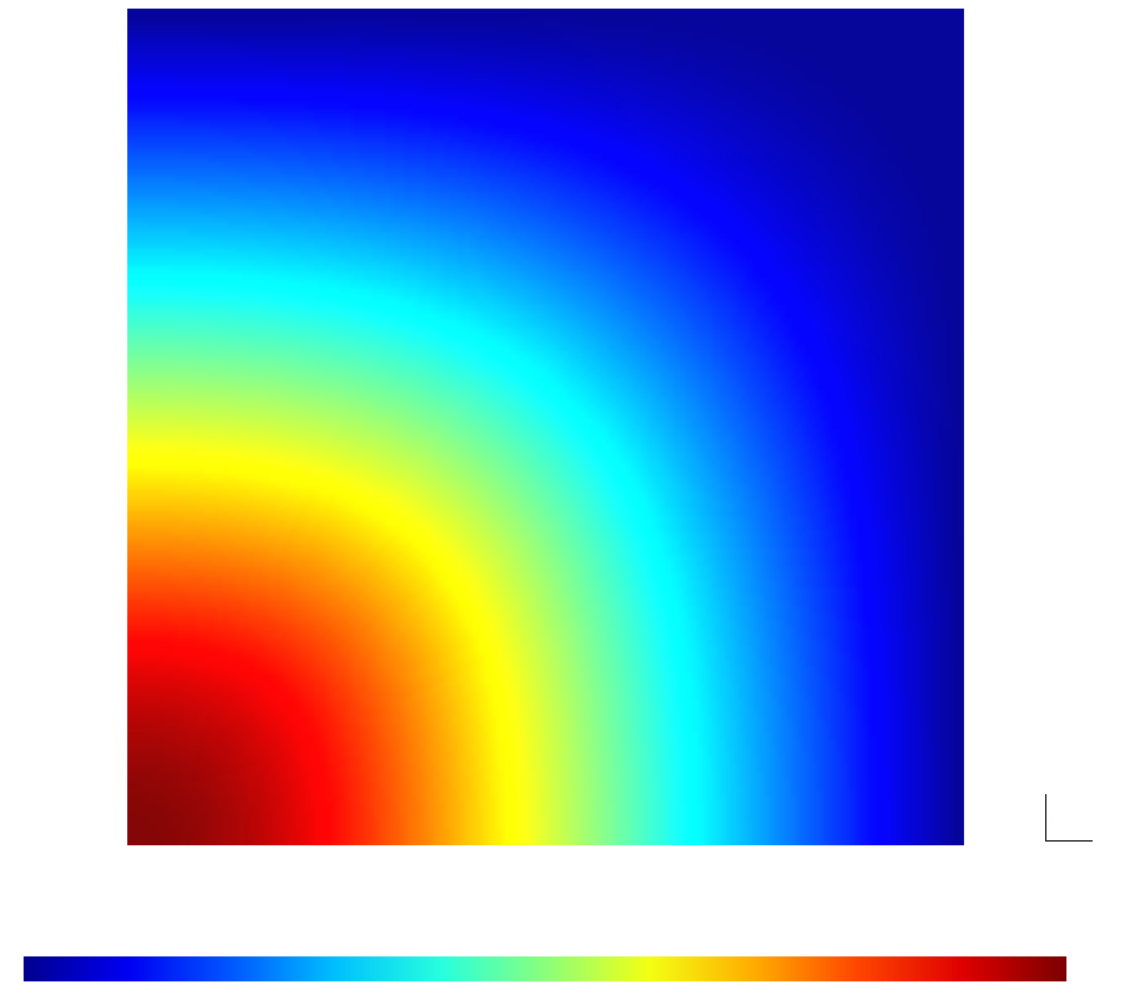} 
		\put(94,8){$X$}
		\put(87,15.5){$Y$}
		\put(3,-5){$0$}
		\put(42.5,-5){$0.0605$}
		\put(87,-5){$0.121$}
		\put(46,5.5){$|\dot{\bm{u}}|$}
	\end{overpic}
	\caption{Clamped pressure loaded square plate -- absolute value of velocity at collapse.}
	\label{ClampedVelocity}
\end{figure}

\FloatBarrier

\subsection{Cook's membrane in plane stress}
\label{b2}

	\begin{figure}[H]
	\centering
	\begin{overpic}[width=0.33\textwidth]{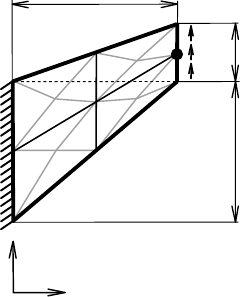} 
		\put(29,102.7){$L_1$}
		\put(81.5,49){$L_2$}
		\put(81.5,81){$L_3$}
		\put(68,81){$F_Z$}
		\put(53,83){$A$}
		\put(23,-2){$X$}
		\put(-2,14){$Z$}
		\put(35,42){thickness $h$}
	\end{overpic}
	\caption{Cook's membrane with mesh $2\times 2$. Parameters: $L_1=4.8\,\text{m},\,L_2=4.4\,\text{m},\,L_3=1.6\,\text{m},\,h=1\,\text{m}$, yield strength $\sigma_0=\sqrt{3}\,\text{Pa}$, loading force $F_Z=1\,\text{N}$, test point $A$.}
	\label{Cook}
\end{figure}
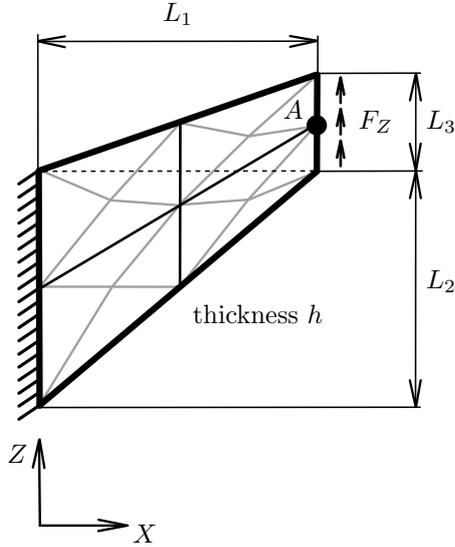
\FloatBarrier

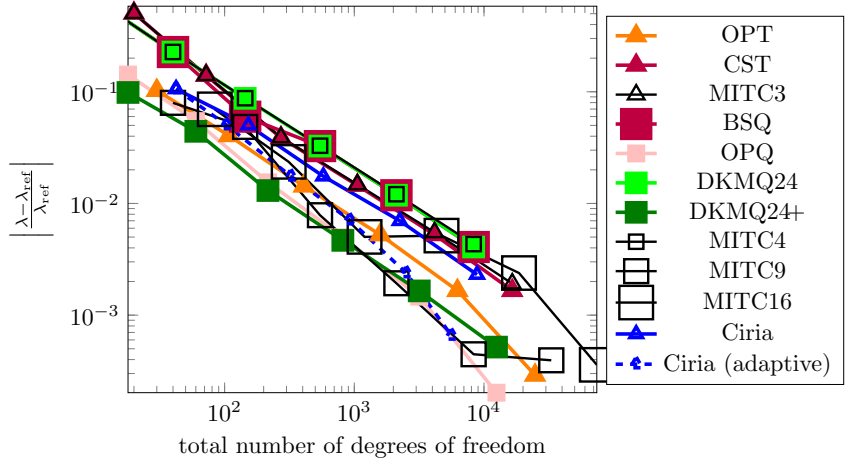
\begin{figure}
	\centering
	\scalebox{.9}
	{
		\begin{tikzpicture}
			\begin{loglogaxis}[
				xlabel=\text{total number of degrees of freedom},
				ylabel= $\left|\frac{\lambda-\lambda_{\text{ref}}}{\lambda_{\text{ref}}}\right|$,
				x tick label style={
					/pgf/number format/1000 sep={}  
				},
				y tick label style={
					/pgf/number format/1000 sep={},
					/pgf/number format/fixed,       
					/pgf/number format/precision=5  
				},
				scaled y ticks=false, 
				xmin=18,
				xmax=74112,
				ymin=0.000202123963564139,
				ymax=0.583719906866752,
				legend style={at={(1.02,0.02)},anchor=south west}
				]
				\addplot[line width=1.5pt, color=orange, mark=triangle*, mark size=4pt] coordinates {
					(30, 0.10279330404064)
					(108, 0.0404708411187668)
					(408, 0.0143096819163844)
					(1584, 0.00524112095644053)
					(6240, 0.00166057310522015)
					(24768, 0.000289788917435573)
				};
				\addplot[line width=1.5pt, color=purple, mark=triangle*, mark size=4pt] coordinates {
					(20, 0.504817273443886)
					(72, 0.14026064463389)
					(272, 0.038737770057223)
					(1056, 0.0145128014422515)
					(4160, 0.00525259692864651)
					(16512, 0.00166044610533696)
				};
				\addplot[line width=1pt, color=black, mark=triangle, mark size=4pt] coordinates {
					(20, 0.505040275020437)
					(72, 0.140521335250496)
					(272, 0.038999200776597)
					(1056, 0.0147690893962397)
					(4160, 0.00550681128109318)
					(16512, 0.00191640853088876)
				};
				\addplot[line width=1.5pt, color=purple, mark=square*, mark size=6pt] coordinates {
					(12, 0.583568229592433)
					(40, 0.226911332184982)
					(144, 0.0598309887013897)
					(544, 0.0326275691930398)
					(2112, 0.0118453469870373)
					(8320, 0.00407367452995448)
				};
				\addplot[line width=1.5pt, color=pink, mark=square*, mark size=3pt] coordinates {
					(18, 0.142555750029195)
					(60, 0.0546500729884386)
					(216, 0.0154538355424501)
					(816, 0.00498252218848537)
					(3168, 0.0014637466425318)
					(12480, 0.000202123963564139)
				};
				\addplot[line width=1.5pt, color=green, mark=square*, mark size=4pt] coordinates {
					(12, 0.583568229592433)
					(40, 0.226911332184982)
					(144, 0.0873140823893496)
					(544, 0.0326276779458134)
					(2112, 0.011845352096228)
					(8320, 0.00407377233446223)
				};
				\addplot[line width=1.5pt, color=col1, mark=square*, mark size=4pt] coordinates {
					(18, 0.0985445506831718)
					(60, 0.0442150721125774)
					(216, 0.0130438164194791)
					(816, 0.00466931347074624)
					(3168, 0.00164811689828328)
					(12480, 0.000517480731052297)
				};
				\addplot[line width=1pt, color=black, mark=square, mark size=3pt] coordinates {
					(12, 0.583719906866752)
					(40, 0.227148154122387)
					(144, 0.0875662187609483)
					(544, 0.0328800252539999)
					(2112, 0.0120976096286349)
					(8320, 0.00432636707345554)
				};
				\addplot[line width=1pt, color=black, mark=square, mark size=5pt] coordinates {
					(40, 0.0795893998890576)
					(144, 0.0484488905757328)
					(544, 0.00784931025925491)
					(2112, 0.00192951360504494)
					(8320, 0.000446371744715663)
					(33024, 0.000394137)
				};
				\addplot[line width=1pt, color=black, mark=square, mark size=7pt] coordinates {
					(84, 0.0696796712600724)
					(312, 0.0235676303573514)
					(1200, 0.00501729023122727)
					(4704, 0.00515256335396468)
					(18624, 0.00238231709097273)
					(74112, 0.000359042537662)
				};
				\addplot[line width=1.5pt, color=blue, mark=triangle, mark size=3pt] coordinates {
					(42.24621125123532, 0.10590622445404649)
					(152.49242250247065, 0.05009926427653862)
					(576.9848450049412, 0.017517225271517008)
					(2241.9696900098825, 0.006992292420880637)
					(8835.939380019765, 0.002306434660749775)					
				};
				\addplot[line width=1.5pt, color=blue, mark=triangle, mark size=3pt, dashed] coordinates {
					(42.24621125123532, 0.10590622445404649)
					(103.41640786499875, 0.05031822959243256)
					(324.4948974278318, 0.017517225271517008)
					(924, 0.007079878547238148)
					(2520, 0.002379423099381088)
					(5610.938076978759, 0.000656895947681977)
				};
				\legend{OPT,CST,MITC3,BSQ,OPQ,DKMQ24,DKMQ24\!+~,MITC4, MITC9,MITC16,Ciria,Ciria (adaptive)}
			\end{loglogaxis}
		\end{tikzpicture}
	}
	\caption{B02: Cook's membrane. The results for the Ciria element are given in \cite{Ciria} in terms of the number of elements $N_e$. The number of dofs $\Nd$ is approximated in this case by the formula valid for the uniform mesh $\Nd=N_e+\sqrt{2N_e}$. }	
	\label{fig:B02}
\end{figure}

The standard Cook's membrane problem tests an element's plastic behaviour for the case of in-plane shear loading. The problem also tests the effect of mesh distortion. The membrane is fully fixed at $X=0$ and loaded by a distributed force at $X=L_1$ (see Fig. \ref{Cook}). The $Z$-component of the displacement is evaluated at the test point A. The approximate reference solution for this example is taken from \cite{Ciria}: $\lambda=0.68504$. We note that in this in-plane test problem, we eliminate drilling rotations $\varphi_Z$ in all elements that do not use the drilling rotations to increase their numerical efficiency (CST, MITC3, BSQ, DKMQ24, MITC3, MITC9, MICT16). The numerical results are summarized in Tables \ref{Tab:B02TRI}, \ref{Tab:B02QUAD}, and \ref{Tab:B02HQUAD}. 
Let us discuss the results depicted in Fig. \ref{fig:B02}. We again observe that most efficient elements are quads (DKMQ24+, MITC9, and OPQ in particular), which even, in the case of OPQ element, overcome results by the Ciria element obtained with the help of adaptive meshing. We also observe that in this pure membrane case, the optimal triangle OPT and quadrangle OPQ by Felippa overcome constant strain triangle CST and bilinear strain quadrangle BSQ. The highest efficiency is exhibited by the OPQ quadrangle element. The DKMQ24+ element performs better than the DKMQ24 element in this example. All tested elements, except for MITC16 which converges from below, converge to the reference solution monotonically from above.

 

\begin{center}
	\begin{table}
		\centering
		\begin{tabularx}{\tabulkywidth}{c *{6}{Y}}
			\hline
			\#elements&\#dofs& \scalebox{0.75}[1.0]{OPT}&\#dofs& \scalebox{0.75}[1.0]{CST}& \scalebox{0.75}[1.0]{MITC3} \rule[-2mm]{0mm}{6mm} \\
			&& $\lambda$ [-] && $\lambda$ [-] & $\lambda$ [-]  \\
			\hline \hline
			$2\!\!\times\!\!2\!\!\times\!\!4    $ & $30   $ & $0.75546$ & $20   $ & $1.03086$ & $1.03101$  \\
			$4\!\!\times\!\!4\!\!\times\!\!4    $ & $108  $ & $0.71276$ & $72   $ & $0.78112$ & $0.78130$  \\
			$8\!\!\times\!\!8\!\!\times\!\!4    $ & $408  $ & $0.69484$ & $272  $ & $0.71158$ & $0.71176$  \\
			$16\!\!\times\!\!16\!\!\times\!\!4  $ & $1584 $ & $0.68863$ & $1056 $ & $0.69498$ & $0.69516$  \\
			$32\!\!\times\!\!32\!\!\times\!\!4  $ & $6240 $ & $0.68618$ & $4160 $ & $0.68864$ & $0.68881$  \\
			$64\!\!\times\!\!64\!\!\times\!\!4  $ & $24768$ & $0.68524$ & $16512$ & $0.68618$ & $0.68635$  \\
			\hline
			\multicolumn{6}{c}{Reference solution: $\lambda=0.68504$ (\cite{Ciria})} \\
			\hline
		\end{tabularx}
		\captionof{table}{Cook's membrane (triangle elements).  \label{Tab:B02TRI}}
	\end{table}
\end{center}

\begin{center}
	\begin{table}
		\centering
		\begin{tabularx}{\tabulkywidth}{c *{8}{Y}}
			\hline
			\#elements&\#dofs& 
			\scalebox{0.75}[1.0]{OPQ}& 
			\scalebox{0.75}[1.0]{DKMQ24\!+~}& 
			\#dofs& 
			\scalebox{0.75}[1.0]{BSQ}& 
			\scalebox{0.75}[1.0]{DKMQ24}& 
			\scalebox{0.75}[1.0]{MITC4} \rule[-2mm]{0mm}{6mm} \\
			&& $\lambda$ [-] & $\lambda$ [-] && $\lambda$ [-] & $\lambda$ [-] & $\lambda$ [-]  \\
			\hline \hline
			$2\!\!\times\!\!2    $ & $18   $ & $0.78270$ & $0.75255$ & $12   $ & $1.08481$ & $1.08481$ & $1.08491$  \\
			$4\!\!\times\!\!4    $ & $60   $ & $0.72248$ & $0.71533$ & $40   $ & $0.84048$ & $0.84048$ & $0.84065$  \\
			$8\!\!\times\!\!8    $ & $216  $ & $0.69563$ & $0.69398$ & $144  $ & $0.64405$ & $0.74485$ & $0.74503$  \\
			$16\!\!\times\!\!16  $ & $816  $ & $0.68845$ & $0.68824$ & $544  $ & $0.70739$ & $0.70739$ & $0.70756$  \\
			$32\!\!\times\!\!32  $ & $3168 $ & $0.68604$ & $0.68617$ & $2112 $ & $0.69315$ & $0.69315$ & $0.69333$  \\
			$64\!\!\times\!\!64  $ & $12480$ & $0.68518$ & $0.68539$ & $8320 $ & $0.68783$ & $0.68783$ & $0.68800$  \\
			\hline
			\multicolumn{8}{c}{Reference solution: $\lambda=0.68504$ (\cite{Ciria})} \\
			\hline
		\end{tabularx}
		\captionof{table}{Cook's membrane (quad elements).  \label{Tab:B02QUAD}}
	\end{table}
\end{center}

\begin{center}
	\begin{table}
		\centering
		\begin{tabularx}{110mm}{c *{5}{Y}}
			\hline
			\#elements&\#dofs&\scalebox{0.75}[1.0]{MITC9}&\#dofs&\scalebox{0.75}[1.0]{MITC16}  \rule[-2mm]{0mm}{6mm} \\
			&& $\lambda$ [-] && $\lambda$ [-]  \\
			\hline \hline
			$2\!\!\times\!\!2    $ & $40     $ & $0.73956$& $84      $ & $0.63731$  \\
			$4\!\!\times\!\!4    $ & $144    $ & $0.71823$& $312     $ & $0.66890$  \\
			$8\!\!\times\!\!8    $ & $544    $ & $0.69042$& $1200    $ & $0.68160$  \\
			$16\!\!\times\!\!16  $ & $2112   $ & $0.68636$& $4704    $ & $0.68151$  \\
			$32\!\!\times\!\!32  $ & $8320   $ & $0.68535$& $18624   $ & $0.68341$  \\
			$64\!\!\times\!\!64  $ & $33024  $ & $0.68531$& $74112   $ & $0.68479$  \\
			\hline
			\multicolumn{5}{c}{Reference solution: $\lambda=0.68504$ (\cite{Ciria})} \\
			\hline
		\end{tabularx}
		\captionof{table}{Cook's membrane (high order quad elements).  \label{Tab:B02HQUAD}}
	\end{table}
\end{center}

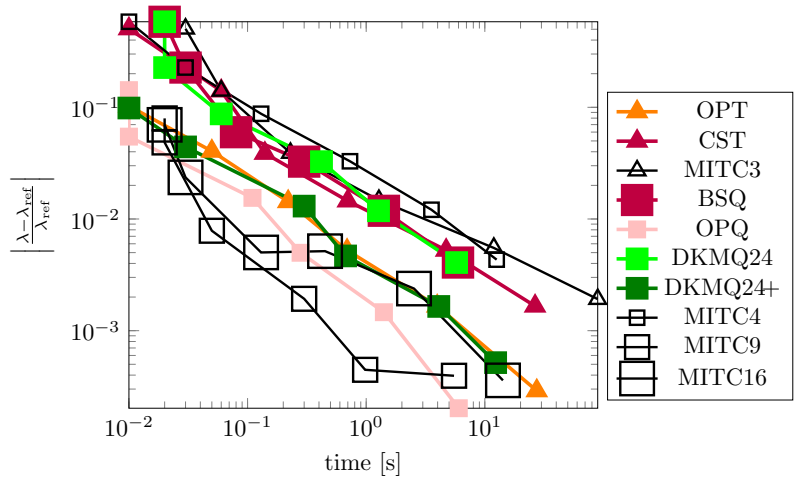
\begin{figure}
	\centering
	\scalebox{.9}
	{
		\begin{tikzpicture}
			\begin{loglogaxis}[
				xlabel=\text{time [s]},
				ylabel= $\left|\frac{\lambda-\lambda_{\text{ref}}}{\lambda_{\text{ref}}}\right|$,
				x tick label style={
					/pgf/number format/1000 sep={}  
				},
				y tick label style={
					/pgf/number format/1000 sep={},
					/pgf/number format/fixed,       
					/pgf/number format/precision=5  
				},
				scaled y ticks=false, 
				xmin=0.01,
				xmax=89.44,
				ymin=0.000202123963564139,
				ymax=0.583719906866752,
				legend style={at={(1.02,0.02)},anchor=south west}
				]
				\addplot[line width=1.5pt, color=orange, mark=triangle*, mark size=4pt] coordinates {
					(0.01, 0.10279330404064)
					(0.05, 0.0404708411187668)
					(0.22, 0.0143096819163844)
					(0.69, 0.00524112095644053)
					(3.97, 0.00166057310522015)
					(27.41, 0.000289788917435573)
				};
				\addplot[line width=1.5pt, color=purple, mark=triangle*, mark size=4pt] coordinates {
					(0.01, 0.504817273443886)
					(0.06, 0.14026064463389)
					(0.14, 0.038737770057223)
					(0.7, 0.0145128014422515)
					(4.72, 0.00525259692864651)
					(26.52, 0.00166044610533696)
				};
				\addplot[line width=1pt, color=black, mark=triangle, mark size=4pt] coordinates {
					(0.03, 0.505040275020437)
					(0.06, 0.140521335250496)
					(0.23, 0.038999200776597)
					(1.28, 0.0147690893962397)
					(11.91, 0.00550681128109318)
					(89.44, 0.00191640853088876)
				};
				\addplot[line width=1.5pt, color=purple, mark=square*, mark size=6pt] coordinates {
					(0.02, 0.583568229592433)
					(0.03, 0.226911332184982)
					(0.08, 0.0598309887013897)
					(0.3, 0.0326275691930398)
					(1.41, 0.0118453469870373)
					(5.94, 0.00407367452995448)
				};
				\addplot[line width=1.5pt, color=pink, mark=square*, mark size=3pt] coordinates {
					(0.01, 0.142555750029195)
					(0.01, 0.0546500729884386)
					(0.11, 0.0154538355424501)
					(0.28, 0.00498252218848537)
					(1.41, 0.0014637466425318)
					(6.03, 0.000202123963564139)
				};
				\addplot[line width=1.5pt, color=green, mark=square*, mark size=4pt] coordinates {
					(0.02, 0.583568229592433)
					(0.02, 0.226911332184982)
					(0.06, 0.0873140823893496)
					(0.42, 0.0326276779458134)
					(1.27, 0.011845352096228)
					(5.78, 0.00407377233446223)
				};
				\addplot[line width=1.5pt, color=col1, mark=square*, mark size=4pt] coordinates {
					(0.01, 0.0985445506831718)
					(0.03, 0.0442150721125774)
					(0.3, 0.0130438164194791)
					(0.67, 0.00466931347074624)
					(4.09, 0.00164811689828328)
					(12.33, 0.000517480731052297)
				};
				\addplot[line width=1pt, color=black, mark=square, mark size=3pt] coordinates {
					(0.01, 0.583719906866752)
					(0.03, 0.227148154122387)
					(0.13, 0.0875662187609483)
					(0.73, 0.0328800252539999)
					(3.55, 0.0120976096286349)
					(12.58, 0.00432636707345554)
				};
				\addplot[line width=1pt, color=black, mark=square, mark size=5pt] coordinates {
					(0.02, 0.0795893998890576)
					(0.02, 0.0484488905757328)
					(0.05, 0.00784931025925491)
					(0.3, 0.00192951360504494)
					(0.98, 0.000446371744715663)
					(5.5, 0.000394878401261265)
				};
				\addplot[line width=1pt, color=black, mark=square, mark size=7pt] coordinates {
					(0.02, 0.0696796712600724)
					(0.03, 0.0235676303573514)
					(0.13, 0.00501729023122727)
					(0.45, 0.00515256335396468)
					(2.53, 0.00238231709097273)
					(14.22, 0.000359042537662)
				};
				\legend{OPT,CST,MITC3,BSQ,OPQ,DKMQ24,DKMQ24\!+~,MITC4,MITC9,MITC16}
			\end{loglogaxis}
		\end{tikzpicture}
	}
	\caption{B02: Cook's membrane }
	\label{fig:B02TIME}
\end{figure}

\begin{figure}[H]
	\centering
	\begin{overpic}[width=0.4\textwidth]{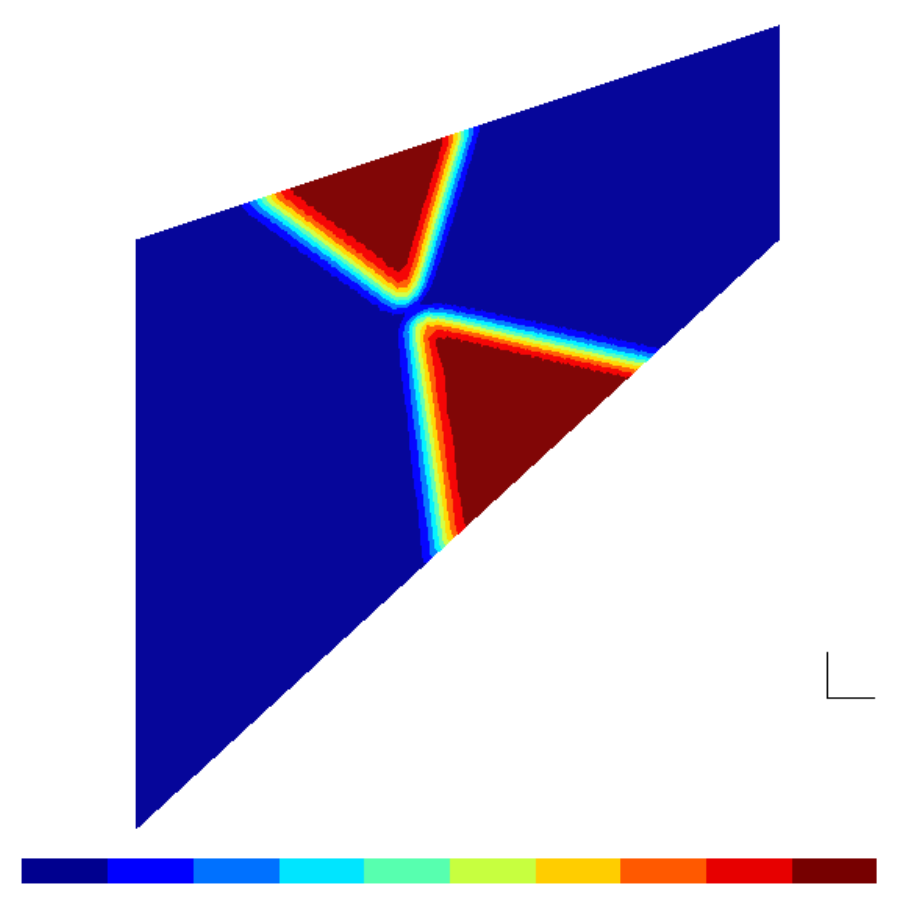} 
		\put(94,18){$X$}
		\put(85.3,26){$Y$}
		\put(1,-3){$0$}
		\put(44.5,-3){$0.149$}
		\put(88,-3){$0.297$}
		\put(45.5,8.5){$D_{\p}$}
	\end{overpic}
	\caption{Cook's membrane -- dissipation power at collapse.}
	\label{CookFELADissipation}
\end{figure}

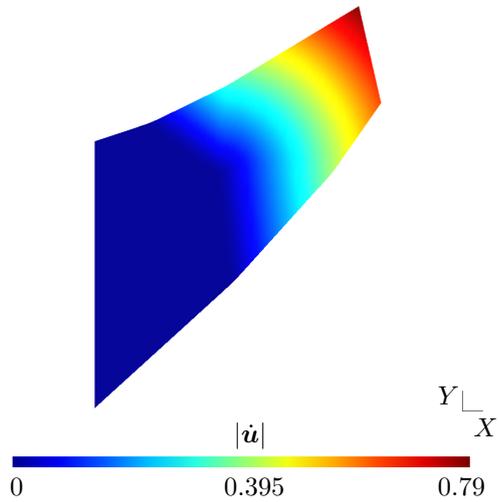
\begin{figure}[H]
	\centering
	\begin{overpic}[width=0.4\textwidth]{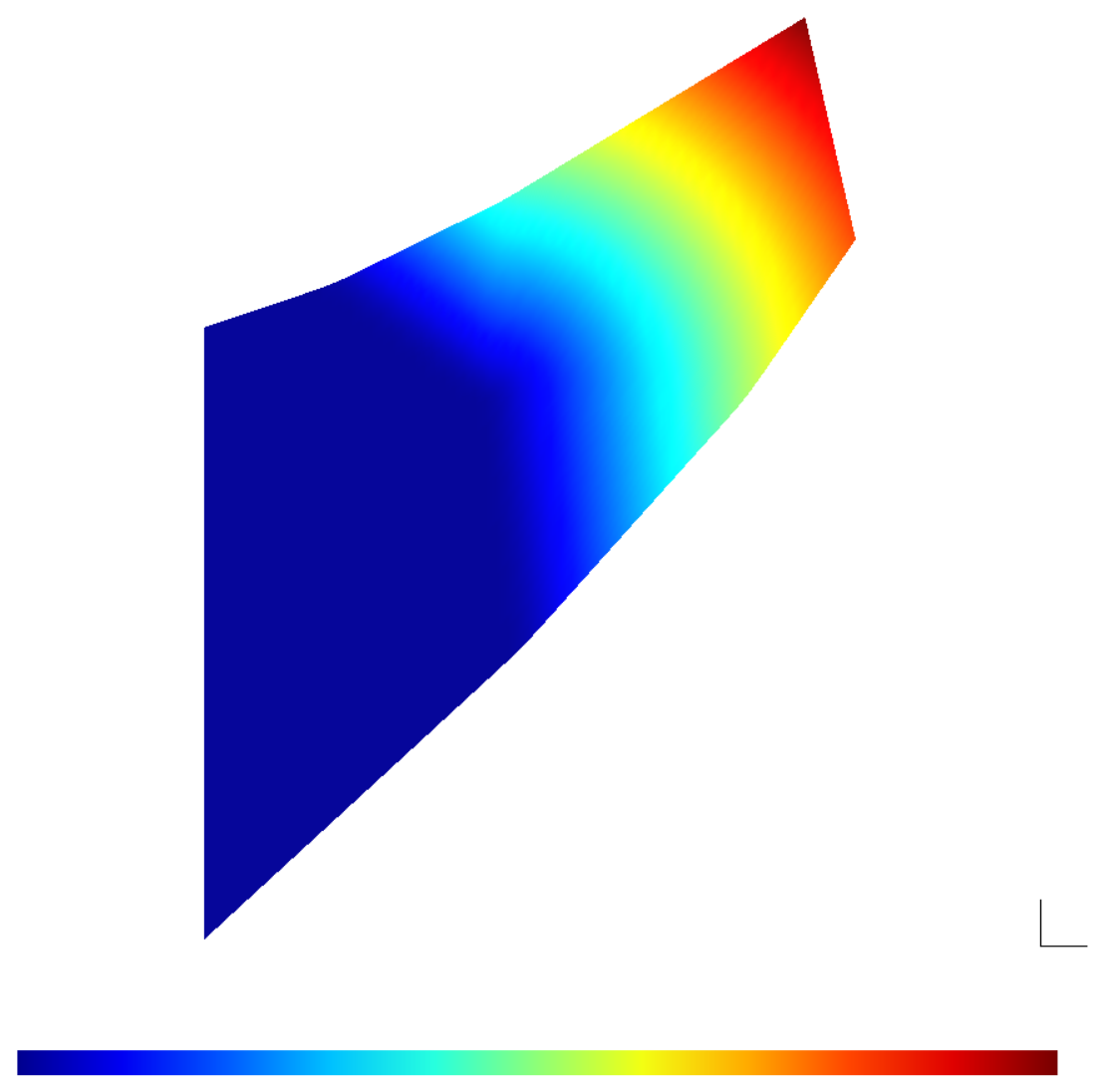} 
		\put(95.3,7.5){$X$}
		\put(88,14){$Y$}
		\put(1,-4.5){$0$}
		\put(44.5,-4.5){$0.395$}
		\put(88,-4.5){$0.79$}
		\put(46.5,6.3){$|\dot{\bm{u}}|$}
	\end{overpic}
	\caption{Cook's membrane -- absolute value of velocity at collapse.}
	\label{CookFELAVelocity}
\end{figure}

\FloatBarrier

\subsection{Simply supported spherical cap under uniform pressure}
\label{b3}

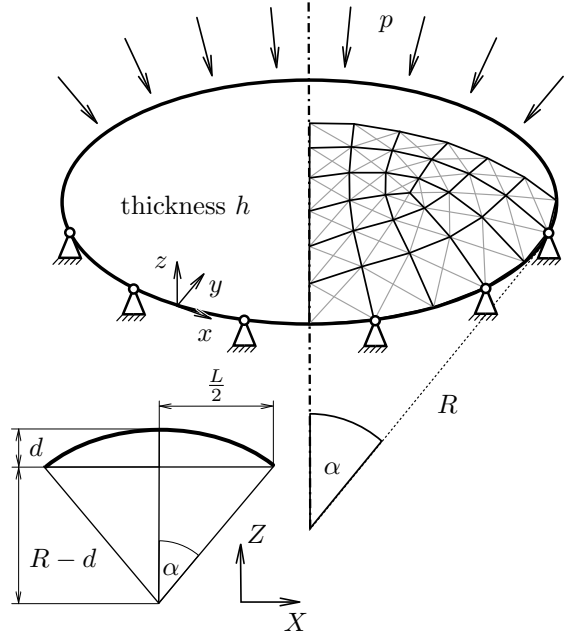
\begin{figure}[H]
	\centering
	\begin{overpic}[width=0.45\textwidth]{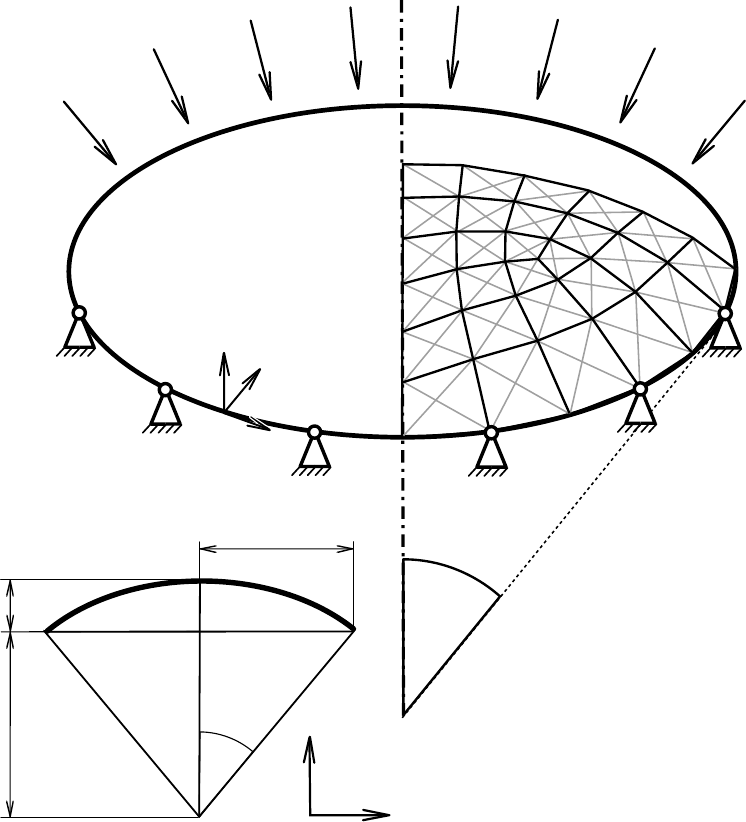} 
		\put(60.5,96){$p$}
		\put(18,65){thickness $h$}
		\put(70,32){$R$}
		\put(24.9,5.2){$\alpha$}
		\put(51.4,22){$\alpha$}
		\put(32,35.5){$\frac{L}{2}$}
		\put(3,25){$d$}
		\put(3,6.5){$R-d$}
		\put(44.8,-4){$X$}
		\put(38.7,10.5){$Z$}
		\put(30.6,44){$x$}
		\put(32.5,52){$y$}
		\put(23.5,56){$z$}
	\end{overpic}
	\caption{Simply supported spherical cap under uniform pressure with quad mesh $3\times3\times 3$, a triangular mesh depicted in a grey colour. Parameters: $R=1.25\,\text{m},\,h=0.001\,\text{m},\,L=2\,\text{m},\,d=R-\sqrt{R^2-\left(\frac{L}{2}\right)^2}=0.5\,\text{m},\,\alpha=\arcsin\frac{L}{2R}\approx53.1^{\circ}$, yield strength $\sigma_0=240\,\text{MPa}$, loading pressure $p=1\,\text{Pa}$.}
	\label{SphericalCap}
\end{figure}

A simply supported spherical cap of thickness $h$ loaded by uniform pressure $p$ is considered (see Fig. \ref{SphericalCap}). The Tresca failure criterion is considered and the 3-point Gauss-Lobatto quadrature is used for the integration across the shell thickness. This benchmark is taken from \cite{Bleyer2016}. Only one quarter is computed due to symmetry. Boundary and symmetry conditions are as follows: $u_X=\varphi_Y=\varphi_Z=0$ at $X=0$, $u_Y=\varphi_X=\varphi_Z=0$ at $Y=0$, $u_X=u_Y=u_Z=0$ at $Z=0$. Moreover, $\varphi_y=\varphi_z=0$ is set in the local coordinate system at $Z=0$. The analytical solution $p=\frac{2h\sigma_0}{R}=384.0\,\text{kPa}$ is taken from \cite{Bleyer2021}. The numerical results are summarized in Tables \ref{Tab:B03TRI}, \ref{Tab:B03QUAD}, and \ref{Tab:B03HQUAD}.

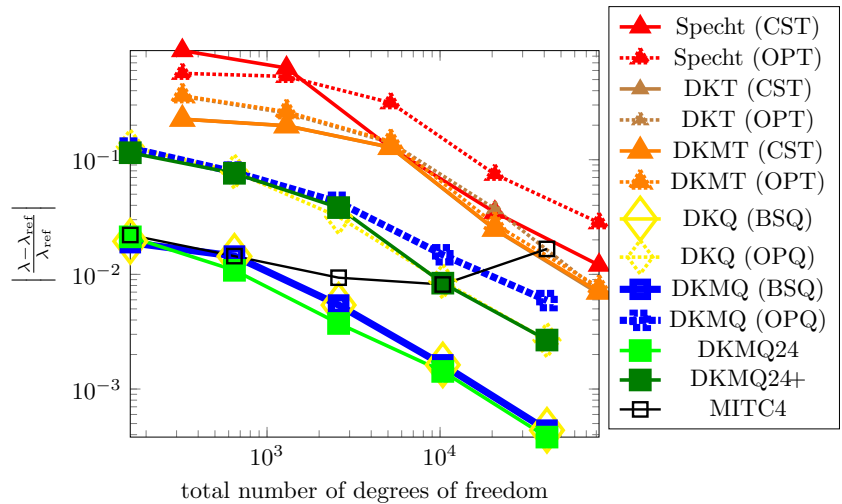
\begin{figure}
	\centering
	\scalebox{.9}
	{
		\begin{tikzpicture}
			\begin{loglogaxis}[
				xlabel=\text{total number of degrees of freedom},
				ylabel= $\left|\frac{\lambda-\lambda_{\text{ref}}}{\lambda_{\text{ref}}}\right|$,
				x tick label style={
					/pgf/number format/1000 sep={}  
				},
				y tick label style={
					/pgf/number format/1000 sep={},
					/pgf/number format/fixed,       
					/pgf/number format/precision=5  
				},
				scaled y ticks=false, 
				xmin=162,
				xmax=82944,
				ymin=0.00038095572916673,
				ymax=0.894259770572917,
				legend style={at={(1.02,0.02)},anchor=south west}
				]
				\addplot[line width=1.5pt, color=red, mark=triangle*, mark size=4pt, solid] coordinates {
					(324, 0.894259770572917)
					(1296, 0.63071941796875)
					(5184, 0.12690225)
					(20736, 0.0347791458333334)
					(82944, 0.0120278802083334)
				};
				\addplot[line width=1.5pt, color=red, mark=triangle*, mark size=4pt, densely dotted] coordinates {
					(324, 0.56660873828125)
					(1296, 0.532642317708333)
					(5184, 0.31376878125)
					(20736, 0.074599828125)
					(82944, 0.0281870091145833)
				};
				\addplot[line width=1.5pt, color=brown, mark=triangle*, mark size=3pt, solid] coordinates {
					(324, 0.225929434895833)
					(1296, 0.197799376302083)
					(5184, 0.128098509114583)
					(20736, 0.0245691666666667)
					(82944, 0.00678264192708335)
				};
				\addplot[line width=1.5pt, color=brown, mark=triangle*, mark size=3pt, densely dotted] coordinates {
					(324, 0.35659196875)
					(1296, 0.255566009114583)
					(5184, 0.140996455729167)
					(20736, 0.0373223385416666)
					(82944, 0.00716183072916669)
				};
				\addplot[line width=1.5pt, color=orange, mark=triangle*, mark size=5pt, solid] coordinates {
					(324, 0.225935109375)
					(1296, 0.197830467447917)
					(5184, 0.128295674479167)
					(20736, 0.0248045625)
					(82944, 0.00701608984375007)
				};
				\addplot[line width=1.5pt, color=orange, mark=triangle*, mark size=5pt, densely dotted] coordinates {
					(324, 0.3619703125)
					(1296, 0.262525053385417)
					(5184, 0.14424908984375)
					(20736, 0.0274023020833334)
					(82944, 0.00771385026041672)
				};
				\addplot[line width=1.5pt, color=yellow, mark=diamond, mark size=9pt, solid] coordinates {
					(162, 0.01943287890625)
					(648, 0.0143732656250001)
					(2592, 0.00539754036458332)
					(10368, 0.00162344921875005)
					(41472, 0.000436917968750066)
				};
				\addplot[line width=1.5pt, color=yellow, mark=diamond, mark size=7pt, densely dotted] coordinates {
					(162, 0.129709071614583)
					(648, 0.0787331744791667)
					(2592, 0.0317650260416666)
					(10368, 0.00868300130208337)
					(41472, 0.00266751171874999)
				};
				\addplot[line width=3pt, color=blue, mark=square, mark size=3pt, solid] coordinates {
					(162, 0.0189023580729167)
					(648, 0.0144517591145834)
					(2592, 0.00540477864583331)
					(10368, 0.00162675260416673)
					(41472, 0.000438197916666701)
				};
				\addplot[line width=3pt, color=blue, mark=square, mark size=3pt, densely dotted] coordinates {
					(162, 0.126375912760417)
					(648, 0.0780199986979166)
					(2592, 0.0425289049479167)
					(10368, 0.01481973828125)
					(41472, 0.0059070638020834)
				};
				\addplot[line width=1.5pt, color=green, mark=square*, mark size=4pt] coordinates {
					(162, 0.0212477552083333)
					(648, 0.0108500208333333)
					(2592, 0.00373063671874999)
					(10368, 0.00142447135416671)
					(41472, 0.00038095572916673)
				};
				\addplot[line width=1.5pt, color=col1, mark=square*, mark size=4pt] coordinates {
					(162, 0.11506812109375)
					(648, 0.0761313736979167)
					(2592, 0.0382261744791667)
					(10368, 0.00838111718750002)
					(41472, 0.00267416406250006)
				};
				\addplot[line width=1pt, color=black, mark=square, mark size=3pt] coordinates {
					(162, 0.0220456054687501)
					(648, 0.0144302421875001)
					(2592, 0.00935649739583341)
					(10368, 0.00816069140625)
					(41472, 0.0166721302083334)
				};
				\legend{Specht (CST),Specht (OPT),DKT (CST),DKT (OPT),DKMT (CST),DKMT (OPT),DKQ (BSQ),DKQ (OPQ),DKMQ (BSQ),DKMQ (OPQ),DKMQ24,DKMQ24\!+~,MITC4}
			\end{loglogaxis}
		\end{tikzpicture}
	}
	\caption{B03: Spherical cap }
	\label{fig:B03}
\end{figure}

Let us discuss the results depicted in Fig. \ref{fig:B03}. None of the MITC elements converge. All remaining quad elements are clearly superior to all triangle elements. The most efficient elements are DKMQ24, DKMQ (BSQ) and DKQ (BSQ). The optimal membrane triangle OPT and optimal quadrangle OPQ perform worse than the CST and BSQ elements, possibly due to the bending-dominant situation. The DKMQ24 element performs better than DKMQ24+ in this example. Contrary to other benchmark problems in this paper, all the tested elements converge to the reference solution from below.

\begin{center}
	\begin{table}
		\centering
		\begin{tabularx}{\tabulkywidth}{c *{9}{Y}}
			\hline
			\#elements&\#dofs& \scalebox{0.75}[1.0]{Specht} \newline  \scalebox{0.75}[1.0]{(CST)}& \scalebox{0.75}[1.0]{Specht} \newline  \scalebox{0.75}[1.0]{(OPT)}& \scalebox{0.75}[1.0]{DKT} \newline  \scalebox{0.75}[1.0]{(CST)}& \scalebox{0.75}[1.0]{DKT} \newline  \scalebox{0.75}[1.0]{(OPT)}& \scalebox{0.75}[1.0]{DKMT} \newline  \scalebox{0.75}[1.0]{(CST)}& \scalebox{0.75}[1.0]{DKMT} \newline  \scalebox{0.75}[1.0]{(OPT)}& \scalebox{0.75}[1.0]{MITC3} \rule[-2mm]{0mm}{6mm} \\
			&& $\lambda$ [-] & $\lambda$ [-] & $\lambda$ [-] & $\lambda$ [-] & $\lambda$ [-] & $\lambda$ [-] & $\lambda$ [-]  \\
			\hline \hline
			$3\!\!\times\!\!3\!\!\times\!\!3\!\!\times\!\!4$ & $324      $ & $40.60$ & $166.42$ & $297.24$ & $247.07$ & $297.24$ & $245.00$ & $291.21$  \\
			$6\!\!\times\!\!6\!\!\times\!\!3\!\!\times\!\!4$ & $1296     $ & $141.80$ & $179.47$ & $308.05$ & $285.86$ & $308.03$ & $283.19$ & $289.34$  \\
			$12\!\!\times\!\!12\!\!\times\!\!3\!\!\times\!\!4$ & $5184     $ & $335.27$ & $263.51$ & $334.81$ & $329.86$ & $334.73$ & $328.61$ & $288.17$  \\
			$24\!\!\times\!\!24\!\!\times\!\!3\!\!\times\!\!4$ & $20736    $ & $370.64$ & $355.35$ & $374.57$ & $369.67$ & $374.48$ & $373.48$ & $284.64$  \\
			$48\!\!\times\!\!48\!\!\times\!\!3\!\!\times\!\!4$ & $82944    $ & $379.38$ & $373.18$ & $381.40$ & $381.25$ & $381.31$ & $381.04$ & $284.04$  \\
			\hline
			\multicolumn{9}{c}{Reference solution: $\lambda=384.00$ (\cite{Bleyer2021})} \\
			\hline
		\end{tabularx}
		\captionof{table}{Spherical cap (triangle elements). Note, that the MITC3 element does not converge. \label{Tab:B03TRI}}
	\end{table}
\end{center}

\begin{center}
	\begin{table}
		\centering
		\begin{tabularx}{\tabulkywidth}{c *{9}{Y}}
			\hline
			\#elements&\#dofs& \scalebox{0.75}[1.0]{DKQ} \newline  \scalebox{0.75}[1.0]{(BSQ)}& \scalebox{0.75}[1.0]{DKQ} \newline  \scalebox{0.75}[1.0]{(OPQ)}& \scalebox{0.75}[1.0]{DKMQ} \newline  \scalebox{0.75}[1.0]{(BSQ)}& \scalebox{0.75}[1.0]{DKMQ} \newline  \scalebox{0.75}[1.0]{(OPQ)}& \scalebox{0.75}[1.0]{DKMQ24}& \scalebox{0.75}[1.0]{DKMQ24\!+~}& \scalebox{0.75}[1.0]{MITC4} \rule[-2mm]{0mm}{6mm} \\
			&& $\lambda$ [-] & $\lambda$ [-] & $\lambda$ [-] & $\lambda$ [-] & $\lambda$ [-] & $\lambda$ [-] & $\lambda$ [-]  \\
			\hline \hline
			$3\!\!\times\!\!3\!\!\times\!\!3    $ & $162      $ & $376.54$ & $334.19$ & $376.74$ & $335.47$ & $375.84$ & $339.81$ & $375.53$  \\
			$6\!\!\times\!\!6\!\!\times\!\!3    $ & $648      $ & $378.48$ & $353.77$ & $378.45$ & $354.04$ & $379.83$ & $354.77$ & $378.46$  \\
			$12\!\!\times\!\!12\!\!\times\!\!3  $ & $2592     $ & $381.93$ & $371.80$ & $381.92$ & $367.67$ & $382.57$ & $369.32$ & $380.41$  \\
			$24\!\!\times\!\!24\!\!\times\!\!3  $ & $10368    $ & $383.38$ & $380.67$ & $383.38$ & $378.31$ & $383.45$ & $380.78$ & $380.87$  \\
			$48\!\!\times\!\!48\!\!\times\!\!3  $ & $41472    $ & $383.83$ & $382.98$ & $383.83$ & $381.73$ & $383.85$ & $382.97$ & $377.60$  \\
			\hline
			\multicolumn{9}{c}{Reference solution: $\lambda=384.00$ (\cite{Bleyer2021})} \\
			\hline
		\end{tabularx}
		\captionof{table}{Spherical cap (quad elements).  \label{Tab:B03QUAD}}
	\end{table}
\end{center}

\begin{center}
	\begin{table}
		\centering
		\begin{tabularx}{110mm}{c *{5}{Y}}
			\hline
			\#elements&\#dofs&\scalebox{0.75}[1.0]{MITC9}&\#dofs&\scalebox{0.75}[1.0]{MITC16}  \rule[-2mm]{0mm}{6mm} \\
			&& $\lambda$ [-] && $\lambda$ [-]  \\
			\hline \hline
			$3\!\!\times\!\!3\!\!\times\!\!3    $ & $666      $ & $131.74$& $1494     $ & $80.07$  \\
			$6\!\!\times\!\!6\!\!\times\!\!3    $ & $2628     $ & $143.62$& $5832     $ & $325.00$  \\
			$12\!\!\times\!\!12\!\!\times\!\!3  $ & $10368    $ & $333.94$& $23328    $ & $324.90$  \\
			$24\!\!\times\!\!24\!\!\times\!\!3  $ & $41472    $ & $333.66$& $93312    $ & $324.45$  \\
			$48\!\!\times\!\!48\!\!\times\!\!3  $ & $165888   $ & $347.08$& $373248   $ & $340.30$  \\
			\hline
			\multicolumn{5}{c}{Reference solution: $\lambda=384.00$ (\cite{Bleyer2021})} \\
			\hline
		\end{tabularx}
		\captionof{table}{Spherical cap (high order quad elements). Note, that the MITC9 and MITC16 element do not converge. \label{Tab:B03HQUAD}}
	\end{table}
\end{center}

\begin{figure}
	\centering
	\scalebox{.9}
	{
		\begin{tikzpicture}
			\begin{loglogaxis}[
				xlabel=\text{time [s]},
				ylabel= $\left|\frac{\lambda-\lambda_{\text{ref}}}{\lambda_{\text{ref}}}\right|$,
				x tick label style={
					/pgf/number format/1000 sep={}  
				},
				y tick label style={
					/pgf/number format/1000 sep={},
					/pgf/number format/fixed,       
					/pgf/number format/precision=5  
				},
				scaled y ticks=false, 
				xmin=0.23,
				xmax=2532.5,
				ymin=0.00038095572916673,
				ymax=0.894259770572917,
				legend style={at={(1.02,0.02)},anchor=south west}
				]
				\addplot[line width=1.5pt, color=red, mark=triangle*, mark size=4pt, solid] coordinates {
					(1.06, 0.894259770572917)
					(1.86, 0.63071941796875)
					(34.86, 0.12690225)
					(275.8, 0.0347791458333334)
					(2216.2, 0.0120278802083334)
				};
				\addplot[line width=1.5pt, color=red, mark=triangle*, mark size=4pt, densely dotted] coordinates {
					(1.5, 0.56660873828125)
					(2.42, 0.532642317708333)
					(29.61, 0.31376878125)
					(193.19, 0.074599828125)
					(1466.61, 0.0281870091145833)
				};
				\addplot[line width=1.5pt, color=brown, mark=triangle*, mark size=3pt, solid] coordinates {
					(0.94, 0.225929434895833)
					(2, 0.197799376302083)
					(27.11, 0.128098509114583)
					(321.94, 0.0245691666666667)
					(2428.05, 0.00678264192708335)
				};
				\addplot[line width=1.5pt, color=brown, mark=triangle*, mark size=3pt, densely dotted] coordinates {
					(1.13, 0.35659196875)
					(2.39, 0.255566009114583)
					(23.25, 0.140996455729167)
					(209.28, 0.0373223385416666)
					(1566.11, 0.00716183072916669)
				};
				\addplot[line width=1.5pt, color=orange, mark=triangle*, mark size=5pt, solid] coordinates {
					(0.81, 0.225935109375)
					(1.92, 0.197830467447917)
					(25.67, 0.128295674479167)
					(311.48, 0.0248045625)
					(2532.5, 0.00701608984375007)
				};
				\addplot[line width=1.5pt, color=orange, mark=triangle*, mark size=5pt, densely dotted] coordinates {
					(1.06, 0.3619703125)
					(2.47, 0.262525053385417)
					(29.78, 0.14424908984375)
					(279.77, 0.0274023020833334)
					(1859.09, 0.00771385026041672)
				};
				\addplot[line width=1.5pt, color=yellow, mark=diamond, mark size=9pt, solid] coordinates {
					(0.23, 0.01943287890625)
					(1.3, 0.0143732656250001)
					(5.34, 0.00539754036458332)
					(40.36, 0.00162344921875005)
					(307.8, 0.000436917968750066)
				};
				\addplot[line width=1.5pt, color=yellow, mark=diamond, mark size=7pt, densely dotted] coordinates {
					(0.28, 0.129709071614583)
					(1.2, 0.0787331744791667)
					(4.56, 0.0317650260416666)
					(32.16, 0.00868300130208337)
					(224.13, 0.00266751171874999)
				};
				\addplot[line width=3pt, color=blue, mark=square, mark size=3pt, solid] coordinates {
					(0.27, 0.0189023580729167)
					(1.14, 0.0144517591145834)
					(4.75, 0.00540477864583331)
					(33.73, 0.00162675260416673)
					(255, 0.000438197916666701)
				};
				\addplot[line width=3pt, color=blue, mark=square, mark size=3pt, densely dotted] coordinates {
					(0.27, 0.126375912760417)
					(1.13, 0.0780199986979166)
					(4.22, 0.0425289049479167)
					(33.14, 0.01481973828125)
					(247.61, 0.0059070638020834)
				};
				\addplot[line width=1.5pt, color=green, mark=square*, mark size=4pt] coordinates {
					(0.27, 0.0212477552083333)
					(1.42, 0.0108500208333333)
					(5.59, 0.00373063671874999)
					(38.2, 0.00142447135416671)
					(344.25, 0.00038095572916673)
				};
				\addplot[line width=1.5pt, color=col1, mark=square*, mark size=4pt] coordinates {
					(0.27, 0.11506812109375)
					(1.13, 0.0761313736979167)
					(5, 0.0382261744791667)
					(30.2, 0.00838111718750002)
					(234.28, 0.00267416406250006)
				};
				\addplot[line width=1pt, color=black, mark=square, mark size=3pt] coordinates {
					(0.38, 0.0220456054687501)
					(1.58, 0.0144302421875001)
					(6.09, 0.00935649739583341)
					(31.44, 0.00816069140625)
					(210.55, 0.0166721302083334)
				};
				\legend{Specht (CST),Specht (OPT),DKT (CST),DKT (OPT),DKMT (CST),DKMT (OPT),DKQ (BSQ),DKQ (OPQ),DKMQ (BSQ),DKMQ (OPQ),DKMQ24,DKMQ24\!+~,MITC4}
			\end{loglogaxis}
		\end{tikzpicture}
	}
	\caption{B03: Spherical cap }
	\label{fig:B03TIME}
\end{figure}
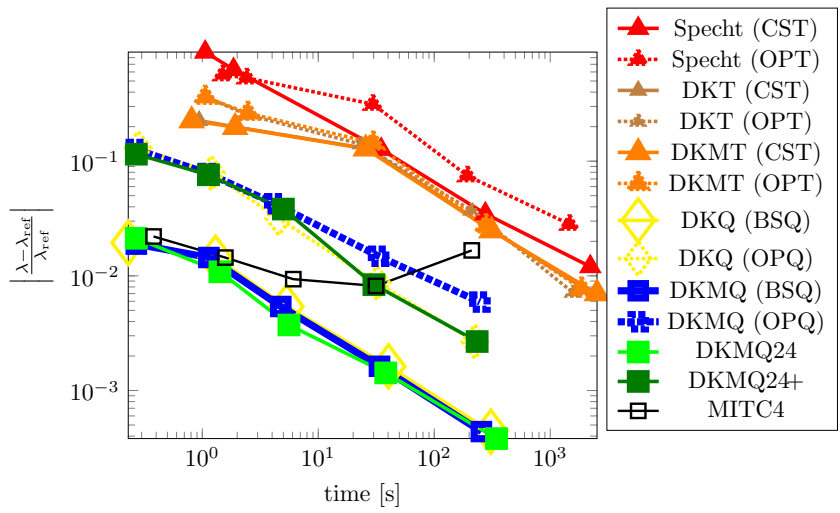

\begin{figure}[H]
	\centering
	\begin{overpic}[width=0.4\textwidth]{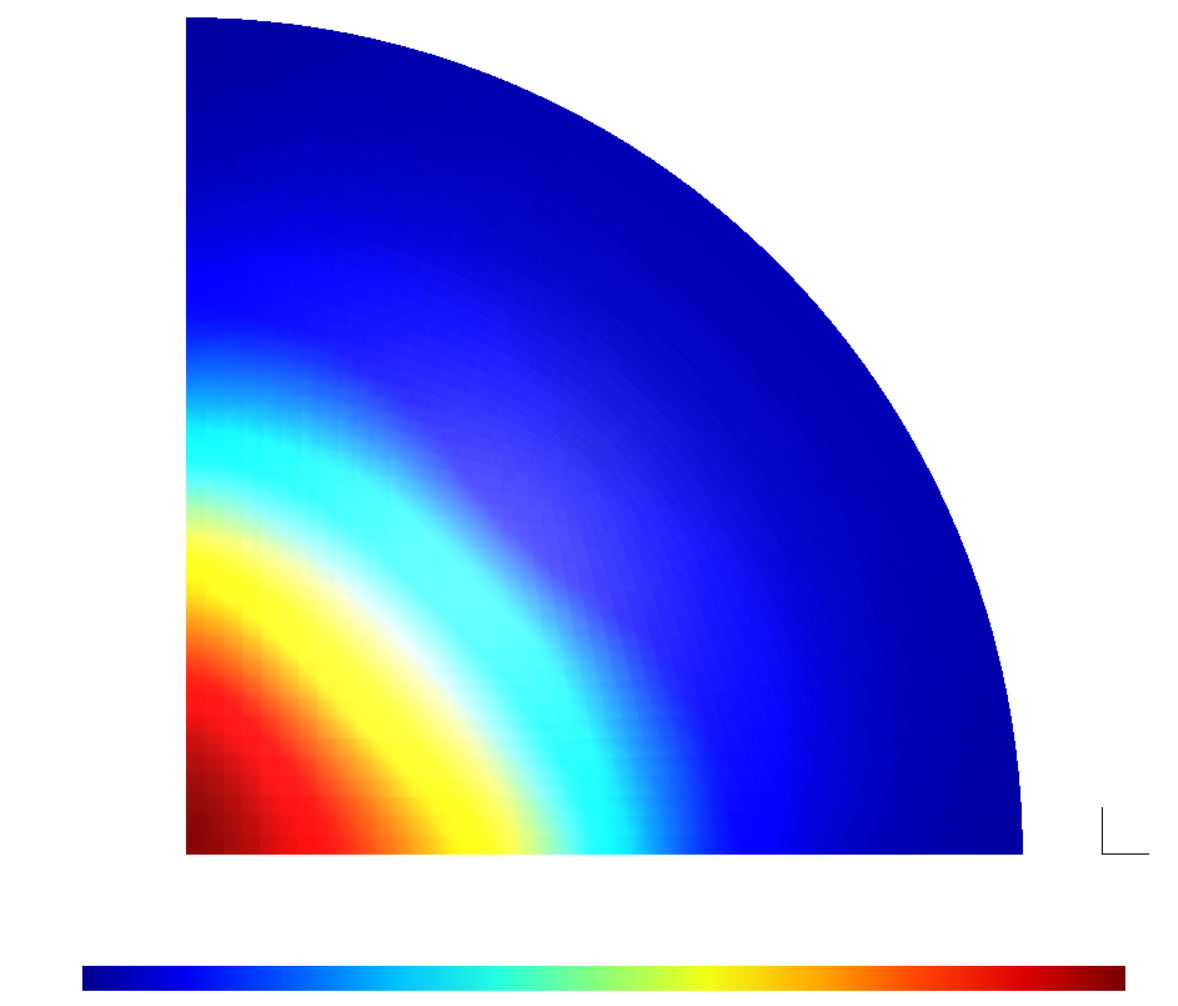} 
		\put(95.5,7.8){$X$}
\put(89.1,15){$Y$}
		\put(4,-4.5){$0$}
		\put(44.5,-4.5){$0.149$}
		\put(88,-4.5){$0.297$}
		\put(46.5,5.5){$D_{\p}$}
	\end{overpic}
	\caption{Simply supported spherical cap under uniform pressure -- dissipation power at collapse (DKMQ24 element).}
	\label{SphericalCapDissipation}
\end{figure}

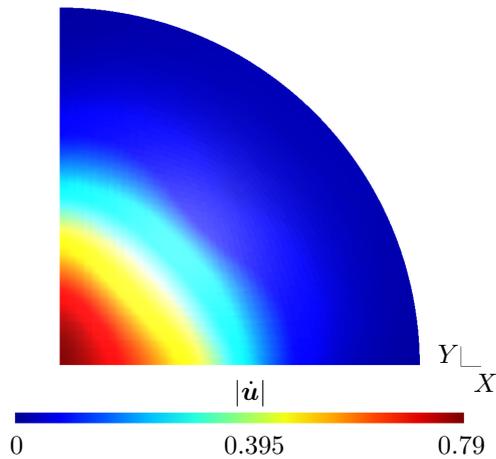
\begin{figure}[H]
	\centering
	\begin{overpic}[width=0.4\textwidth]{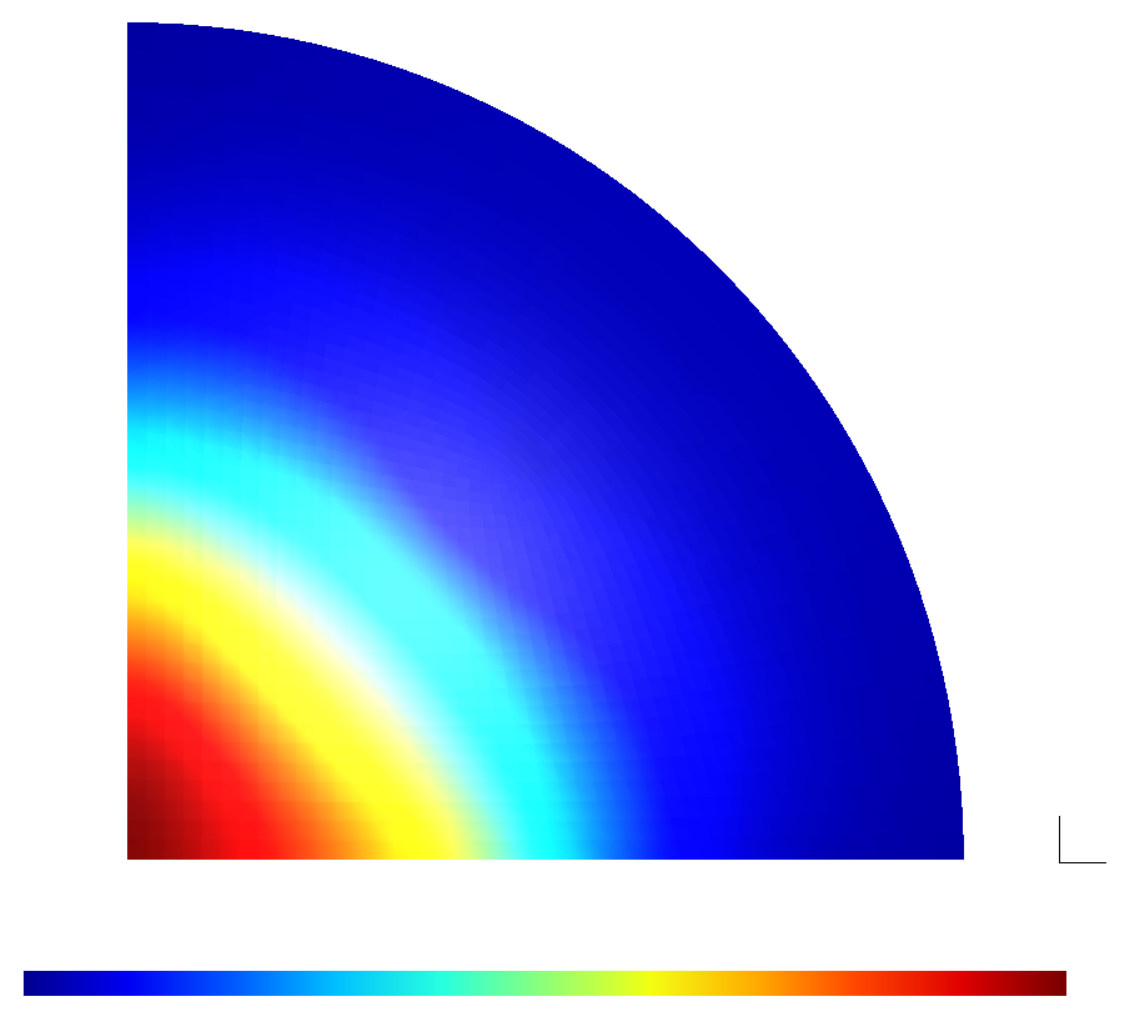} 
		\put(95.3,8){$X$}
		\put(88,14){$Y$}
		\put(1,-4.5){$0$}
		\put(44.5,-4.5){$0.395$}
		\put(88,-4.5){$0.79$}
		\put(46.5,6.5){$|\dot{\bm{u}}|$}
	\end{overpic}
	\caption{Simply supported spherical cap under uniform pressure -- absolute value of velocity at collapse (DKMQ24 element).}
	\label{SphericalCapVelocity}
\end{figure}


\FloatBarrier
\subsection{Clamped cylindrical tube under vertical distributed loading}
\label{b4}

\begin{figure}[H]
	\centering
	\begin{overpic}[width=0.6\textwidth]{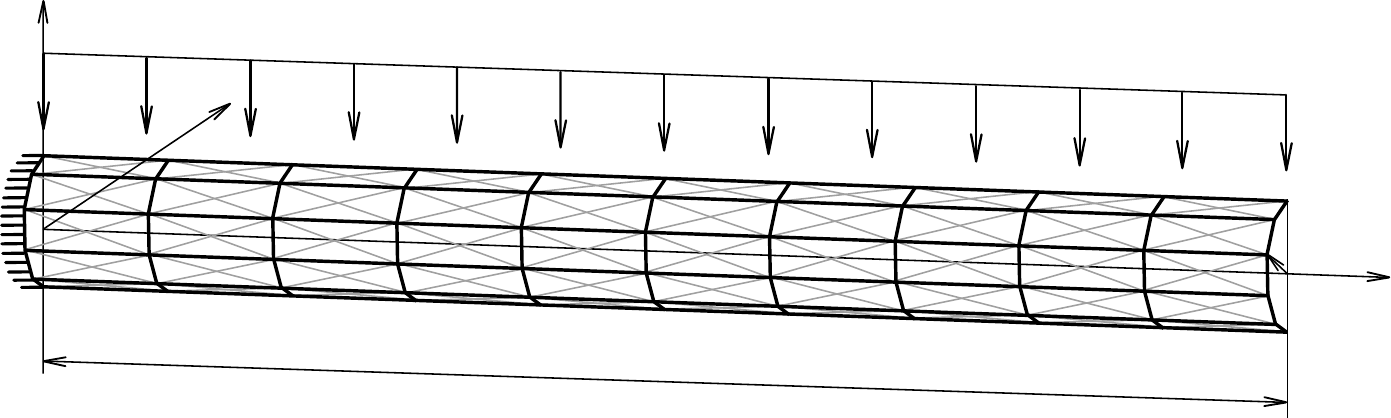} 
		\put(45,27.5){$f_Z$}
		\put(73.5,2.6){thickness $h$}
		\put(93,11.2){$R$}
		\put(49,3.3){$L$}
		\put(97,6){$X$}
		\put(13.5,22.5){$Y$}
		\put(-0.75,27.6){$Z$}
	\end{overpic}
	\caption{Clamped cylindrical tube under vertical distributed loading with quad mesh $10\times 5$, a triangular mesh depicted in a grey colour. Parameters: $R=0.05\,\text{m},\,h=0.001\,\text{m},\,L=2\,\text{m}$, yield strength $\sigma_0=240\,\text{MPa}$, loading force $f_Z=-1\,\text{N}$.}
	\label{clampedTube}
\end{figure}
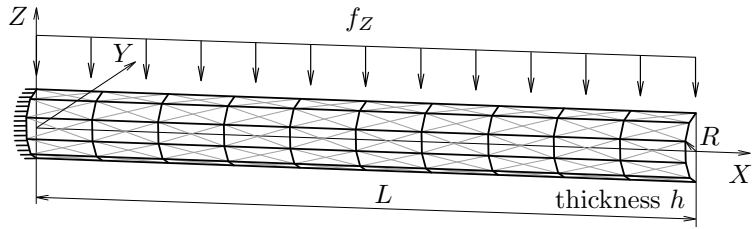

A clamped cylindrical tube of length $2L$, thickness $h$ and radius $R$ under vertical distributed loading $f_Z$ is considered (see Fig. \ref{clampedTube}). The Ilyushin failure criterion with yield stress $\sigma_0$ is considered. This benchmark is taken from \cite{Bleyer2016}. Only one quarter is computed due to symmetry. Boundary and symmetry conditions are as follows: $u_X=u_Y=u_Z=\varphi_X=\varphi_Y=\varphi_Z=0$ at $X=0$, $u_X=\varphi_Y=\varphi_Z=0$ at $X=L=2$, $u_Y=\varphi_X=\varphi_Z=0$ at $Y=0$. The problem has the following approximate analytical solution based on the beam theory (fully valid for $h\ll R\ll2L$): the ultimate line load on the rigid clamped beam of length $2L$ is $q=16\frac{M}{(2L)^2}\,\left[\frac{\text{N}}{\text{m}}\right]$, where the ultimate moment is given by $M \approx 4hR^2\sigma_0\,\left[\text{Nm}\right]$ (\cite{Bleyer2016}). The ultimate load on the simulated quarter of the tube is given by $F_Z=\frac{qL}{2}$, which gives us the final formula
\begin{align}
	F_Z=8\frac{hR^2\sigma_0}{L} = 4800\,\text{N}.
\end{align}
The numerical results are summarized in Tables \ref{Tab:B04TRI}, \ref{Tab:B04QUAD} and \ref{Tab:B04HQUAD}. Let us discuss the results depicted in Fig. \ref{fig:B04}. None of the MITC elements converge. The DKQ element in both variants (DKQ (BSQ/OPQ)) does not perform well, the situation even worsens for smaller tube radius (we have tried $R=0.03\,\text{m}$). The DKMT (OPT) and Specht (OPT) perform surprisingly well in this example. The most efficient elements in this example are quadrangles DKMQ (BSQ), DKMQ24, and DKMQ24+ (both with respect to the number of degrees of freedom and time). DKMQ24+ outperforms the first two elements. All tested elements, except for the DKT (OPT) element at the coarsest mesh, converge monotonically from above to the reference solution.

\FloatBarrier
 
\begin{center}
	\begin{table}
		\centering
		\begin{tabularx}{\tabulkywidth}{c *{9}{Y}}
			\hline
			\#elements&\#dofs& \scalebox{0.75}[1.0]{Specht} \newline  \scalebox{0.75}[1.0]{(CST)}& \scalebox{0.75}[1.0]{Specht} \newline  \scalebox{0.75}[1.0]{(OPT)}& \scalebox{0.75}[1.0]{DKT} \newline  \scalebox{0.75}[1.0]{(CST)}& \scalebox{0.75}[1.0]{DKT} \newline  \scalebox{0.75}[1.0]{(OPT)}& \scalebox{0.75}[1.0]{DKMT} \newline  \scalebox{0.75}[1.0]{(CST)}& \scalebox{0.75}[1.0]{DKMT} \newline  \scalebox{0.75}[1.0]{(OPT)}& \scalebox{0.75}[1.0]{MITC3} \rule[-2mm]{0mm}{6mm} \\
			&& $\lambda$ [-] & $\lambda$ [-] & $\lambda$ [-] & $\lambda$ [-] & $\lambda$ [-] & $\lambda$ [-] & $\lambda$ [-]  \\
			\hline \hline
			$10\!\!\times\!\!5\!\!\times\!\!4   $ & $284      $ & $6587.3$ & $5234.9$ & $6587.3$ & $5726.7$ & $6587.3$ & $5232.6$ & $6172.2$  \\
			$20\!\!\times\!\!10\!\!\times\!\!4  $ & $1169     $ & $5705.6$ & $5151.9$ & $5694.1$ & $5908.3$ & $5694.1$ & $5150.1$ & $5207.6$  \\
			$40\!\!\times\!\!20\!\!\times\!\!4  $ & $4739     $ & $5282.2$ & $5073.7$ & $5287.1$ & $5386.8$ & $5287.0$ & $5065.2$ & $4579.9$  \\
			$80\!\!\times\!\!40\!\!\times\!\!4  $ & $19079    $ & $5085.1$ & $4993.4$ & $5070.9$ & $5147.1$ & $5070.7$ & $5022.1$ & $3754.2$  \\
			$160\!\!\times\!\!80\!\!\times\!\!4 $ & $76559    $ & $4979.5$ & $4939.7$ & $4977.8$ & $5019.8$ & $4977.5$ & $4939.3$ & $3903.3$  \\
			\hline
			\multicolumn{9}{c}{Reference solution: $\lambda=4800.0$ (DKMQ24 element at mesh $160\!\!\times\!\!80$)} \\
			\hline
		\end{tabularx}
		\captionof{table}{Clamped cylindrical tube (triangle elements). Note, that the MITC3 element does not converge. \label{Tab:B04TRI}}
	\end{table}
\end{center}

\begin{center}
	\begin{table}
		\centering
		\begin{tabularx}{\tabulkywidth}{c *{9}{Y}}
			\hline
			\#elements&\#dofs& \scalebox{0.75}[1.0]{DKQ} \newline  \scalebox{0.75}[1.0]{(BSQ)}& \scalebox{0.75}[1.0]{DKQ} \newline  \scalebox{0.75}[1.0]{(OPQ)}& \scalebox{0.75}[1.0]{DKMQ} \newline  \scalebox{0.75}[1.0]{(BSQ)}& \scalebox{0.75}[1.0]{DKMQ} \newline  \scalebox{0.75}[1.0]{(OPQ)}& \scalebox{0.75}[1.0]{DKMQ24}& \scalebox{0.75}[1.0]{DKMQ24\!+~}& \scalebox{0.75}[1.0]{MITC4} \rule[-2mm]{0mm}{6mm} \\
			&& $\lambda$ [-] & $\lambda$ [-] & $\lambda$ [-] & $\lambda$ [-] & $\lambda$ [-] & $\lambda$ [-] & $\lambda$ [-]  \\
			\hline \hline
			$10\!\!\times\!\!5                  $ & $284      $ & $6022.8$ & $5529.1$ & $5909.6$ & $5759.0$ & $5909.6$ & $5365.9$ & $5622.0$  \\
			$20\!\!\times\!\!10                 $ & $1169     $ & $5530.0$ & $5322.1$ & $5407.5$ & $5637.5$ & $5407.5$ & $5177.5$ & $4854.2$  \\
			$40\!\!\times\!\!20                 $ & $4739     $ & $5255.6$ & $5189.7$ & $5122.4$ & $5201.4$ & $5122.4$ & $5046.4$ & $4369.7$  \\
			$80\!\!\times\!\!40                 $ & $19079    $ & $5116.7$ & $5104.2$ & $4980.5$ & $5026.3$ & $4980.5$ & $4980.1$ & $4110.7$  \\
			$160\!\!\times\!\!80                $ & $76559    $ & $5067.1$ & $5072.5$ & $4926.6$ & $4977.7$ & $4926.6$ & $4925.7$ & $3973.9$  \\
			\hline
			\multicolumn{9}{c}{Reference solution: $\lambda=4800.0$ (DKMQ24 element at mesh $160\!\!\times\!\!80$)} \\
			\hline
		\end{tabularx}
		\captionof{table}{Clamped cylindrical tube (quad elements). Note, that the MITC4 element does not converge. \label{Tab:B04QUAD}}
	\end{table}
\end{center}

\begin{center}
	\begin{table}
		\centering
		\begin{tabularx}{110mm}{c *{5}{Y}}
			\hline
			\#elements&\#dofs&\scalebox{0.75}[1.0]{MITC9}&\#dofs&\scalebox{0.75}[1.0]{MITC16}  \rule[-2mm]{0mm}{6mm} \\
			&& $\lambda$ [-] && $\lambda$ [-]  \\
			\hline \hline
			$10\!\!\times\!\!5                  $ & $1169     $ & $0.0$& $2654     $ & $0.0$  \\
			$20\!\!\times\!\!10                 $ & $4739     $ & $0.0$& $10709    $ & $0.0$  \\
			$40\!\!\times\!\!20                 $ & $19079    $ & $0.1$& $43019    $ & $0.0$  \\
			$80\!\!\times\!\!40                 $ & $76559    $ & $0.0$& $172439   $ & $0.0$  \\
			$160\!\!\times\!\!80                $ & $306719   $ & $0.0$& $690479   $ & $0.0$  \\
			\hline
			\multicolumn{5}{c}{Reference solution: $\lambda=4800.0$ (DKMQ24 element at mesh $160\!\!\times\!\!80$)} \\
			\hline
		\end{tabularx}
		\captionof{table}{Clamped cylindrical tube (high order quad elements). Note, that the MITC9 and MITC16 element do not converge. \label{Tab:B04HQUAD}}
	\end{table}
\end{center}

\begin{figure}
	\centering
	\scalebox{.9}
	{
		\begin{tikzpicture}
			\begin{loglogaxis}[
				xlabel=\text{total number of degrees of freedom},
				ylabel= $\left|\frac{\lambda-\lambda_{\text{ref}}}{\lambda_{\text{ref}}}\right|$,
				x tick label style={
					/pgf/number format/1000 sep={}  
				},
				y tick label style={
					/pgf/number format/1000 sep={},
					/pgf/number format/fixed,       
					/pgf/number format/precision=5  
				},
				scaled y ticks=false, 
				xmin=284,
				xmax=76559,
				ymin=0.0261918364583335,
				ymax=0.372348413541667,
				legend style={at={(1.02,0.02)},anchor=south west}
				]
				\addplot[line width=1.5pt, color=red, mark=triangle*, mark size=4pt, solid] coordinates {
					(284, 0.372347917708333)
					(1169, 0.188671255208333)
					(4739, 0.100464672916667)
					(19079, 0.0593859572916668)
					(76559, 0.037392025)
				};
				\addplot[line width=1.5pt, color=red, mark=triangle*, mark size=4pt, densely dotted] coordinates {
					(284, 0.0906142458333334)
					(1169, 0.0733133552083334)
					(4739, 0.0570167729166668)
					(19079, 0.0402973531250001)
					(76559, 0.0291117364583334)
				};
				\addplot[line width=1.5pt, color=brown, mark=triangle*, mark size=3pt, solid] coordinates {
					(284, 0.372348413541667)
					(1169, 0.186281189583333)
					(4739, 0.1014762875)
					(19079, 0.0564346250000001)
					(76559, 0.0370361437500001)
				};
				\addplot[line width=1.5pt, color=brown, mark=triangle*, mark size=3pt, densely dotted] coordinates {
					(284, 0.19306578125)
					(1169, 0.230905655208334)
					(4739, 0.122259069791667)
					(19079, 0.0723227843750001)
					(76559, 0.0457967020833335)
				};
				\addplot[line width=1.5pt, color=orange, mark=triangle*, mark size=5pt, solid] coordinates {
					(284, 0.372346161458334)
					(1169, 0.186276082291667)
					(4739, 0.101450247916667)
					(19079, 0.0564011520833335)
					(76559, 0.0369855177083335)
				};
				\addplot[line width=1.5pt, color=orange, mark=triangle*, mark size=5pt, densely dotted] coordinates {
					(284, 0.0901256197916668)
					(1169, 0.0729288375000001)
					(4739, 0.0552451375000001)
					(19079, 0.0462644218750001)
					(76559, 0.0290217104166668)
				};
				\addplot[line width=1.5pt, color=yellow, mark=diamond, mark size=9pt, solid] coordinates {
					(284, 0.254759140625)
					(1169, 0.152082715625)
					(4739, 0.0949144739583334)
					(19079, 0.0659799614583334)
					(76559, 0.055650871875)
				};
				\addplot[line width=1.5pt, color=yellow, mark=diamond, mark size=7pt, densely dotted] coordinates {
					(284, 0.151896970833333)
					(1169, 0.108773222916667)
					(4739, 0.0811828875000001)
					(19079, 0.0633736041666667)
					(76559, 0.0567648791666668)
				};
				\addplot[line width=3pt, color=blue, mark=square, mark size=3pt, solid] coordinates {
					(284, 0.231158597916667)
					(1169, 0.126564373958333)
					(4739, 0.0671741083333335)
					(19079, 0.0376071989583334)
					(76559, 0.0263723437500002)
				};
				\addplot[line width=3pt, color=blue, mark=square, mark size=3pt, densely dotted] coordinates {
					(284, 0.199794016666667)
					(1169, 0.174472721875)
					(4739, 0.0836185614583334)
					(19079, 0.0471494364583334)
					(76559, 0.0370208416666667)
				};
				\addplot[line width=1.5pt, color=green, mark=square*, mark size=4pt] coordinates {
					(284, 0.231158597916667)
					(1169, 0.126564373958333)
					(4739, 0.0671741083333335)
					(19079, 0.0376071979166667)
					(76559, 0.0263723406250001)
				};
				\addplot[line width=1.5pt, color=col1, mark=square*, mark size=4pt] coordinates {
					(284, 0.117903359375)
					(1169, 0.0786470677083334)
					(4739, 0.0513421781250002)
					(19079, 0.0375197468750001)
					(76559, 0.0261918364583335)
				};
				\legend{Specht (CST),Specht (OPT),DKT (CST),DKT (OPT),DKMT (CST),DKMT (OPT),DKQ (BSQ),DKQ (OPQ),DKMQ (BSQ),DKMQ (OPQ),DKMQ24,DKMQ24\!+~}
			\end{loglogaxis}
		\end{tikzpicture}
	}
	\caption{B04: Clamped cylindrical tube }
	\label{fig:B04}
\end{figure}

\begin{figure}
	\centering
	\scalebox{.9}
	{
		\begin{tikzpicture}
			\begin{loglogaxis}[
				xlabel=\text{time [s]},
				ylabel= $\left|\frac{\lambda-\lambda_{\text{ref}}}{\lambda_{\text{ref}}}\right|$,
				x tick label style={
					/pgf/number format/1000 sep={}  
				},
				y tick label style={
					/pgf/number format/1000 sep={},
					/pgf/number format/fixed,       
					/pgf/number format/precision=5  
				},
				scaled y ticks=false, 
				xmin=0.11,
				xmax=439.47,
				ymin=0.0261918364583335,
				ymax=0.372348413541667,
				legend style={at={(1.02,0.02)},anchor=south west}
				]
				\addplot[line width=1.5pt, color=red, mark=triangle*, mark size=4pt, solid] coordinates {
					(0.33, 0.372347917708333)
					(0.75, 0.188671255208333)
					(5.08, 0.100464672916667)
					(50.95, 0.0593859572916668)
					(393.25, 0.037392025)
				};
				\addplot[line width=1.5pt, color=red, mark=triangle*, mark size=4pt, densely dotted] coordinates {
					(0.28, 0.0906142458333334)
					(0.78, 0.0733133552083334)
					(7.61, 0.0570167729166668)
					(51.05, 0.0402973531250001)
					(359.63, 0.0291117364583334)
				};
				\addplot[line width=1.5pt, color=brown, mark=triangle*, mark size=3pt, solid] coordinates {
					(0.27, 0.372348413541667)
					(0.67, 0.186281189583333)
					(4.69, 0.1014762875)
					(56.02, 0.0564346250000001)
					(357.92, 0.0370361437500001)
				};
				\addplot[line width=1.5pt, color=brown, mark=triangle*, mark size=3pt, densely dotted] coordinates {
					(0.25, 0.19306578125)
					(0.83, 0.230905655208334)
					(7.2, 0.122259069791667)
					(47.11, 0.0723227843750001)
					(333.41, 0.0457967020833335)
				};
				\addplot[line width=1.5pt, color=orange, mark=triangle*, mark size=5pt, solid] coordinates {
					(0.3, 0.372346161458334)
					(0.73, 0.186276082291667)
					(4.36, 0.101450247916667)
					(54.92, 0.0564011520833335)
					(439.47, 0.0369855177083335)
				};
				\addplot[line width=1.5pt, color=orange, mark=triangle*, mark size=5pt, densely dotted] coordinates {
					(0.3, 0.0901256197916668)
					(0.75, 0.0729288375000001)
					(7.39, 0.0552451375000001)
					(47.56, 0.0462644218750001)
					(364.76, 0.0290217104166668)
				};
				\addplot[line width=1.5pt, color=yellow, mark=diamond, mark size=9pt, solid] coordinates {
					(0.2, 0.254759140625)
					(0.81, 0.152082715625)
					(4.49, 0.0949144739583334)
					(28.59, 0.0659799614583334)
					(225.52, 0.055650871875)
				};
				\addplot[line width=1.5pt, color=yellow, mark=diamond, mark size=7pt, densely dotted] coordinates {
					(0.19, 0.151896970833333)
					(0.78, 0.108773222916667)
					(4.39, 0.0811828875000001)
					(31.25, 0.0633736041666667)
					(164.63, 0.0567648791666668)
				};
				\addplot[line width=3pt, color=blue, mark=square, mark size=3pt, solid] coordinates {
					(0.14, 0.231158597916667)
					(0.83, 0.126564373958333)
					(4.75, 0.0671741083333335)
					(30.45, 0.0376071989583334)
					(235.88, 0.0263723437500002)
				};
				\addplot[line width=3pt, color=blue, mark=square, mark size=3pt, densely dotted] coordinates {
					(0.19, 0.199794016666667)
					(0.97, 0.174472721875)
					(4.11, 0.0836185614583334)
					(27.8, 0.0471494364583334)
					(161.44, 0.0370208416666667)
				};
				\addplot[line width=1.5pt, color=green, mark=square*, mark size=4pt] coordinates {
					(0.16, 0.231158597916667)
					(0.78, 0.126564373958333)
					(4.5, 0.0671741083333335)
					(30.59, 0.0376071979166667)
					(245.72, 0.0263723406250001)
				};
				\addplot[line width=1.5pt, color=col1, mark=square*, mark size=4pt] coordinates {
					(0.11, 0.117903359375)
					(0.77, 0.0786470677083334)
					(4.5, 0.0513421781250002)
					(25.91, 0.0375197468750001)
					(180.03, 0.0261918364583335)
				};
				\legend{Specht (CST),Specht (OPT),DKT (CST),DKT (OPT),DKMT (CST),DKMT (OPT),DKQ (BSQ),DKQ (OPQ),DKMQ (BSQ),DKMQ (OPQ),DKMQ24,DKMQ24\!+~}
			\end{loglogaxis}
		\end{tikzpicture}
	}
	\caption{B04: Clamped cylindrical tube }
	\label{fig:B04TIME}
\end{figure}

\begin{figure}[H]
	\centering
	\begin{overpic}[width=0.6\textwidth]{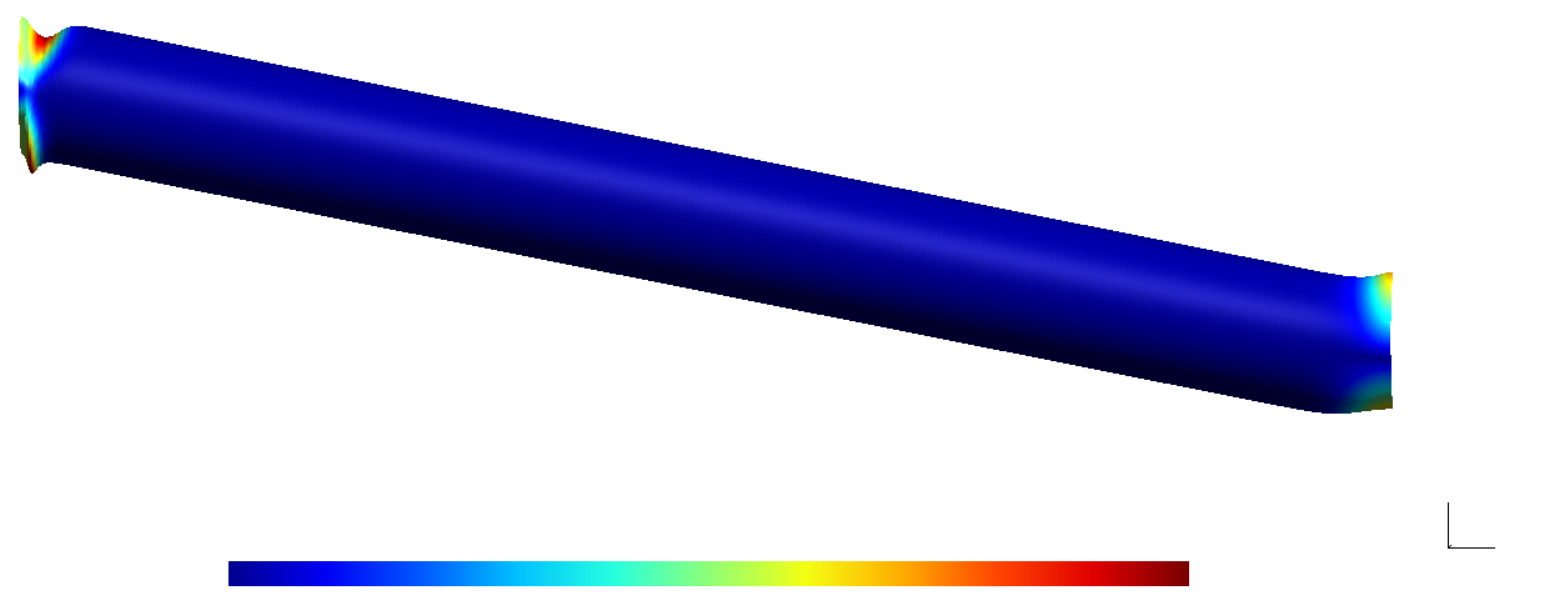} 
		\put(95.5,0.8){$X$}
		\put(91,5.5){$Z$}
		\put(17,-2.6){$0$}
		\put(76,-2.6){$17.2$}
		\put(46.5,4.5){$D_{\p}$}
	\end{overpic}
	\caption{Clamped cylindrical tube -- dissipation power at collapse (DKMQ24 element).}
	\label{CylindricalTubDissipation}
\end{figure}

\begin{figure}[H]
	\centering
	\begin{overpic}[width=0.6\textwidth]{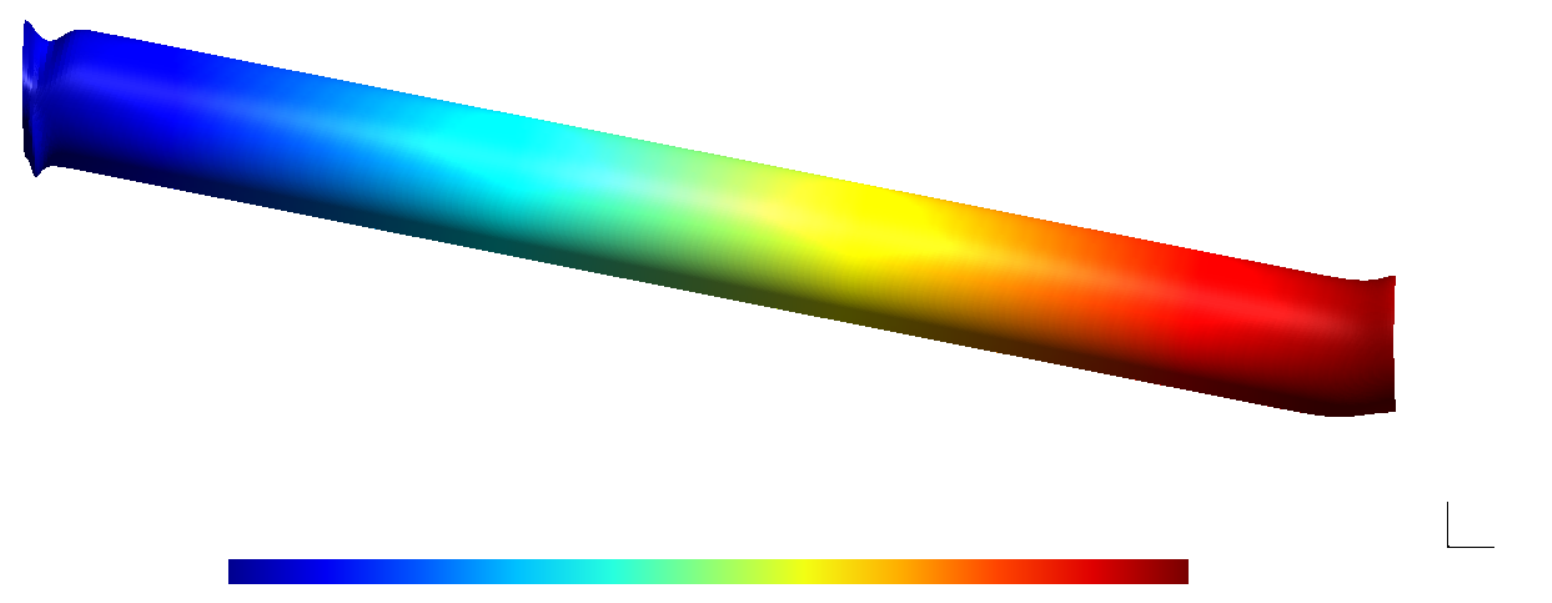} 
		\put(95.8,0.75){$X$}
		\put(91,6){$Z$}
		\put(17,-2.25){$0$}
		\put(74,-2.25){$1.96$}
		\put(46.5,5){$|\dot{\bm{u}}|$}
	\end{overpic}
	\caption{Clamped cylindrical tube -- absolute value of velocity at collapse (DKMQ24 element).} 
	\label{CylindricalTubVelocity}
\end{figure}
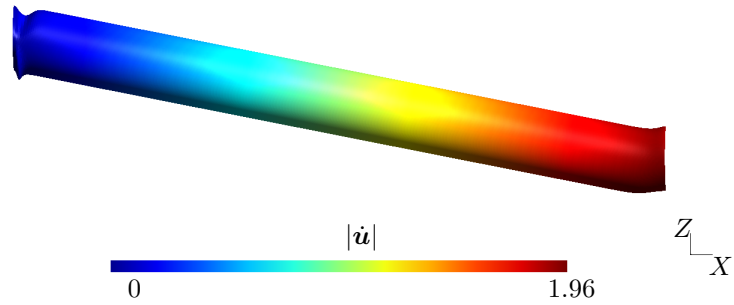
 
\FloatBarrier

\subsection{Block with asymmetric holes} 
\label{b5}

\begin{figure}[H]
	\centering
	\begin{overpic}[width=0.44\textwidth]{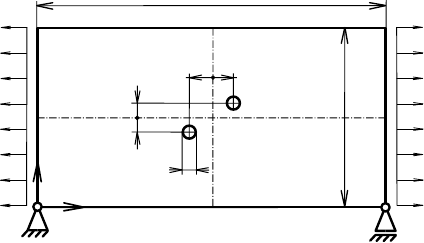} 
		\put(-6,28){$p$}
		\put(102.5,28){$p$}
		\put(43,11.8){$h_1$}
		\put(25,25.4){$h_2$}
		\put(25,31.0){$h_2$}
		\put(44.5,41.3){$h_3$}
		\put(50.9,41.3){$h_3$}
		\put(48,59){$l_X$}
		\put(75.5,24){$l_Y$}
		\put(15,2.6){$X$}
		\put(10.6,15){$Y$}
	\end{overpic}
	\caption{Block with asymmetric holes. Parameters: $l_X=0.1\,\text{m},\, l_Y=0.05\,\text{m},\, h_1=4\,\text{mm},\, h_2=2.5\,\text{mm},\, h_3=5.5\,\text{mm}, p=1\,\text{MPa}$, strengths (recalculated from \cite{Ming2}) $f_{\text{t}x}=3.3\,\text{MPa}$, $f_{\text{c}x}=4.6\,\text{MPa}$, $f_{\text{t}y}=79\,\text{MPa}$, $f_{\text{c}y}=52\,\text{MPa}$, $f_{\text{v}xy}=9.2\,\text{MPa}$.}
	\label{blockWithAsymmetricHoles}
\end{figure}

The block with asymmetric holes under tensile loading, a commonly-used benchmark for FELA, was first studied by Zouain et al. in 2002 (\cite{Zouain}). The block of length $l_X$, height $h_Y$ and thickness $h$ is subjected to pressure $p$ in $X$, see Figure \ref{blockWithAsymmetricHoles}. Boundary and symmetry conditions are as follows: $u_X=u_Y=u_Z=\varphi_X=\varphi_Y=\varphi_Z=0$ at $[X,Y]=[0,0]$, $u_Y=u_Z=\varphi_X=\varphi_Y=\varphi_Z=0$ at $[X,Y]=[l_X,0]$ and $u_Z=\varphi_X=\varphi_Y=0$ everywhere. The plane strain orthotropic Tsai-Wu failure criterion with fibers oriented in the $Y$ direction is considered, with strengths values taken from \cite{Ming1}. Since the analytical solution is not available, we use Aitken's extrapolation of the three finest results of the MITC9 element as the reference ($\lambda=3.1067$) which seems to be most probable limit value in this example. It should be noted that the solution given in \cite{Ming2} for this plane strain problem is calculated approximately on solid elements with double-symmetric conditions in the $Y$-direction and therefore converges to a slightly smaller value ($\lambda=3.082$). We also note that in this in-plane test problem, we eliminate drilling rotations $\varphi_Z$ in all elements that do not use the drilling rotations to increase their numerical efficiency (CST, MITC3, BSQ, DKMQ24, MITC3, MITC9, MICT16). Results are summarized in Tables \ref{Tab:B05TRI}, \ref{Tab:B05QUAD} and \ref{Tab:B05HQUAD}. Let us discuss the results depicted in Fig. \ref{fig:B05}. Quadrangle elements again outperform all tested triangle elements, in particular the elements OPQ, MITC9, and MITC16. The DKMQ24 element performs better than the DKMQ24+ element in this example. All tested elements converge monotonically from above to the reference solution.
 
\FloatBarrier

\begin{figure}[H]
	\centering
	\begin{overpic}[width=0.29\textwidth]{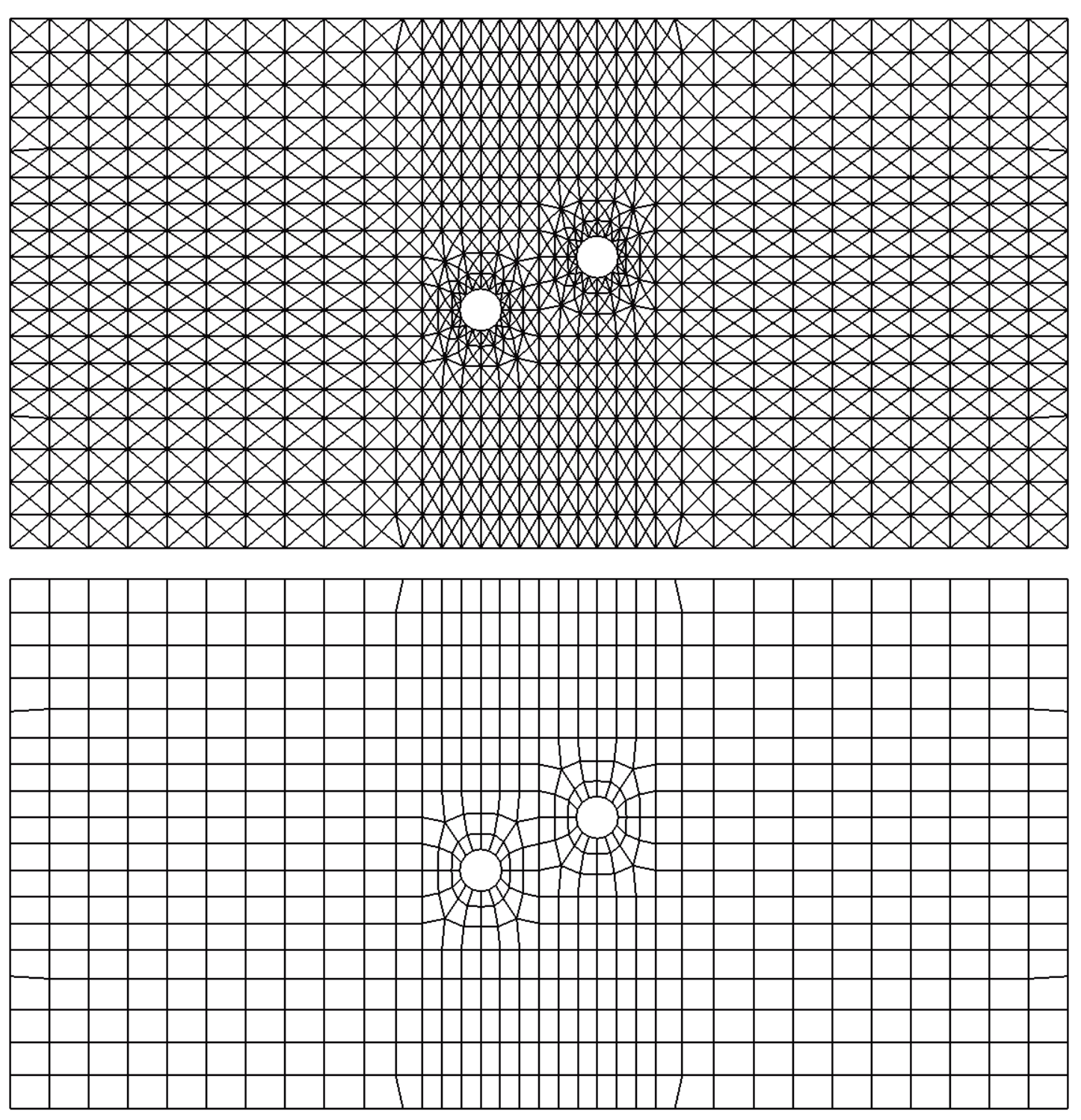} 
	\end{overpic}
	\caption{Block with asymmetric holes -- above:  triangular mesh with 2576 elements, below: quad mesh with 644 elements.}
	\label{blockWithAsymmetricHolesMesh}
\end{figure}
\FloatBarrier

\begin{center}
	\begin{table}
		\centering
		\begin{tabularx}{\tabulkywidth}{c *{6}{Y}}
			\hline
			\#elements&\#dofs& \scalebox{0.75}[1.0]{OPT}&\#dofs& \scalebox{0.75}[1.0]{CST}& \scalebox{0.75}[1.0]{MITC3} \rule[-2mm]{0mm}{6mm} \\
			&& $\lambda$ [-] && $\lambda$ [-] & $\lambda$ [-]  \\
			\hline \hline
			$644   $ & $1060     $ & $3.2862$ & $707      $ & $3.2847$ & $3.3912$  \\
			$2576  $ & $4060     $ & $3.1789$ & $2707     $ & $3.1779$ & $3.2397$  \\
			$10304 $ & $15856    $ & $3.1424$ & $10571    $ & $3.1398$ & $3.1737$  \\
			$41216 $ & $62632    $ & $3.1238$ & $41755    $ & $3.1221$ & $3.1372$  \\
			\hline
			\multicolumn{6}{c}{Reference solution: $\lambda=3.1067$ (Aitken's extrapolation of the MITC9 results)} \\
			\hline
		\end{tabularx}
		\captionof{table}{Block with asymmetric holes (triangle elements).  \label{Tab:B05TRI}}
	\end{table}
\end{center}

\begin{center}
	\begin{table}
		\centering
		\begin{tabularx}{\tabulkywidth}{c *{8}{Y}}
			\hline
			\#elements 
			& \#dofs 
			& \scalebox{0.75}[1.0]{OPQ} 
			& \scalebox{0.75}[1.0]{DKMQ24\!+~} 
			& \#dofs 
			& \scalebox{0.75}[1.0]{BSQ} 
			& \scalebox{0.75}[1.0]{DKMQ24}
			& \scalebox{0.75}[1.0]{MITC4}
			\rule[-2mm]{0mm}{6mm} \\
			&& $\lambda$ [-]  & $\lambda$ [-] && $\lambda$ [-] & $\lambda$ [-] & $\lambda$ [-] \\
			\hline \hline
			$161  $ & $577  $ & $3.3402$ & $3.3809$ & $385  $ & $3.3696$ & $3.3709$ & $3.3867$ \\
			$644  $ & $2128 $ & $3.1558$ & $3.1833$ & $1419 $ & $3.1800$ & $3.1848$ & $3.1845$ \\
			$2576 $ & $8128 $ & $3.1244$ & $3.1454$ & $5419 $ & $3.1371$ & $3.1389$ & $3.1383$ \\
			$10304$ & $31720$ & $3.1129$ & $3.1398$ & $21147$ & $3.1192$ & $3.1188$ & $3.1197$ \\
			\hline
			\multicolumn{8}{c}{Reference solution: $\lambda=3.1067$ (Aitken's extrapolation of the MITC9 results)} \\
			\hline
		\end{tabularx}
		\captionof{table}{Block with asymmetric holes (quad elements).  \label{Tab:B05QUAD}}
	\end{table}
\end{center}

\begin{center}
	\begin{table}
		\centering
		\begin{tabularx}{123mm}{c *{5}{Y}}
			\hline
			\#elements & \#dofs & \scalebox{0.75}[1.0]{MITC9} & \#dofs & \scalebox{0.75}[1.0]{MITC16}  \rule[-2mm]{0mm}{6mm} \\
			&     & $\lambda$ [-]               &        & $\lambda$ [-] \\
			\hline \hline
			$161  $  & $1419     $ & $3.3113$ & $3097      $ & $3.2163$  \\
			$644  $  & $5419     $ & $3.1266$ & $11995     $ & $3.1200$  \\
			$2576 $  & $21147    $ & $3.1137$ & $47179     $ & $3.1085$  \\
			$10304$  & $83515    $ & $3.1092$ & $187099    $ & $3.1083$  \\
			\hline
			\multicolumn{5}{c}{Reference solution: $\lambda=3.1067$ (Aitken's extrapolation of the MITC9 results)} \\
			\hline
		\end{tabularx}
		\captionof{table}{Block with asymmetric holes (high order quad elements).  \label{Tab:B05HQUAD}}
	\end{table}
\end{center}

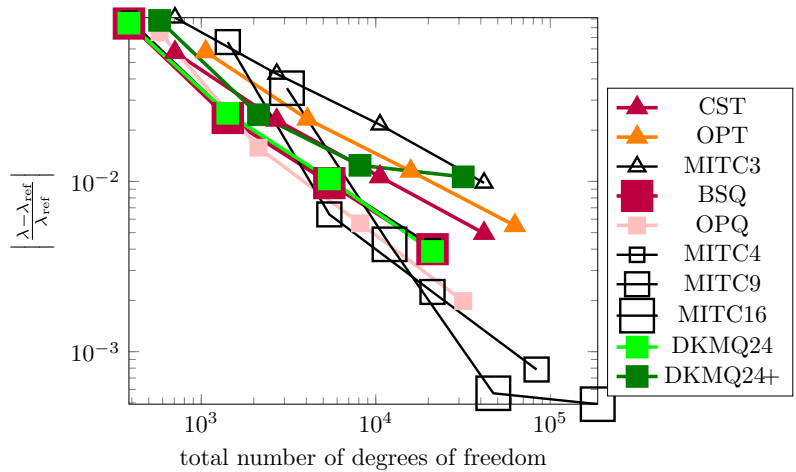
\begin{figure}
	\centering
	\scalebox{.9}
	{
		\begin{tikzpicture}
			\begin{loglogaxis}[
				xlabel=\text{total number of degrees of freedom},
				ylabel= $\left|\frac{\lambda-\lambda_{\text{ref}}}{\lambda_{\text{ref}}}\right|$,
				x tick label style={
					/pgf/number format/1000 sep={}  
				},
				y tick label style={
					/pgf/number format/1000 sep={},
					/pgf/number format/fixed,       
					/pgf/number format/precision=5  
				},
				scaled y ticks=false, 
				xmin=385,
				xmax=187099,
				ymin=0.0004915108954801527,
				ymax=0.09154559415362762,
				legend style={at={(1.02,0.02)},anchor=south west}
				]
				\addplot[line width=1.5pt, color=purple, mark=triangle*, mark size=4pt] coordinates {
					(707, 0.05726665222498009)
					(2707, 0.022888893109094498)
					(10571, 0.010649402735401452)
					(41755, 0.0049521251650038465)
				};
				\addplot[line width=1.5pt, color=orange, mark=triangle*, mark size=4pt] coordinates {
					(1060, 0.057757197480193936)
					(4060, 0.023220429600479247)
					(15856, 0.011488544126702322)
					(62632, 0.005503827637075042)
				};
				\addplot[line width=1pt, color=black, mark=triangle, mark size=4pt] coordinates {
					(707, 0.09154559415362762)
					(2707, 0.042788807714306994)
					(10571, 0.02153796237672395)
					(41755, 0.009809617583923045)
				};
				\addplot[line width=1.5pt, color=purple, mark=square*, mark size=6pt] coordinates {
					(385, 0.0845926623571215)
					(1419, 0.023583832220647593)
					(5419, 0.009771635733453708)
					(21147, 0.003997750701940035)
				};
				\addplot[line width=1.5pt, color=pink, mark=square*, mark size=3pt] coordinates {
					(577, 0.07515352875531871)
					(2128, 0.015775665028404353)
					(8128, 0.00567571160445315)
					(31720, 0.001981171946090278)
				};
				\addplot[line width=1pt, color=black, mark=square, mark size=3pt] coordinates {
					(385, 0.09011515903934252)
					(1419, 0.025012979814577764)
					(5419, 0.01015885748018746)
					(21147, 0.004159978266656389)
				};
				\addplot[line width=1pt, color=black, mark=square, mark size=5pt] coordinates {
					(1419, 0.06583960650802922)
					(5419, 0.006399941804079989)
					(21147, 0.002241894817955923)
					(83515, 0.0007853874164841868)
				};
				\addplot[line width=1pt, color=black, mark=square, mark size=7pt] coordinates {
					(3097, 0.03527709208377765)
					(11995, 0.0042703831370883814)
					(47179, 0.0005697277570397749)
					(187099, 0.0004915108954801527)
				};
				\addplot[line width=1.5pt, color=green, mark=square*, mark size=4pt] coordinates {
					(385, 0.085031384918051)
					(1419, 0.02511244076199322)
					(5419, 0.01033943220996116)
					(21147, 0.0038783331890238242)
				};
				\addplot[line width=1.5pt, color=col1, mark=square*, mark size=4pt] coordinates {
					(577, 0.08824954204510381)
					(2128, 0.02463959891165904)
					(8128, 0.012426502705241252)
					(31720, 0.010636849411941214)
				};
				\legend{CST,OPT,MITC3,BSQ,OPQ,MITC4,MITC9,MITC16,DKMQ24,DKMQ24\!+~}
			\end{loglogaxis}
		\end{tikzpicture}
	}
	\caption{B05: Block with asymmetric holes }
	\label{fig:B05}
\end{figure}

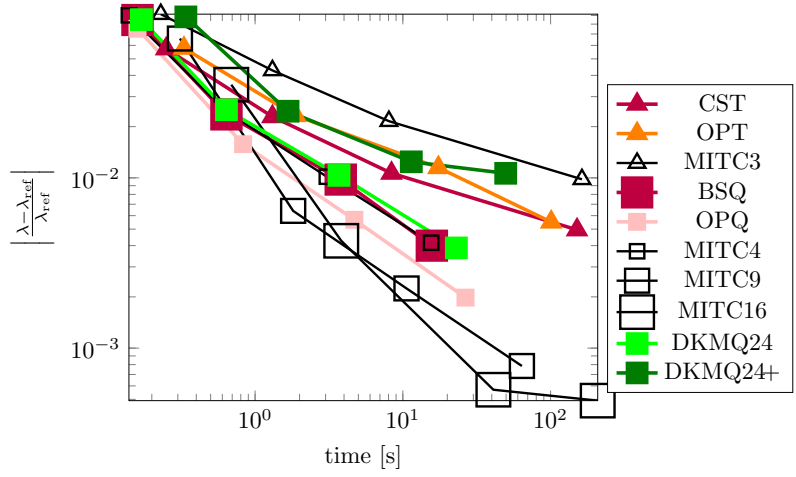
\begin{figure}
	\centering
	\scalebox{.9}
	{
		\begin{tikzpicture}
			\begin{loglogaxis}[
				xlabel=\text{time [s]},
				ylabel= $\left|\frac{\lambda-\lambda_{\text{ref}}}{\lambda_{\text{ref}}}\right|$,
				x tick label style={
					/pgf/number format/1000 sep={}  
				},
				y tick label style={
					/pgf/number format/1000 sep={},
					/pgf/number format/fixed,       
					/pgf/number format/precision=5  
				},
				scaled y ticks=false, 
				xmin=0.14,
				xmax=208.45,
				ymin=0.0004915108954801527,
				ymax=0.09154559415362762,
				legend style={at={(1.02,0.02)},anchor=south west}
				]
				\addplot[line width=1.5pt, color=purple, mark=triangle*, mark size=4pt] coordinates {
					(0.25, 0.05726665222498009)
					(1.31, 0.022888893109094498)
					(8.38, 0.010649402735401452)
					(151.02, 0.0049521251650038465)
				};
				\addplot[line width=1.5pt, color=orange, mark=triangle*, mark size=4pt] coordinates {
					(0.33, 0.057757197480193936)
					(1.88, 0.023220429600479247)
					(17.41, 0.011488544126702322)
					(100.95, 0.005503827637075042)
				};
				\addplot[line width=1pt, color=black, mark=triangle, mark size=4pt] coordinates {
					(0.23, 0.09154559415362762)
					(1.31, 0.042788807714306994)
					(8.05, 0.02153796237672395)
					(163.08, 0.009809617583923045)
				};
				\addplot[line width=1.5pt, color=purple, mark=square*, mark size=6pt] coordinates {
					(0.16, 0.0845926623571215)
					(0.64, 0.023583832220647593)
					(3.8, 0.009771635733453708)
					(15.72, 0.003997750701940035)
				};
				\addplot[line width=1.5pt, color=pink, mark=square*, mark size=3pt] coordinates {
					(0.16, 0.07515352875531871)
					(0.83, 0.015775665028404353)
					(4.7, 0.00567571160445315)
					(26.56, 0.001981171946090278)
				};
				\addplot[line width=1pt, color=black, mark=square, mark size=3pt] coordinates {
					(0.14, 0.09011515903934252)
					(0.62, 0.025012979814577764)
					(3.06, 0.01015885748018746)
					(15.55, 0.004159978266656389)
				};
				\addplot[line width=1pt, color=black, mark=square, mark size=5pt] coordinates {
					(0.31, 0.06583960650802922)
					(1.81, 0.006399941804079989)
					(10.59, 0.002241894817955923)
					(63.97, 0.0007853874164841868)
				};
				\addplot[line width=1pt, color=black, mark=square, mark size=7pt] coordinates {
					(0.69, 0.03527709208377765)
					(3.83, 0.0042703831370883814)
					(41.02, 0.0005697277570397749)
					(208.45, 0.0004915108954801527)
				};
				\addplot[line width=1.5pt, color=green, mark=square*, mark size=4pt] coordinates {
					(0.17, 0.085031384918051)
					(0.64, 0.02511244076199322)
					(3.7, 0.01033943220996116)
					(23.03, 0.0038783331890238242)
				};
				\addplot[line width=1.5pt, color=col1, mark=square*, mark size=4pt] coordinates {
					(0.34, 0.08824954204510381)
					(1.67, 0.02463959891165904)
					(11.41, 0.012426502705241252)
					(49.7, 0.010636849411941214)
				};
				\legend{CST,OPT,MITC3,BSQ,OPQ,MITC4,MITC9,MITC16,DKMQ24,DKMQ24\!+~}
			\end{loglogaxis}
		\end{tikzpicture}
	}
	\caption{B05: Block with asymmetric holes }
	\label{fig:B05TIME}
\end{figure}


\begin{figure}[H]
	\centering
	\begin{overpic}[width=0.5\textwidth]{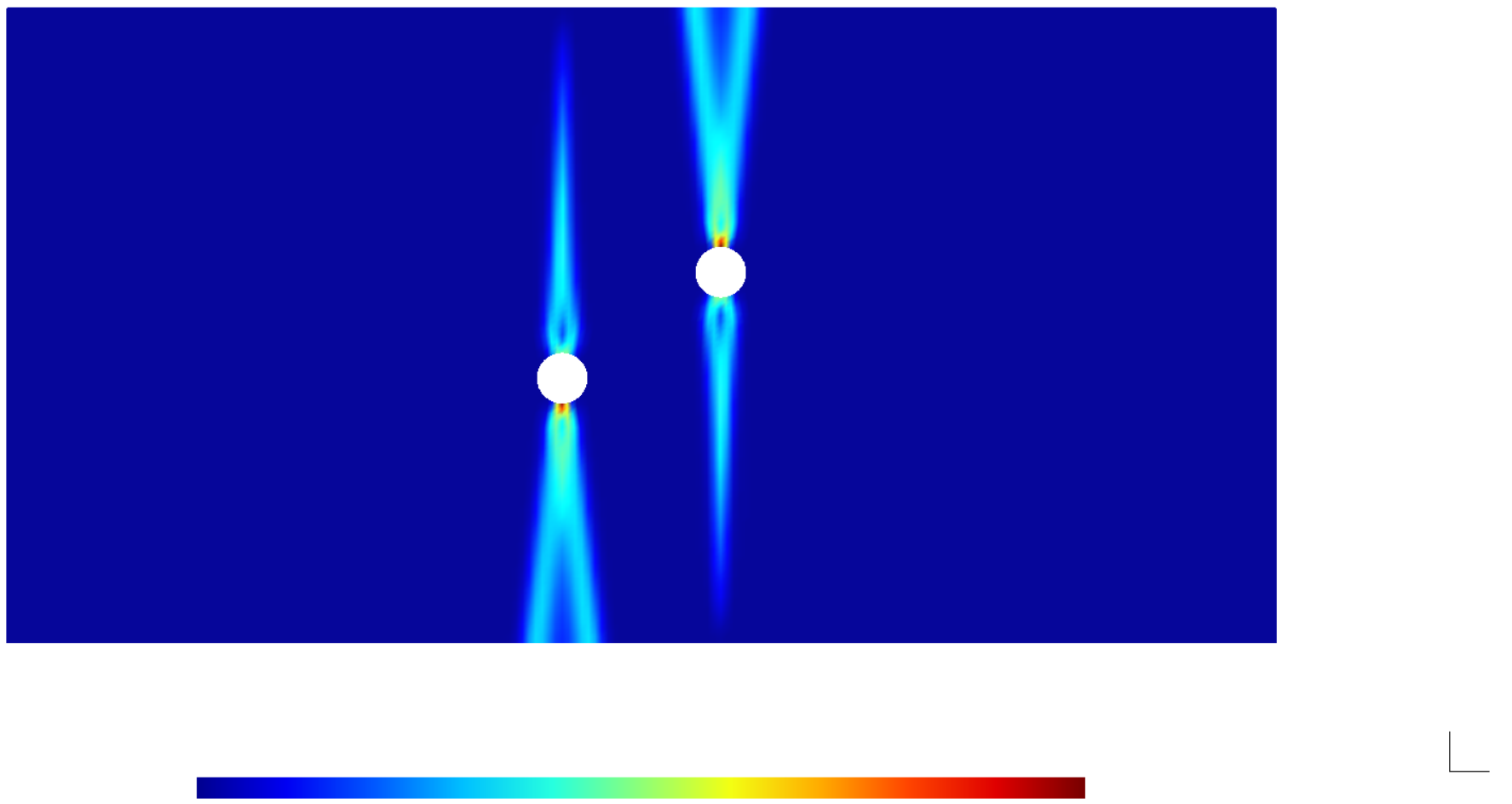} 
		\put(99,-2){$X$}
		\put(92.6,4.5){$Y$}
		\put(13,-4.3){$0$}
		\put(40,-4.3){$13.55$}
		\put(67,-4.3){$27.1$}
		\put(42,4.5){$D_{\p}$}
	\end{overpic}
	\caption{Block with asymmetric holes -- dissipation power at collapse.}
	\label{PlateWithHoles}
\end{figure}	

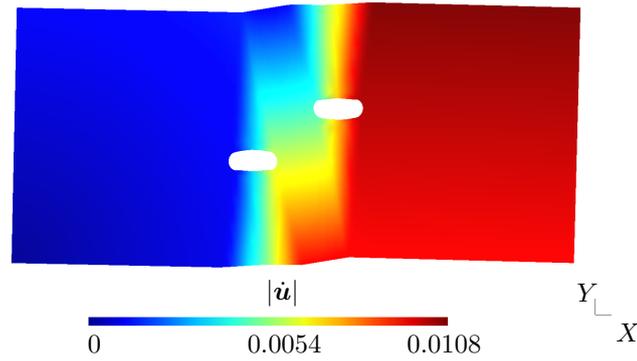
\begin{figure}[H]
	\centering
	\begin{overpic}[width=0.5\textwidth]{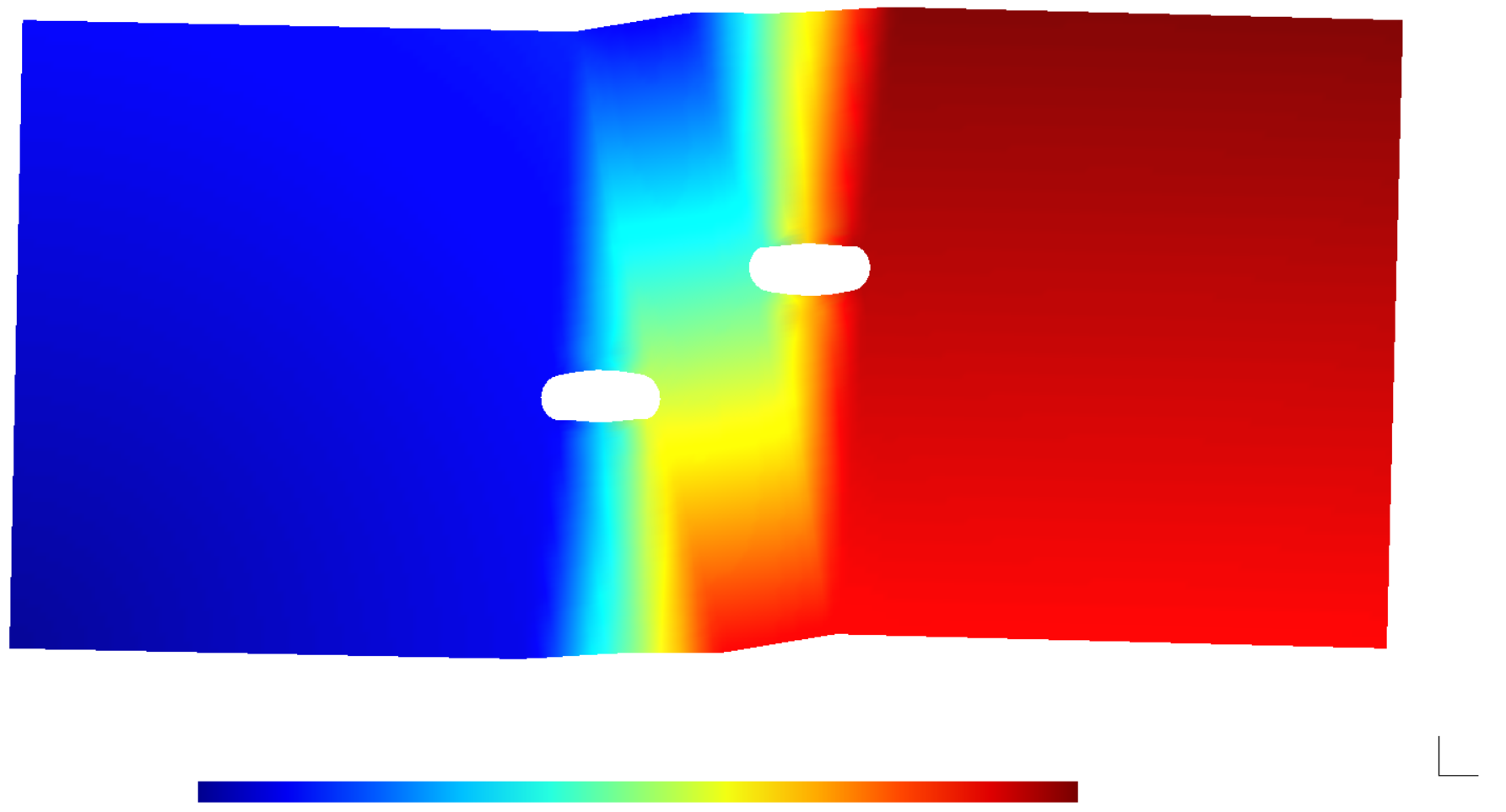} 
		\put(99,-2){$X$}
		\put(92.6,4.5){$Y$}
		\put(13,-3.8){$0$}
		\put(39,-3.8){$0.0054$}
		\put(65,-3.8){$0.0108$}
		\put(42,4.75){$|\dot{\bm{u}}|$}
	\end{overpic}
	\caption{Block with asymmetric holes -- velocity field at collapse.}
	\label{PlateWithHolesDiss}
\end{figure}

\subsection{Simply supported plate in the sequential limit analysis} 
\label{b6}

\begin{figure}[H]
	\centering
	\begin{overpic}[width=0.35\textwidth]{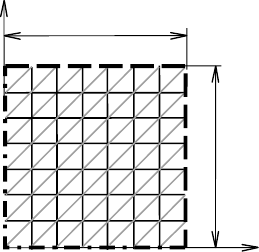} 
		\put(35,87.5){$\frac{L}{2}$}
\put(86,35){$\frac{L}{2}$}
\put(25,74){thickness\,$h$}
\put(95.5,5){$X$}
\put(-6,92.5){$Y$}
	\end{overpic}
	\caption{Simply supported pressure loaded square plate. The triangular mesh is depicted in a grey colour. Parameters: $L=20\,\text{mm},\,h=0.2\,\text{mm}$, yield strength $\sigma_0=200\,\text{MPa}$, loading pressure $p=20\,\text{kPa}$.}
	\label{SimplySupportedPlate}
\end{figure}
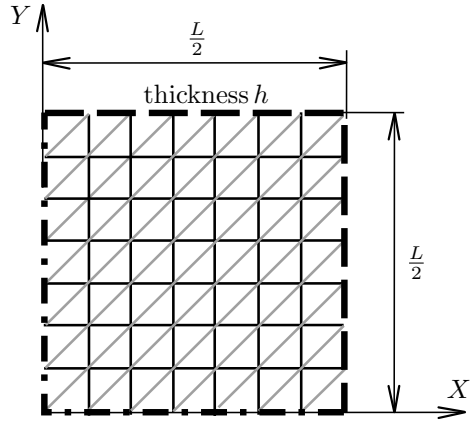

\begin{figure}
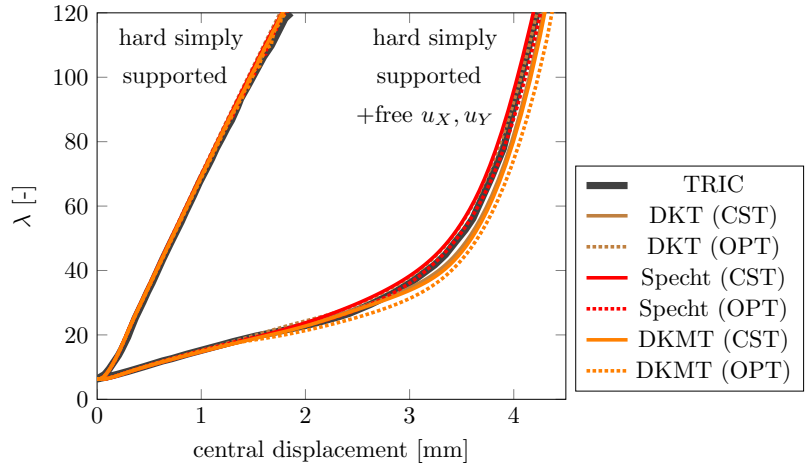

	\centering
	\scalebox{.9}
	{

	}
	\caption{SFELA of a square pressure loaded plate, two types of boundary conditions are tested: 1. hard simply supported, and 2. hard simply supported with free $u_X$ and $u_Y$. The step is $0.0325$~mm. Comparison of triangle elements.}
	\label{fig:B06T}
\end{figure}

\begin{figure}
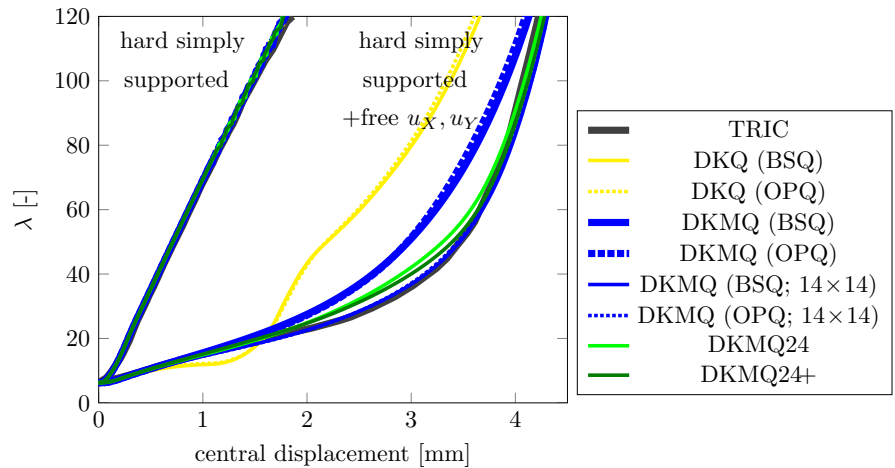

	\centering
	\scalebox{.9}
	{
		%
	}
	\caption{SFELA of a square pressure loaded plate, two types of boundary conditions are tested: 1. hard simply supported, and 2. hard simply supported with free $u_X$ and $u_Y$. The step is $0.0325$~mm. Comparison of the TRIC and quadrangle elements. The DKMQ (BSQ) / DKMQ (OPQ) elements require a $14\!\times\!14$ mesh to achieve comparable accuracy to other shell elements under boundary conditions 2. This is due to the reduced accuracy of the rigid links correction (warping matrix). The DKQ (BSQ) / DKQ (OPQ) elements exhibit higher error compared to other elements under the same boundary conditions. This error is not reduced on finer meshes, even when we tried an $84\!\times\!84$ mesh.}
	\label{fig:B06Q}
\end{figure}

A simply supported, pressure-loaded square plate is considered in the SFELA (see Fig. \ref{SimplySupportedPlate}). Only one quarter is computed due to its symmetry. Symmetry conditions are as follows:
$u_X=\varphi_Y=\varphi_Z=0$ at $X=0$ and $u_Y=\varphi_X=\varphi_Z=0$ at $Y=0$. Two different boundary conditions are used: 
\begin{enumerate}
	\item $u_X=u_Y=u_Z=\varphi_X=\varphi_Z=0$ at $X=\frac{L}{2}$ and $u_X=u_Y=u_Z=\varphi_Y=\varphi_Z=0$ at $Y=\frac{L}{2}$ (hard simply supported conditions)
	\item $u_Z=\varphi_X=\varphi_Z=0$ at $X=\frac{L}{2}$ and $u_Z=\varphi_Y=\varphi_Z=0$ at $Y=\frac{L}{2}$ (hard simply supported conditions together with free $u_X=u_Y=0$, which represents in-plane sliding of supports)
\end{enumerate} 
Triangle elements are calculated with a $7\!\times\!7\!\times\!2$ mesh, while quad elements are calculated with a $7\!\times\!7$ mesh (see Fig. \ref{SimplySupportedPlate}). The Ilyushin stress-resultant plasticity model is considered. Since the analytical solution for this example is not available, our simulations are compared with results by Corradi (\cite{Corradi}, Fig. 5b), obtained with the help of the TRIC shell element (\cite{TRIC}). Numerical results are depicted in Figures \ref{fig:B06T} and \ref{fig:B06Q}. Under boundary conditions 2, the DKMQ (BSQ) / DKMQ (OPQ) elements require a 14×14 mesh to achieve accuracy comparable to other shell elements. This is attributed to the reduced accuracy of the rigid links correction (warping matrix). Under the same boundary conditions, the DKQ (BSQ) / DKQ (OPQ) elements exhibit higher error compared to other elements. This error is not reduced on finer meshes, even when the $84\!\times\!84$ mesh was tried. At this benchmark, none of the MITC elements converge. For both variants of the boundary conditions, all other elements provide practically identical results in comparison with Corradi’s TRIC element. 

\FloatBarrier
\FloatBarrier
\section{Conclusions}


 


The DK, DKM, and MITC classes of continuous shell elements underwent testing in a variety of FELA and SFELA benchmark problems. These benchmark problems incorporated different geometric models (including plates, membranes in the plane stress or plane strain approximation, and general shells) and different isotropic and orthotropic material models. The performance of these elements was then compared with other elements from literature that are used in the context of FELA and SFELA. The main results are: 

\begin{itemize}	
	\item The DK and DKM classes of continuous shell elements perform well in the context of the FELA. These elements do not offer a rigid upper bound, on the other hand offer a significantly higher convergence rate when compared to the triangle elements designed for the strict upper bound (\cite{Bleyer2013}, \cite{Makrodimopoulos}, \cite{Ciria}). Among all the tested shell elements, the quadrangles DKMQ24 and DKMQ24+ appeared to be numerically most efficient for general shell problems.
	\item The higher numerical efficiency of quadrangle elements over triangle elements is observed.
	\item The MITC shell elements work in the context of FELA only in the pure membrane case.
	\item In the SFELA benchmark all the DK and DKM shell elements except for the DKQ quadrangle element coincide with the reference results given by the TRIC triangle. To replicate the reference results calculated with the TRIC triangle, the DKMQ quadrangle requires a finer mesh compared to other tested elements. 
	\item However, for the DKM element class and the clamped plate problem, we encountered a loss of convergence on finer meshes if the ratio of the shell thickness to the shell characteristic length exceeds a certain limit. In this case, a numerical penalization was proposed to recover convergence.
	\item In all our test problems, convergence is monotonic: from below in one problem (the spherical cap problem) and from above in all other problems.
\end{itemize}
We can conclude that continuous shell elements of the DK and DKM class perform well not only in the context of FEM but also in FELA and SFELA contexts. They have lower discretization error compared to triangle elements, which guarantee a strict upper bound.

 
\section*{Appendix -- Derivation of the dissipation power density for the considered stress-based failure criteria}
\subsection*{Plane-Stress Tresca Failure Criterion}
The Tresca in-plane plane stress failure criterion is given by:
\begin{align}
	\max(|\sigma_2-\sigma_1|,|\sigma_1|,|\sigma_2|)\le \sigma_0,
\end{align}
where $\sigma_1, \sigma_2$ are eigenstresses. The dissipation power density is defined by:
\begin{align}
	d_{\p}^{z}(\dot{\bm{\epsilon}})&=\max_{\forall \bm{\sigma}:~\!\max(|\sigma_2-\sigma_1|,|\sigma_1|,|\sigma_2|)\le \sigma_0}~\sigma_1\dot{\epsilon}_1+\sigma_2\dot{\epsilon}_2. \label{eq546}
\end{align}
The expression $\sigma_1\dot{\epsilon}_1+\sigma_2\dot{\epsilon}_2$ is linear (for given strain rates $\dot{\epsilon}_1,\dot{\epsilon}_2$) and therefore the maximum can be found only at boundary points: $[0,\pm \sigma_0],~[\pm \sigma_0, 0],~\pm[\sigma_0, \sigma_0]$. This immediately yields:
\begin{align}
	d_{\p}^{z}(\dot{\bm{\epsilon}})&=\sigma_0\max\left(|\dot{\epsilon}_1|, |\dot{\epsilon}_2|, |\dot{\epsilon}_1+\dot{\epsilon}_2|\right). 
\end{align}
By using standard results for eigenvectors in 2D, the final expression of the dissipation power density is obtained:
\begin{align}
	d_{\p}^z(\dot{\bm{\epsilon}})&=\sigma_0\max\left[\left|\frac{\dot{\epsilon}_x+\dot{\epsilon}_y-\sqrt{A}}{2}\right|, \left|\frac{\dot{\epsilon}_x+\dot{\epsilon}_y+\sqrt{A}}{2}\right|,|\dot{\epsilon}_x+\dot{\epsilon}_y|\right],  \label{eq547} \\
	A&=(\dot{\epsilon}_x - \dot{\epsilon}_y)^2+\dot{\curlyvee}_{xy}^2. \label{eq547b}
\end{align}
The expression (\ref{eq547}) can be easily implemented in the SOCP framework, where $A$ is calculated as a cone and each absolute value is prescribed by two linear inequalities (without an explicit use of absolute values).
\subsection*{Plane-strain Tsai-Wu failure criterion}
The plane-strain Tsai-Wu failure criterion is given by the yield function:
\begin{align}
	&f=\sigma_x\left[\frac{1}{f_{\t x}}-\frac{1}{f_{\c x}}\right]+
	(\sigma_y+\sigma_z)\left[\frac{1}{f_{\t y}}-\frac{1}{f_{\c y}}\right]+ \frac{\sigma_x^2}{f_{\t x}f_{\c x}}+\frac{\sigma_y^2+\sigma_z^2}{f_{\t y}f_{\c y}} +\frac{\tau_{xy}^2}{f_{\vv xy}^2} \le 1. 
	\label{functionf2}
\end{align}
The dissipation power density is defined by:
\begin{align}
	d_{\p}^z(\dot{\bm{\epsilon}})=&\max_{f(\bm{\sigma})\le 1}~\sigma_x\dot{\epsilon}_x+\sigma_y\dot{\epsilon}_y+\tau_{xy}\dot{\curlyvee}_{xy}.
\end{align}
The Lagrangian and its derivatives yield:
\begin{align}
	\mathcal{L}&=\sigma_x\dot{\epsilon}_x+\sigma_y\dot{\epsilon}_y+\tau_{xy}\dot{\curlyvee}_{xy}+\lambda f,\\
	\PART{\mathcal{L}}{\tau_{xz}}&=\frac{2\tau_{xz}}{f_{\vv xy}^2}\lambda=0~\implies~\tau_{xz}=0,~~~~
	\PART{\mathcal{L}}{\tau_{yz}}=\frac{2\tau_{yz}}{f_{\vv xy}^2}\lambda=0~\implies~\tau_{yz}=0,\\
	\PART{\mathcal{L}}{\tau_{xy}}&=\dot{\curlyvee}_{xy}+\frac{2\tau_{xy}}{f_{\vv xy}^2}\lambda=0~\implies~\tau_{xy}=-\frac{\dot{\curlyvee}_{xy}f_{\vv xy}^2}{2\lambda},\\
	\PART{\mathcal{L}}{\sigma_x}&=\dot{\epsilon}_x+\lambda\left(\frac{1}{f_{\t x}}-\frac{1}{f_{\c x}}+\frac{2\sigma_x}{f_{\t x}f_{\c x}}\right)=0~\implies~\sigma_x=\frac{f_{\t x}-f_{\c x}}{2}-\frac{f_{\t x}f_{\c x}}{2\lambda}\dot{\epsilon}_x,\\
	\PART{\mathcal{L}}{\sigma_y}&=\dot{\epsilon}_y+\lambda\left(\frac{1}{f_{\t y}}-\frac{1}{f_{\c y}}+\frac{2\sigma_y}{f_{\t y}f_{\c y}}\right)=0~\implies~\sigma_y=\frac{f_{\t y}-f_{\c y}}{2}-\frac{f_{\t y}f_{\ cy}}{2\lambda}\dot{\epsilon}_y,\\
	\PART{\mathcal{L}}{\sigma_z}&=                    \lambda\left(\frac{1}{f_{\t y}}-\frac{1}{f_{\c y}}+\frac{2\sigma_z}{f_{\t y}f_{\c y}}\right)=0~\implies~\sigma_z=\frac{f_{\ ty}-f_{\c y}}{2}. 
\end{align}
By inserting all stresses into a function $f$ (Eq.(\ref{functionf2})), we get the expression for $\lambda$:
\begin{align}
	2\lambda=\pm\sqrt{
		\frac{f_{\t x}f_{\c x}(\dot{\epsilon}_x)^2+
			f_{\t y}f_{\c y}(\dot{\epsilon}_y)^2+
			(f_{\vv xy}\dot{\curlyvee}_{xy})^2}{\chi}},~\chi=1+\frac{1}{4}\frac{(f_{\t x}-f_{\c x})^2}{f_{\t x}f_{\c x}}+\frac{1}{2}\frac{(f_{\t y}-f_{\c y})^2}{f_{\t y}f_{\c y}},
\end{align}
where the minus sign is chosen according to the positiveness of $d_{\text{p}}^z(\dot{\bm{\epsilon}})$. The final expression of the dissipation power density has the form: 
\begin{align}
	d_{\p}^z(\dot{\bm{\epsilon}})=&\frac{f_{\t x}-f_{\c x}}{2}\dot{\epsilon}_x+\frac{f_{\t y}-f_{\c y}}{2}\dot{\epsilon}_y + \sqrt{(\bm{R}^{-\T}\dot{\bm{\epsilon}})^{\T}\bm{R}^{-\T}\dot{\bm{\epsilon}}},
\end{align}
where
\begin{align}
	\bm{R}^{-\T}=&\left[
	\begin{array}{ccc}
		\sqrt{\chi f_{\t x}f_{\c x}} & 0 & 0\\
		0 & \sqrt{\chi f_{\t y}f_{\c y}} & 0\\
		0 & 0 &\sqrt{\chi}f_{\vv xy}\\
	\end{array}\right],~\chi=1+\frac{1}{4}\frac{(f_{\t x}-f_{\c x})^2}{f_{\t x}f_{\c x}}+\frac{1}{2}\frac{(f_{\t y}-f_{\c y})^2}{f_{\t y}f_{\c y}}.
\end{align}
 


\end{document}